\documentclass[12pt]{article}
\usepackage{amsmath,amssymb,enumerate, color}
\usepackage{graphicx,xcolor}
\usepackage[all]{xy}
\usepackage{bm}
\newtheorem{tm}{Theorem}[section]
\newtheorem{lm}[tm]{Lemma}

\newtheorem{pr}[tm]{Proposition}

\makeatletter
\newcommand{\subscripts}[3]{%
  \@mathmeasure\z@\displaystyle{#2}%
  \global\setbox\@ne\vbox to\ht\z@{}\dp\@ne\dp\z@
  \setbox\tw@\box\@ne
  \@mathmeasure4\displaystyle{\copy\tw@_{#1}}%
  \@mathmeasure6\displaystyle{{#2}_{#3}}%
  \dimen@-\wd6 \advance\dimen@\wd4 \advance\dimen@\wd\z@
  \hbox to\dimen@{}\mathop{\kern-\dimen@\box4\box6}%
}
\makeatother
\newcommand{\qed}{~~\hbox{\rule{4pt}{8pt}}}

\newcommand{\A}{\mathcal{A}}
\newcommand{\T}{\mathcal{T}}
\newcommand{\g}{\frak{g}}

\newcommand{\bg}{\bm{\mathrm{g}}}
\newcommand{\nn}{\nonumber}
\newcommand{\Prob}{\bm{\mathrm{P}}}

\newcommand{\h}{\mathrm{H}}
\newcommand{\e}{\mathrm{e}}
\newcommand{\x}{\bm{\mathrm{x}}}
\newcommand{\y}{\bm{\mathrm{y}}}
\newcommand{\ve}{\varepsilon}
\newcommand{\dis}{\displaystyle}

\newcommand{\Ker}{\mathrm{Ker}\,}
\newcommand{\im}{\mathrm{Im}\,}
\newcommand{\del}{\partial}
\newcommand{\ol}{\overline}
\newcommand{\la}{\langle}
\newcommand{\La}{\langle \! \langle}
\newcommand{\ra}{\rangle}
\newcommand{\Ra}{\rangle\!\rangle}
\newcommand{\LA}{\to}
\newcommand{\Span}{\mathrm{span}_{\mathbb{R}}}
\newcommand{\G}{G_{\Gamma}}

\newcommand{\Hom}{\mathrm{Hom}}
\newcommand{\Dom}{\mathrm{Dom}}

\newcommand{\Log}{\text{\rm{\,log\,}}}
\newcommand{\Exp}{\text{\rm{exp\,}}}
\newcommand{\hol}{\text{-H\"ol}}

\allowdisplaybreaks[4]

\makeatletter
 
 \@addtoreset{equation}{section}
\makeatother
\begin{document}
\setlength{\baselineskip}
{15.5pt}
\title{
Central limit theorems for non-symmetric random walks 
on nilpotent covering graphs: Part II
}
\author{\Large
{Satoshi Ishiwata\footnote{
Department of Mathematical Sciences, Faculty of Science, Yamagata University,
1-4-12, Kojirakawa, Yamagata 990-8560, Japan
(e-mail: {\tt ishiwata@sci.kj.yamagata-u.ac.jp})}
\footnote{Partially supported by Grant-in-Aid for Young Scientists (B)(No. 25800034) 
and Grant-in-Aid for Scientific Research (C)(No. 17K05215), JSPS}, 
Hiroshi Kawabi\footnote{
Department of Mathematics, Faculty of Economics, Keio University, 4-1-1, 
Hiyoshi, Kohoku-ku, Yokohama 223-8521, Japan (e-mail: {\tt kawabi@keio.jp})}
\footnote{Partially supported by Grant-in-Aid for Scientific Research 
(C)(No. 26400134, No. 17K05300), JSPS}
}
{\Large{and Ryuya Namba\footnote{
Graduate School of Natural Sciences, Okayama University, 3-1-1, 
Tsushima-Naka, Kita-ku, Okayama 700-8530, Japan (e-mail: {\tt{rnamba@s.okayama-u.ac.jp}})}
\hspace{0.1mm} \footnote{Partially supported by Grant-in-Aid JSPS Research Fellow 
(No. 18J10225), JSPS}
}}
}
\date{\today}
%
\maketitle 
%
%
\begin{abstract}

In the present paper, as a continuation of our preceding paper \cite{IKN}, we study another kind of 
central limit theorems (CLTs) for non-symmetric random walks on nilpotent covering
graphs 
from a view point of discrete geometric analysis developed by Kotani and Sunada.
We introduce a one-parameter 
family of random walks which interpolates between the original
non-symmetric random walk and the symmetrized one.
We first prove a semigroup CLT for the family of random walks by realizing 
the nilpotent covering graph into a nilpotent Lie group via discrete harmonic maps.  
The limiting diffusion semigroup is generated by
the homogenized sub-Laplacian with a constant drift of the asymptotic direction on the 
nilpotent Lie group, which is equipped with the Albanese metric associated with the 
symmetrized random walk. 
We next prove a functional CLT (i.e., Donsker-type
invariance principle) in a H\"older space over the nilpotent Lie group
by combining
the semigroup CLT, standard martingale techniques, and a novel pathwise 
argument inspired by rough path theory.
Applying the corrector method, 
we finally extend these CLTs to the case where the realizations are not necessarily 
harmonic. 

\vspace{3mm}

\vspace{3mm}
\noindent
{\bf AMS Classification (2010):} 
60F17, 
60G50, 
60J10, 
22E25
\end{abstract}


\section{Introduction}

The long time behavior of random walks on graphs is one of principal themes 
among probability theory, harmonic analysis, geometry, graph theory, group theory and so on.
In particular, the central limit theorem (CLT), a generalization of the Laplace--de Moivre theorem, 
has been studied intensively and extensively in various settings. 
Our main concern is such a long time behavior of random walks with periodic structures by 
using ideas from homogenization theory. Generally speaking, homogenization theory is 
a method which relates a periodic system to the corresponding homogenized system 
through a scaling relation between
the time and the underlying space 
(cf.~Bensoussan--Lions--Papanicolaou \cite{BLP}). 
However, since the notion of the scale change on graphs
is not defined, it is not possible to apply this method directly to the case where the underlying space is 
a graph. Therefore, it is necessary to find a realization of the graph, preserving
the periodicity, in a space on which a scaling is defined. 

In a series of papers
\cite{Kotani, Kotani contemp,
KS00-CMP, KS00-TAMS, KS06, KSS},
Kotani and Sunada studied 
long time asymptotics of symmetric random walks on a {\it{crystal lattice}} $X$, a covering graph of a
finite graph $X_{0}$ whose covering transformation group $\Gamma$ is abelian, by 
placing a special emphasis on the geometric feature.
We remark that the crystal lattice has
inhomogeneous local structures though it has a periodic global structure.
Especially, in \cite{KS00-TAMS}, 
they introduced the notion of {\it{standard realization}}, which is a discrete
harmonic map $\Phi_{0}$ from a crystal lattice $X$ into the Euclidean space 
$\Gamma \otimes {\mathbb R}$
equipped with the {\it{Albanese metric}},
to characterize an equilibrium configuration of crystals. 
In \cite{KS00-CMP}, they proved the CLT by applying 
the homogenization method mentioned above 
through the standard realization $\Phi_0$.
As the scaling limit, they captured 
a homogenized 
Laplacian
on $\Gamma \otimes \mathbb R$. 
In terms of probability theory, it means that, for fixed $0\leq t \leq 1$,
a sequence of $\Gamma \otimes \mathbb{R}$-valued random variables 
$\big\{ n^{-1/2} \Phi_0(w_{[nt]}) \big\}_{n=1}^\infty$ 
converges to $B_t$  as $n \to \infty$ in law. 
Here $\{w_n\}_{n=0}^\infty$ is the given random walk on $X$ and 
$(B_t)_{0\leq t \leq 1}$ is a standard Brownian motion on $\Gamma \otimes \mathbb{R}$
equipped with the Albanese metric.
In their proof, both the symmetry of the given random walk $\{w_n\}_{n=0}^\infty$
and the harmonicity of the realization
$\Phi_{0}$ play an important role to show the convergence of the 
sequence of infinitesimal generators associated with 
$\big\{ n^{-1/2} \Phi_0(w_{[nt]}) ;~0\leq t\leq 1 \big\}_{n=1}^\infty$. 
Indeed, these properties are effectively used to delete a 
diverging drift term as $n\to \infty$ from the homogenized Laplacian.
%
See also Kotani \cite{Kotani}. 
Ishiwata \cite{Ishiwata} discussed a similar problem to \cite{Kotani, KS00-CMP} for symmetric
random walks on a $\Gamma$-{\it{nilpotent covering graph}} $X$, 
a covering graph of a finite graph $X_{0}$ 
whose covering transformation group $\Gamma$ is 
a finitely generated 
nilpotent group of step $r$ ($r\in \mathbb N$).
See \cite[Section 6]{IKN} for several examples of $\Gamma$-nilpotent 
covering graphs in the case where $\Gamma$ is the 3-dimensional 
discrete Heisenberg group $\mathbb{H}^3(\mathbb{Z})$. 
Needless to say, $X$ is a crystal lattice in the case $r=1$.
In the nilpotent case, $X$ is properly realized into a nilpotent Lie group $G$ such that
$\Gamma$ is isomorphic to a cocompact lattice of $G$.
Thanks to a basic property of the canonical dilation on $G$, the diverging drift 
term appears only in $\g^{(1)}$-direction,
where $\g^{(1)}$ is the generating part of the Lie algebra $\g$ of $G$. 
Hence, it is sufficient to introduce the notion of harmonicity of the realization 
$\Phi_{0}: X \LA G$ only on $\g^{(1)}$ for proving the CLT in the nilpotent case.

If we consider a non-symmetric case, the above method
does not work well because the diverging drift term from the non-symmetry of the given random walk
does not vanish. To overcome this difficulty, in the case of crystal lattices, Ishiwata, Kawabi and Kotani
\cite{IKK} introduced two kinds of schemes.  
One is to replace the usual transition operator by the {\it transition-shift operator},
which ``deletes'' the diverging drift term.
Combining this scheme with
a modification of the harmonicity of the realization $\Phi_{0}$,
they proved that a sequence 
$\big\{ n^{-1/2} \big( \Phi_0(w_{[nt]})-[nt] \rho_{\mathbb R}(\gamma_{p}) \big); 0\leq t \leq 1 \big\}_{n=1}^\infty$
converges in law to 
a  $\Gamma \otimes \mathbb R$-valued standard Brownian motion 
$(B_{t})_{t\geq 0}$ 
as $n \to \infty$.
Here $\rho_{\mathbb{R}}(\gamma_p)\in \Gamma \otimes \mathbb R$ is the so-called 
{\it asymptotic direction} which appears in the law of large numbers for the random walk
$\{\Phi_0(w_n)\}_{n=0}^\infty$
on $\Gamma \otimes \mathbb R$. See \cite{IKK, KS06} for details. 
The other is to introduce a one-parameter family of 
$\Gamma \otimes \mathbb R$-valued 
random walks 
$(\xi^{(\varepsilon)})_{0\leq \varepsilon \leq 1}$ 
which ``weakens'' the diverging drift term, where
this family interpolates 
the original non-symmetric random walk 
$\xi^{(1)}_{n}:=\Phi_{0}(w_{n})$ ($n=0,1,2,\ldots$)
and the symmetrized one $\xi^{(0)}$.
Putting $\varepsilon=n^{-1/2}$ and letting $n\to \infty$, we capture
a drifted Brownian motion $(B_{t}+\rho_{\mathbb{R}}(\gamma_p)t)_{0\leq t\leq 1}$ 
as the limit of a sequence
$\big\{ n^{-1/2} \xi_{[nt]}^{(n^{-1/2})}; 0\leq t \leq 1 \big\}_{n=1}^\infty$.
See Trotter \cite{Trotter} for related early works. We note that this scheme 
is well-known in the study of the hydrodynamic limit of weakly asymmetric exclusion processes. 
See e.g., Kipnis--Landim \cite{KL}, Tanaka \cite{Tanaka} and references therein. 

In our preceding paper \cite{IKN}, we proved a CLT for 
a non-symmetric random walk 
$\{w_n\}_{n=0}^\infty$
on the $\Gamma$-nilpotent covering graph $X$ by applying
the former transition-shift scheme with 
a notion of {\it modified harmonic realization} $\Phi_0: X \to G$, which is a generalization of
\cite{KS06, Ishiwata}. More precisely,  as the CLT-scaling limit mentioned above, we captured 
a diffusion process on $G$ generated by a homogenized sub-Laplacian 
with a non-trivial $\g^{(2)}$-valued drift $\beta(\Phi_0)$ arising from the non-symmetry 
of the given random walk, where $\g^{(2)}:=[\g^{(1)}, \g^{(1)}]$. 

The main purpose of the present paper is to prove 
another kind of CLTs 
for the non-symmetric random walk $\{w_n\}_{n=0}^\infty$ on 
the $\Gamma$-nilpotent covering graph $X$ by 
applying the latter scheme to the nilpotent setting. 
We first introduce a one-parameter family of transition probabilities 
$(p_\ve)_{0 \leq \ve \leq 1}$ on $X$ as the linear interpolation between 
the given transition probability $p_{1}:=p$ and the symmetrized one $p_{0}$, that is,
$p_{\varepsilon}:=p_{0}+\varepsilon (p-p_{0})$ ($0\leq \varepsilon \leq 1$).
For each $\varepsilon$, we next take a 
$\Gamma$-periodic realization 
$\Phi^{(\ve)} : X \LA G$ 
associated with the transition probability $p_\ve$, and define 
a $G$-valued random walk $\{ {\overline{\xi}}^{(\varepsilon)}_{n}\}_{n=0}^{\infty}$ by
${\overline{\xi}}^{(\varepsilon)}_{n}:=\Phi^{(\varepsilon)}(w_{n})$ ($n=0,1,2,\ldots$).
Note that $\Phi^{(\ve)}$ is not necessarily modified harmonic as above.
Then by putting $\varepsilon=n^{-1/2}$
and letting $n\to \infty$, 
we obtain a CLT
for the family of $G$-valued random walks $\{ {\overline{\xi}}^{(n^{-1/2})} \}_{n=1}^{\infty}$. 
More precisely, under suitable assumptions 
{\bf (A1)} and {\bf (A2)}
on the family of $\Gamma$-periodic realizations $(\Phi^{(\ve)})_{0\leq \varepsilon \leq 1}$, 
we prove that 
a sequence $\big\{ \tau_{n^{-1/2}} \big( {\overline{\xi}}^{(n^{-1/2})}_{[nt]} \big);
0\leq t \leq 1 \big\}_{n=1}^{\infty}$ converges 
in law 
to a $G$-valued diffusion process $(Y_{t})_{0\leq t \leq 1}$ as $n \to \infty$, where
$\tau_{n^{-1/2}}:G \to G$ is the canonical dilation 
whose scale is $n^{-1/2}$, and the diffusion process 
$(Y_{t})_{0\leq t \leq 1}$ is generated by 
the homogenized 
sub-Laplacian
with the $\mathfrak g^{(1)}$-valued drift $\rho_{\mathbb{R}}(\gamma_p)$
defined on $G$ equipped with the Albanese metric
$g_{0}^{(0)}$
associated with the symmetrized transition probability $p_{0}$.

Here we would like to emphasize that, to our best knowledge, the present paper provides 
the first result on CLTs in the nilpotent
setting in which a $\mathfrak g^{(1)}$-valued drift appears in the infinitesimal
generator of the limiting diffusion. On the other hand, as mentioned in \cite{IKN},
there are many papers on CLTs in which $\mathfrak g^{(2)}$-valued drift like
$\beta(\Phi_0)$
appears 
in the infinitesimal generator of the limiting diffusion. 
In view of these circumstances, the study of the 
long time asymptotics of random walks on 
more general graphs by applying our latter scheme would be an interesting
problem. 

The rest of the present paper is organized as follows:
We introduce our framework and state the main results 
(Theorems \ref{CLT2} and \ref{FCLT2-general}) in Section 2.
We make a preparation from nilpotent Lie groups and 
discrete geometric analysis in Section 3. 
We devote ourselves to prove the main results in Section 4. 
Since the realization map
$\Phi^{(\varepsilon)}$ is not necessarily modified harmonic, 
several technical difficulties arise in the proof. 
To overcome them, we take a modified harmonic realization $\Phi^{(\varepsilon)}_{0}: X \to G$
and show that the
(${\mathfrak g}^{(1)}$-){\it{corrector}},
the difference between 
$\Phi^{(\varepsilon)}$ and $\Phi^{(\varepsilon)}_{0}$ in 
the ${\mathfrak g}^{(1)}$-direction, is not so big.
This approach is the so-called
{\it{corrector method}} in the context of stochastic homogenization theory,
and it is effectively used in the study of random walks in random environments
(see e.g., Biskup \cite{Biskup}, Papanicolaou--Varadhan \cite{PV}, Kozlov \cite{Kozlov} 
and Kumagai \cite{Kumagai}).
In Section 4.1, 
we present a key property (Proposition \ref{key-beta}) of the family of modified harmonic
realizations $(\Phi_0^{(\ve)})_{0 \leq \ve \leq 1}$.
Combining
this property with Trotter's approximation theorem, we prove a semigroup CLT (Theorem \ref{CLT2}) under
only {\bf (A1)}
in Section 4.2.
In Section 4.3, we prove a functional CLT (Theorem \ref{FCLT2})
for a random walk $\big \{ \xi_{n}^{(n^{-1/2})}:=\Phi_{0}^{(n^{-1/2})}(w_{n}) \big \}_{n=1}^{\infty}$
under not only {\bf (A1)} but {\bf (A2)}.   
To prove such kind of CLT which is much stronger than the semigroup CLT, 
we show 
in Lemma \ref{tightness2}
tightness of the family of probability measures 
induced by the $G$-valued stochastic processes 
given by the geodesic interpolation of 
$\big\{ \tau_{n^{-1/2}} \big( \xi^{(n^{-1/2})}_{[nt]} \big); 0\leq t \leq 1
\big\}_{n=1}^\infty$.
In the case $r=2$, we can easily prove it
by combining the modified harmonicity of each $\Phi_0^{(\ve)}$ with standard martingale techniques.
On the other hand, the same argument as in the case $r=2$ is insufficient 
in the case $r \geq 3$. To handle the higher-step terms,
we employ a pathwise argument similar to \cite{IKN}, which is inspired by
the Lyons extension theorem (cf.~Lyons \cite{Lyons}, Lyons--Qian \cite{LQ}, 
Lyons--Caruana--L\'evy \cite{LCL}
and Friz--Victoir \cite{FV}) in rough path theory. 
Since rough path theory is built on free nilpotent Lie groups and our
nilpotent Lie group $G$ is not necessarily free, we need a careful examination of the proof of
the Lyons extension theorem. Combining Theorem \ref{FCLT2} with
several nice properties of the (${\mathfrak g}^{(1)}$-)corrector,
we then prove that a functional CLT (Theorem \ref{FCLT2-general}) also holds for the family of
random walks $\{ {\overline{\xi}}^{(n^{-1/2})} \}_{n=1}^{\infty}$.
Finally, in Section 5, we give an example of non-symmetric random walks on nilpotent covering graphs 
with explicit computations. 

Throughout the present paper, 
$C$ denotes a positive constant that may change from line to line 
and $O(\cdot)$ stands for the Landau symbol. 
If the dependence of $C$ and $O(\cdot)$ are significant, 
we denote them like $C(N)$ and $O_N(\cdot)$, respectively.


\section{Framework and Results}

Let $\Gamma$ be a torsion free, finitely generated nilpotent group of step $r$ ($r \in \mathbb{N}$) 
and $X=(V, E)$ a $\Gamma$-nilpotent covering graph, where $V$ is 
the set of its vertices and $E$ is the set of all oriented edges. 
For $e\in E$, we denote by $o(e)$, $t(e)$ and $\overline e$ the origin, 
the terminus and the inverse edge of $e$, respectively. 
We set $E_{x}=\{ e\in E \vert~o(e)=x\}$ for $x \in V$.
Let $\Omega_{x,n}(X)$ be the set of all paths  $c=(e_{1}, e_2, \ldots, e_{n})$ 
of length $n \in \mathbb{N} \cup\{\infty\}$ starting from $x \in V$. 
For simplicity, we write $\Omega_{x}(X):=\Omega_{x,\infty}(X)$.
Denote by $X_0=(V_0, E_0):=\Gamma \backslash X$ the quotient graph of $X$, 
which is  finite by definition. 

Let  $p:E \to [0,1]$ be
a {\it{$\Gamma$-invariant transition probability}} 
satisfying
$$ \sum_{e\in E_{x}} p(e)=1 \qquad (x\in V) \quad \mbox{ and } \quad p(e)+p({\overline{e}})>0 \qquad
(e\in E).$$
The transition probability $p$ induces, in a natural manner, the probability measure $\mathbb P_{x}$
on the set $\Omega_{x}(X) \, (x \in V)$. The random walk associated with 
$p$ is the time-homogeneous Markov chain
$(\Omega_{x}(X), {\mathbb P}_{x}, \{w_{n} \}_{n=0}^{\infty})$ with values in $X$ defined by 
$w_{n}(c):=o(e_{n+1})$
for $n=0,1,2,\ldots$ and $c=(e_{1},e_2, \ldots, e_{n}, \ldots)\in \Omega_{x}(X)).$
We define the {\it{transition operator}} $L$ associated with the random walk
by
$$Lf(x):=\sum_{e\in E_{x}} p(e) f\big(t(e)\big) \qquad (x\in V,~ f:V \to \mathbb R) $$
and the the $n$-step transition probability by
$p(n,x,y):=L^{n}\delta_{y}(x)$ for $n\in \mathbb N$ and $x,y \in V$,
where $\delta_{y}$ stands for the Dirac delta at $y$.
Since $p$ is $\Gamma$-invariant, the random on $X$ associated with $p$ 
induces a random walk on the quotient graph $X_{0}$ through the covering map $\pi: X\to X_{0}$,
and vice versa. By abuse of the notation, we denote the transition probability of the random walk on $X_{0}$ by 
$p: E_{0} \to [0,1]$. We define the $n$-step transition probability $p(n,x,y)$ ($n\in \mathbb N, \, x,y\in V_{0}$) 
as above.

Throughout the present paper, we assume that the random walk on $X_{0}$ is {\it{irreducible}}, that is,
for every $x,y\in V_{0}$, there exists some $n=n(x,y)\in \mathbb N$ such that $p(n,x,y)>0$. Note that the 
irreducibility of the random walk on $X$ implies this condition. Conversely, this condition does not
imply the irreducibility on $X$ in general. Under the irreducibility condition on $X_{0}$, by applying
the Perron-Frobenius theorem, we find a
unique positive function $m: V_{0} \to (0,1]$, called the {\it{invariant probability 
measure}} on $X_{0}$, satisfying

$$
\sum_{x \in V_0}m(x)=1 \quad\text{and}\quad
 m(x)=\sum_{e \in (E_0)_x}p(\ol{e})m\big(t(e)\big) \qquad (x \in V_0).
$$
Set $\widetilde{m}(e):=p(e)m\big(o(e)\big)$ for $e \in E_0$. 
The random walk on $X_0$ is called ($m$-){\it symmetric} 
if $\widetilde{m}(e)=\widetilde{m}(\ol{e})$ holds for every $e \in E_0$.
Otherwise, it is called ($m$-){\it non-symmetric}. 
We also write $m : V \LA (0, 1]$ for the $\Gamma$-invariant lift of $m : V_0 \LA (0, 1]$. 
Let $\h_1(X_0, \mathbb{R})$  be the first homology group 
of $X_0$. 
We define the {\it homological direction} of the random walk on $X_0$ by 
$$
\gamma_p:=\sum_{e \in E_0}\widetilde{m}(e)e \in \h_1(X_0, \mathbb{R}).
$$
We easily see that the random walk on $X_0$ is $(m$-)symmetric
if and only if $\gamma_p=0$. 
For more details, see Section 3.3.

We introduce a continuous state space 
in which the $\Gamma$-nilpotent covering graph $X$ is properly embedded.
We find a connected and simply connected nilpotent Lie group $\G=(G, \cdot)$ 
such that $\Gamma$ is isomorphic to a cocompact lattice in $G$
due to Malc\'ev's theorem (cf.~Malc\'ev \cite{Malcev}). 
We call a piecewise smooth $\Gamma$-equivariant map 
$\Phi : X \LA G$ a {\it periodic realization} of $X$.
We write $(\g, [\cdot, \cdot])$ for the Lie algebra of $G$. 

We construct a new product $*$ on $G$ in the following way. 
Set
$\frak{n}_1:=\frak{g}$ and $\frak{n}_{k+1}:=[\frak{g}, \frak{n}_k]$ for $k \in \mathbb{N}$.
Then we have
$\frak{g}=\frak{n}_1 \supset \dots 
\supset \frak{n}_r \supsetneq \frak{n}_{r+1}=\{\bm{0}_{\frak{g}}\}$ by definition. 
Define the subspace $\frak{g}^{(k)}$ of $\frak{g}$ by
$\frak{n}_k=\frak{g}^{(k)} \oplus \frak{n}_{k+1}$ for $k=1, 2, \dots, r$.
Then the Lie algebra $\frak{g}$ is decomposed as 
$\frak{g}=\frak{g}^{(1)} \oplus \g^{(2)} \oplus \dots \oplus \frak{g}^{(r)}$ 
and each  $Z \in \frak{g}$ is uniquely written as $Z=Z^{(1)} +Z^{(2)}+ \dots + Z^{(r)}$, 
where 
$Z^{(k)} \in \frak{g}^{(k)}$ for $k=1, 2, \dots, r$.
Then we define a map $\tau^{(\g)}_\ve : \g \LA \g$ by
$\tau^{(\g)}_\ve(Z):=\ve Z^{(1)} + \ve^2 Z^{(2)} + \dots + \ve^r Z^{(r)}$
for $\ve \geq 0$
and also define a Lie bracket product $[\![ \cdot, \cdot ]\!]$ on $\frak{g}$ by
$$
[\![ Z_1, Z_2 ] \!]:=\lim_{\ve \searrow 0} \tau^{(\g)}_\ve 
\big[\tau^{(\g)}_{1/\ve}(Z_1), \tau^{(\g)}_{1/\ve}(Z_2)\big] 
\qquad (Z_1, Z_2 \in \frak{g}).
$$ 
We introduce a map $\tau_\ve : G \LA G$, called the {\it dilation operator} on $G$, by
\begin{equation}\label{deformed product}
\tau_\ve (g):=\Exp \big( \tau^{(\g)}_\ve \big( \Log(g)\big)\big) \qquad (\ve \geq 0, \, g \in G),
\end{equation}
which is the scalar multiplication on $G$. 
Then a Lie group product $*$ on $G$ 
is defined by
$$
g * h :=\lim_{\ve \searrow 0} 
\tau_{\ve}\big( \tau_{1/\ve}(g) \cdot \tau_{1/\ve}(h)\big) \qquad (g, h \in G).
$$
The Lie group $G_\infty=(G, *)$ is called the {\it limit group} of $(G, \cdot)$. 
We note that the Lie algebra of $(G, *)$ coincides with 
$(\g, [\![\cdot, \cdot]\!])$.

We equip $(G, \cdot)$ with the so-called {\it Carnot--Carath\'eodory metric} $d_{\mathrm{CC}}$, 
which is an intrinsic metric on $G$ defined by
\begin{equation}\label{CC-distance}
d_{\mathrm{CC}}(g, h):=\inf\Big\{ \int_0^1 \|\dot{w}_t\|_{g_0^{(0)}} \, dt \, 
\Big| \, \begin{matrix} w \in \mathrm{Lip}([0, 1]; G), \, w_0=g, \, w_1=h, \\ 
\text{$w$ is tangent to $\g^{(1)}$}\end{matrix} \Big\}
\end{equation}
for $g, h \in G$, where $\mathrm{Lip}([0, 1]; G)$ 
denotes the set of all Lipschitz continuous 
paths and $\|\cdot\|_{g_0^{(0)}}$ the Albanese norm associated with $p_0$ on $\g^{(1)}$ defined below. 
See \cite[Section 3.2]{IKN} for more details.  

Let $\pi_1(X_0)$ be the fundamental group of $X_0$. 
Then we take a canonical surjective homomorphism $\rho : \pi_1(X_0) \LA \Gamma$. 
This map then gives rise to a surjective homomorphism 
$\rho : \h_1(X_0, \mathbb{Z}) \LA \Gamma/[\Gamma, \Gamma]$
so that we have
the canonical surjective linear map $\rho_{\mathbb{R}}$
from $\h_1(X_0, \mathbb{R})$ onto $\frak{g}^{(1)}$. 
We call $\rho_{\mathbb{R}}(\gamma_p) \in \frak{g}^{(1)}$
 the {\it asymptotic direction}. 
Note that $\gamma_p=0$ implies 
$\rho_{\mathbb{R}}(\gamma_p)=\bm{0}_{\frak{g}}$.
However, the converse does not always hold.
We induce a flat metric $g_0=g_0^{(1)}$ 
on $\frak{g}^{(1)}$,
called the {\it Albanese metric} associated with  $p$
by the discrete Hodge-Kodaira theorem (cf.~\cite[Lemma 5.2]{KS06}). 


For the given transition probability $p$, 
we introduce a family of $\Gamma$-invariant transition probabilities 
$(p_\ve)_{0 \leq \ve \leq 1}$ on $X$ by
\begin{equation}\label{p_epsilon}
p_\ve(e):=p_0(e)+\ve q(e) \qquad (e \in E),
\end{equation}
where
$$
p_0(e):=\frac{1}{2}\Big( p(e) + \frac{m\big(t(e)\big)}{m\big(o(e)\big)}p(\ol{e})\Big), \quad
q(e):=\frac{1}{2}\Big( p(e) - \frac{m\big(t(e)\big)}{m\big(o(e)\big)}p(\ol{e})\Big).
$$
We note that the family $(p_\ve)_{0 \leq \ve \leq 1}$ is given by the linear interpolation
between the transition probability $p=p_1$ and the $m$-symmetric probability $p_0$.
Moreover, 
the homological direction $\gamma_{p_\ve}$ equals $\ve \gamma_p$ and
the invariant measure $m_\ve$ associated with $p_\ve$ coincides with $m$
for  $0 \leq \ve \leq 1$ (cf.~\cite[Proposition 2.3]{KS06}).

Let $L_{(\ve)}$ be the transition operator associated with $p_\ve$ for $0 \leq \ve \leq 1$. 
We also denote by $g_0^{(\ve)}$ the Albanese metric on $\g^{(1)}$ associated with $p_\ve$. 
We write $G_{(\ve)}$ for the nilpotent Lie group of step $r$ whose Lie algebra is
$\g=(\g^{(1)}, g_0^{(\ve)}) \oplus \g^{(2)} \oplus \cdots \oplus \g^{(r)}$.

We now take the family of periodic realizations  $(\Phi^{(\ve)})_{0 \leq \ve \leq 1}$ satisfying

\vspace{2mm}
\noindent
{\bf (A1)}: {\it For every $0 \leq \ve \leq 1$, 
\begin{equation}\label{ass-mean}
\sum_{x \in \mathcal{F}}m(x) 
 \Log \big( \Phi^{(\ve)}(x)^{-1} \cdot \Phi^{(0)}(x)\big) \big|_{\g^{(1)}}
=0,
\end{equation}
 where $\mathcal{F}$ denotes a fundamental domain of $X$. 
 }
 
 \vspace{2mm}
  \noindent
We note that
it is always possible to take 
$(\Phi^{(\ve)})_{0 \leq \ve \leq 1}$ satisfying {\bf (A1)}.

For a metric space $\T$, we denote by $C_\infty(\T)$ 
the space of continuous functions on $\T$ vanishing at infinity
with the uniform topology $\|\cdot \|_\infty^{\T}$. 
We define an {\it approximation operator} 
$P_\ve : C_\infty(G_{(0)}) \LA C_\infty(X)$ by
$P_\ve f(x):=f\big( \tau_\ve \Phi^{(\ve)}(x)\big)$
for $0 \leq \ve \leq 1$ and $x \in V$.
We extend each element $Z \in \frak{g}$ to a left invariant vector field
$Z_*$ on the stratified Lie group $(G, *)$ in the usual manner. 
We take an orthonormal basis $\{V_1, V_2, \dots, V_{d_1}\}$
of $(\g^{(1)}, g_0^{(0)})$.
Then the first main result is stated as follows:

\begin{tm}\label{CLT2}
{\rm(1)} 
For $0 \leq s \leq t$ and $f \in C_\infty(G_{(0)})$, we have
\begin{equation}\label{semigroup CLT2}
\lim_{n \to \infty}\Big\| L_{(n^{-1/2})}^{[nt]-[ns]}P_{n^{-1/2}} f 
- P_{n^{-1/2}} \e^{-(t-s)\A}f\Big\|_\infty^X=0,
\end{equation}
where $(\e^{-t\A})_{t \geq 0}$ is the $C_0$-semigroup whose infinitesimal generator 
$\A$ is given by
\begin{equation}\label{generator2}
\A=-\frac{1}{2}\sum_{i=1}^{d_1} V_{i_*}^2 - \rho_{\mathbb{R}}(\gamma_p)_*.
\end{equation}

\noindent
{\rm (2)} Let $\mu$ be a Haar measure on $G_{(0)}$. 
Then, for any $f \in C_\infty(G_{(0)})$ and for any sequence 
$\{x_n\}_{n=1}^\infty \subset V$ satisfying
$
\lim_{n \to \infty}\tau_{n^{-1/2}}\big(\Phi^{(n^{-1/2})}(x_n)\big)=:g \in G_{(0)},
$
we have
\begin{equation}\label{semigroup CLT2-2}
\lim_{n \to \infty}L_{(n^{-1/2})}^{[nt]}P_{n^{-1/2}}f(x_n)=\e^{-t\A}f(g)
:=\int_{G_{(0)}} \mathcal{H}_t(h^{-1} * g)f(h) \, \mu(dh) \qquad (t \geq 0),
\end{equation}
where $\mathcal{H}_t(g)$ is a fundamental solution to
$\del u/\del t+\A u=0.$
\end{tm}

We now fix a reference point $x_* \in V$ such that 
$\Phi^{(0)}(x_*)=\bm{1}_G$ and put 
$$
\ol{\xi}_n^{(\ve)}(c):=\Phi^{(\ve)}\big(w_n(c)\big) 
\qquad \big(0 \leq \ve \leq 1, \, n=0, 1, 2, \dots, \, c \in \Omega_{x_*}(X)\big).
$$
Note that {\bf (A1)} does not imply that $\Phi^{(\ve)}(x_*)=\bm{1}_G$ for $0 < \ve \leq 1$
in general. 
We then obtain a $G_{(0)}$-valued random walk 
$(\Omega_{x_*}(X), \mathbb{P}_{x_*}^{(\ve)}, \{\ol{\xi}_n^{(\ve)}\}_{n=0}^\infty)$
associated with the transition probability $p_\ve$. 
For $t \geq 0, \, n=1, 2, \dots$ and $0 \leq \ve \leq 1$, 
let $\ol{\mathcal{X}}_t^{(\ve, n)}$ 
be a map from $\Omega_{x_*}(X)$ to $G_{(0)}$ given by
$$
\ol{\mathcal{X}}_t^{(\ve, n)}(c)
:=\tau_{n^{-1/2}}\big( \ol{\xi}_{[nt]}^{(\ve)}(c)\big) \qquad \big(c \in \Omega_{x_*}(X)\big).
$$
We write 
$\mathcal{D}_n$ for the partition $\{t_k=k/n \, | \, k=0, 1, 2, \dots, n\}$
of the time interval $[0, 1]$ for $n \in \mathbb{N}$. 
We define 
$$
\ol{\mathcal{Y}}_{t_k}^{(\ve, n)}(c)
:=\tau_{n^{-1/2}}\big( \ol{\xi}_{nt_k}^{(\ve)}(c)\big)
=\tau_{n^{-1/2}}\big( \Phi^{(\ve)}(w_{k}(c))\big)
 \qquad \big(t_k \in \mathcal{D}_n, \, c \in \Omega_{x_*}(X)\big)
$$
and also define the $G_{(0)}$-valued continuous stochastic process 
$(\ol{\mathcal{Y}}_t^{(\ve, n)})_{0 \leq t \leq 1} $
given by the $d_{\mathrm{CC}}$-geodesic interpolation of 
$\{ \ol{\mathcal{Y}}_{t_k}^{(\ve, n)}\}_{k=0}^n$. 
We consider a stochastic differential equation
\begin{equation}\label{SDE2}
dY_t
= \sum_{i=1}^{d_1}V_{i_*}^{(0)}(Y_t) \circ dB_t^i 
+ \rho_{\mathbb{R}}(\gamma_p)_*(Y_t) \, dt,
\qquad Y_0=\bm{1}_G,
\end{equation}
where $(B_t)_{0 \leq t \leq 1}=(B_t^1, B_t^2, \dots, B_t^{d_1})_{0 \leq t \leq 1}$ 
is an $\mathbb{R}^{d_1}$-valued standard Brownian motion
starting from $B_0=\bm{0}$. 
It is known that the infinitesimal generator of (\ref{SDE2}) coincides with
$-\A$ defined by (\ref{generator2}) 
(see Ishiwata--Kawabi--Namba \cite[Section 5]{IKN}). 
Let $(Y_t)_{0 \leq t \leq 1}$ be the $G_{(0)}$-valued 
diffusion process which is the solution to (\ref{SDE2}). 
For $\alpha<1/2$, we define the $\alpha$-H\"older distance $\rho_\alpha$ 
on $C([0, 1]; G_{(0)})$ by
$$
 \rho_\alpha(w^1, w^2):=\sup_{0 \leq s < t \leq 1}
\frac{d_{\mathrm{CC}}(u_s, u_t)}{|t-s|^\alpha}+d_{\mathrm{CC}}(\bm{1}_G, u_0),
\qquad u_t:=(w_t^1)^{-1}*w_t^2 \quad (0 \leq t \leq 1). 
$$
Then we set
$C^{0, \alpha\hol}([0, 1]; G_{(0)})
:=\ol{\mathrm{Lip}([0, 1]; G_{(0)})}^{\rho_\alpha}$,
which is a Polish space
(cf.~Friz--Victoir \cite[Section 8]{FV}).   
Let $\ol{\Prob}^{(\ve, n)}$ be 
the probability measure on $C^{0, \alpha\hol}([0, 1]; G_{(0)})$
induced by the stochastic process $\ol{\mathcal{Y}}_\cdot^{(\ve, n)}$
for $0 \leq \ve \leq 1$ and $n \in \mathbb{N}$. 

To establish the functional CLT 
for the family of non-symmetric random walks $\{\ol{\xi}_n^{(\ve)}\}_{n=0}^\infty$, 
we need to impose an additional assumption.

\vspace{2mm}
\noindent
{\bf (A2)}: {\it There exists a positive constant $C$ such that, for $k=2, 3, \dots, r$, 
\begin{equation}\label{ass-dphi}
\sup_{0 \leq \ve \leq 1}\max_{x \in \mathcal{F}}
 \big\| \log\big( \Phi^{(\ve)}(x)^{-1} \cdot \Phi^{(0)}(x)\big) \big|_{\g^{(k)}}\big\|_{\g^{(k)}} \leq C,
\end{equation}
where $\|\cdot\|_{\g^{(k)}}$ denotes a Euclidean norm on $\g^{(k)} \cong \mathbb{R}^{d_k}$
for $k=2, 3, \dots, r$. 
}

\vspace{1mm}
\noindent
Intuitively speaking, 
the situations that the distance between $\Phi^{(\ve)}$ 
and $\Phi^{(0)}$ tends to be too big as $\ve \searrow 0$ are removed under {\bf (A2)}. 
By setting 
$$
\log \big(\Phi^{(\ve)}(x)\big)\big|_{\g^{(k)}}
=\log \big(\Phi^{(0)}(x)\big)\big|_{\g^{(k)}} \qquad (x \in \mathcal{F}, \, k=2, 3, \dots, r)
$$
for $\Phi^{(\ve)} : X \LA G$ with (\ref{ass-mean}), 
the family $(\Phi^{(\ve)})_{0 \leq \ve \leq 1}$ satisfies {\bf (A2)}.
This means that it is always possible to take a family $(\Phi^{(\ve)})_{0 \leq \ve \leq 1}$
satisfying {\bf (A2)} as well as {\bf (A1)}.

Then our main theorem is stated as follows: 
\begin{tm}\label{FCLT2-general} 
We assume {\bf (A1)} and {\bf (A2)}. Then the sequence
$\{\ol{\mathcal{Y}}^{(n^{-1/2}, n)}\}_{n= 1}^\infty $
converges in law to the diffusion process
$(Y_t)_{0 \leq t \leq 1}$ in 
$C^{0, \alpha\text{\normalfont{-H\"ol}}}([0, 1]; G_{(0)})$
as $n \to \infty$ for all $\alpha<1/2$. 
\end{tm}


In our preceding paper {\normalfont \cite{IKN}}, 
we captured another $G$-valued diffusion process $(\widehat{Y}_t)_{0 \leq t \leq 1}$
by applying the transition-shift scheme mentioned in Section 1. 
More precisely, the infinitesimal generator of $(\widehat{Y}_t)_{0 \leq t \leq 1}$
is the homogenized sub-Laplacian on $G$ with a non-trivial drift $\beta(\Phi_0) \in \g^{(2)}$
arising from the non-symmetry of the given random walk, 
where the Lie group $G$ is equipped with the Albanese metric $g_0=g_0^{(1)}$. 
In particular, even in the centered case $\rho_{\mathbb{R}}(\gamma_p)=\bm{0}_{\g}$, 
the non-trivial drift $\beta(\Phi_0)$ remains in general.
On the other hand, in this case, the limiting diffusion $(Y_t)_{0 \leq t \leq 1}$ 
is generated by the homogenized sub-Laplacian on $G_{(0)}$ equipped with the Albanese metric $g_0^{(0)}$.
See  \cite[Remark 5.3]{IKN} for explicit expressions of the limiting diffusions $(Y_t)_{0 \leq t \leq 1}$
and $(\widehat{Y}_t)_{0 \leq t \leq 1}$.



\section{Preliminaries}

\subsection{Some properties on nilpotent Lie groups}
Let us review some properties of the limit group. 
For more details, see e.g., Alexopoulos \cite{A1}, Ishiwata \cite{Ishiwata} and 
Ishiwata--Kawabi--Namba \cite{IKN}. 
Recall that $(G, *)$ is the limit group of a connected and simply connected nilpotent Lie group 
$(G, \cdot)$ of step $r$ (see (\ref{deformed product})). 
We note that $(G, *)$ is {\it stratified}. Namely, 
the decomposition of the corresponding Lie algebra
$(\g=\g^{(1)}\oplus \g^{(2)} \oplus \cdots \oplus \g^{(r)}, [\![\cdot, \cdot]\!])$
satisfies that $[\![ \g^{(i)}, \g^{(j)}]\!] \subset \g^{(i+j)}$ when $i+j \leq r$
and $\frak{g}^{(1)}$ generates  $\frak{g}$. 
We also note that, for $\ve \geq 0$, the dilation $\tau_\ve : G \LA G$ 
is a group automorphism on $(G, *)$ (see \cite[Lemma 2.1]{Ishiwata}). 

We introduce several notations 
that will be used throughout the present paper.

\vspace{1mm}
\noindent
{\bf (1)} {\it Global coordinates on $G$ }: 
We set $d_k=\dim_{\mathbb{R}}\g^{(k)}$ for $k=1, 2, \dots, r$ and $d=\dim G$. 
For $k=1, 2, \dots, r$, we
denote by $\{X_1^{(k)}, X_2^{(k)}, \dots, X_{d_k}^{(k)}\}$
a basis of the subspace $\g^{(k)}$.
We introduce two kinds of 
global coordinate systems in $G$ through $\exp : \g \LA G$.
We identify the nilpotent Lie group $G$ with $\mathbb{R}^d$
as a differentiable manifold by
$$
\begin{aligned}
\mathbb{R}^d \ni &(g^{(1)}, g^{(2)}, \dots,g^{(r)})\\
\longmapsto 
&\,\,g=\Exp\big( g_{d_r}^{(r)}X_{d_r}^{(r)}\big) \cdot \Exp\big( g_{d_r-1}^{(r)}X_{d_r-1}^{(r)}\big)\cdot
\cdots \cdot\Exp\big( g_{1}^{(r)}X_{1}^{(r)}\big)\\
&\hspace{1cm} \cdot \Exp\big( g_{d_{r-1}}^{(r-1)}X_{d_{r-1}}^{(r-1)}\big) 
\cdot \Exp\big( g_{d_{r-1}-1}^{(r-1)}X_{d_{r-1}-1}^{(r-1)}\big) \cdot\cdots
\cdot \Exp\big( g_{1}^{(r-1)}X_{1}^{(r-1)}\big) \\
&\hspace{1cm}\cdots \cdot \Exp\big( g_{d_1}^{(1)}X_{d_1}^{(1)}\big) \cdot \Exp\big( g_{d_1-1}^{(1)}X_{d_1-1}^{(1)}\big)
\cdot \cdots \cdot\Exp\big( g_{1}^{(1)}X_{1}^{(1)}\big) \in G,
\end{aligned}
$$
which is called the {\it canonical $(\cdot)$-coordinates of the second kind}. 
We also define the {\it canonical $(*)$-coordinates of the second kind} just by replacing 
$\cdot$ by the deformed product $*$ in the above correspondence. 

\vspace{1mm}
\noindent
{\bf (2)} {\it Campbell--Baker--Hausdorff formula} : 
The relations between the deformed product and the given product 
on $G$ is described by the following formula:
\begin{equation}\label{CBH-formula}
\log\big(\exp(Z_1) \cdot \exp(Z_2)\big)
=Z_1+Z_2+\frac{1}{2}[Z_1, Z_2]+\cdots \qquad(Z_1, Z_2 \in \g). 
\end{equation}

\vspace{1mm}
\noindent
{\bf (3)} {\it Several formulas} :
It follows from the definition of the deformed product on $G$ that
\begin{equation}\label{rel-1}
\Log(g * h)\big|_{\g^{(k)}}=\Log(g \cdot h)\big|_{\g^{(k)}} \qquad (k=1, 2).
\end{equation}
However, this relation does not hold in general if $k=3, 4, \dots, r$. 
The following identities give a comparison between $(\cdot)$-coordinates and $(*)$-coordinates. 
For $g \in G$, we have the following:
\begin{align}\label{rel-2}
g_{i*}^{(k)}&=g_i^{(k)} \quad (i=1, 2, \dots, d_k, \, k=1, 2), \\
g_{i*}^{(k)}&=g_{i}^{(k)}+\sum_{0 < |K| \leq k-1}C_K \mathcal{P}^K(g)
\qquad (i=1, 2, \dots, d_k, \, k=3, 4, \dots, r) \label{rel-3}
\end{align}
for some constant $C_K$, where $K$ denotes a multi-index 
$\big((i_1, k_1), (i_2, k_2), \dots, (i_\ell, k_\ell)\big)$ 
with length $|K|:=k_1+k_2+\dots+k_\ell$
and $\mathcal{P}^K(g):=g_{i_1}^{(k_1)} \cdot g_{i_2}^{(k_2)} \cdots g_{i_\ell}^{(k_\ell)}$. 
For $g, h \in G$, we also have the following:
\begin{align}\label{rel-4}
(g*h)_{i*}^{(k)}&=(g\cdot h)_i^{(k)} \qquad (i=1, 2, \dots, d_k, \, k=1, 2), \\
(g*h)_{i*}^{(k)}&=(g\cdot h)_{i}^{(k)}+\sum_{\substack{ |K_1|+|K_2|\leq k-1 \\ |K_2|>0}}
C_{K_1, K_2} \mathcal{P}^{K_1}_*(g)\mathcal{P}^{K_2}(g \cdot h)\nn\\
&\hspace{4cm}
\qquad (i=1, 2, \dots, d_k, \, k=3, 4, \dots, r). \label{rel-5}
\end{align}

\vspace{1mm}
\noindent
{\bf (4)} {\it Homogeneous norms on $G$} : 
We introduce a norm $\|\cdot\|_{\g^{(k)}}$ on $\g^{(k)}$ by the usual Euclidean norm.
If $Z \in \g$ is uniquely decomposed by $Z=Z^{(1)}+Z^{(2)}+\cdots +Z^{(r)} \, (Z^{(k)} \in \g^{(k)})$, 
we define a function $\|\cdot\|_{\g} : \g \LA [0, \infty)$ by
$\|Z\|_{\g}:=\sum_{k=1}^r \|Z^{(k)}\|_{\g^{(k)}}^{1/k}.$
The following are the typical examples of homogeneous norms on $G$:
$$
\|g\|_{\mathrm{CC}}:=d_{\mathrm{CC}}(\bm{1}_G, g), \qquad 
\|g\|_{\Hom}:=\|\Log(g)\|_{\g} \qquad (g \in G). 
$$
Note that these homogeneous norms are equivalent 
in the sense that there exists a positive constants $C$ such that 
$C^{-1}\|g\|_{\mathrm{CC}} \leq \|g\|_{\mathrm{Hom}} \leq C\|g\|_{\mathrm{CC}}$
for $g \in G$
(cf.~\cite[Proposition 3.1]{IKN}).

\subsection{A quick review on discrete geometric analysis}
This subsection is concerned with some basics of {\it discrete geometric analysis} 
on graphs. We refer to Kotani--Sunada \cite{KS06}, Sunada \cite{S0, S, S2} 
and references therein for details. 

Let $X_0=(V_0, E_0)$ be a finite graph
and consider an irreducible random walk on $X_0$ associated 
with a transition probability 
$p : E_0 \LA [0, 1]$. 
Let $m : V_0 \LA (0, 1]$ be the normalized invariant probability measure on $X_0$. 
Set $\widetilde{m}(e):=p(e)m\big(o(e)\big)$ for $e \in E_0$. 

We define the 0-chain group, 1-chain group, 0-cochain group and 1-cochain group by
$$
\begin{aligned}
C_0(X_0, \mathbb{R})&:=\Big\{ \sum_{x \in V_0}a_x x \, \Big| \, a_x \in \mathbb{R}\Big\}, &
C_1(X_0, \mathbb{R})&:=\Big\{ \sum_{e \in E_0}a_e e \, \Big| \, a_e \in \mathbb{R}, \, \ol{e}=-e\Big\}, \\
C^0(X_0, \mathbb{R})&:=\{f : V_0 \LA \mathbb{R}\}, &
C^1(X_0, \mathbb{R})&:=\{\omega : E_0 \LA \mathbb{R} \, | \, \omega(\ol{e})=-\omega(e)\},
\end{aligned}
$$
respectively. 
An element of $C^1(X_0, \mathbb{R})$ is called a 1-form on $X_0$. 
The boundary operator $\del : C_1(X_0, \mathbb{R}) \LA C_0(X_0, \mathbb{R})$ 
and the difference operator $d : C^0(X_0, \mathbb{R}) \LA C^1(X_0, \mathbb{R})$ 
are defined by $\del(e)=t(e)-o(e)$ for $e \in E_0$ and 
$df(e)=f\big(t(e)\big)-f\big(o(e)\big)$ for $e \in E_0$, respectively. 
Then, the first homology group $\h_1(X_0, \mathbb{R})$ and 
the first cohomology group $\h^1(X_0, \mathbb{R})$ are defined by 
$\Ker(\del) \subset C_1(X_0, \mathbb{R})$ and $C^1(X_0, \mathbb{R})/\im(d)$, respectively. 
We write $L : C^0(X_0, \mathbb{R}) \LA C^0(X_0, \mathbb{R})$ for 
the transition operator associated with $p$. 
We define a special 1-chain  by
$$
\gamma_p:=\sum_{e \in E_0}\widetilde{m}(e)e \in C_1(X_0, \mathbb{R}).
$$
We easily verify that $\del(\gamma_p)=0$, that is, $\gamma_p \in \h_1(X_0, \mathbb{R})$. 
Furthermore, it is clear that 
the random walk on $X_0$ is ($m$-)symmetric if and only if $\gamma_p=0$.
The 1-cycle $\gamma_p$ is called the {\it homological direction} of the given random walk on $X_0$. 
A 1-form $\omega \in C^1(X_0, \mathbb{R})$ is said to be {\it modified harmonic} if
$\sum_{e \in (E_0)_x}p(e)\omega(x)=\la \gamma_p, \omega \ra$ for $x \in V_0$.
Denote by $\mathcal{H}^1(X_0)$ the set of modified harmonic 1-forms
with the inner product 
$$
\La \omega_1, \omega_2 \Ra_p
:=\sum_{e \in E_0}\widetilde{m}(e)\omega_1(e)\omega_2(e)
 - \la \gamma_p, \omega_1 \ra \la \gamma_p, \omega_2 \ra 
 \qquad \big(\omega_1, \omega_2 \in \mathcal{H}^1(X_0)\big)
$$
associated with the transition probability $p$. 
We may identify $\h^1(X_0, \mathbb{R})$
with $\mathcal{H}^1(X_0)$ 
by the discrete Hodge-Kodaira theorem (cf.~\cite[Lemma 5.2]{KS06}).
Then the inner product $\La \cdot, \cdot \Ra_p$ is induced on $\h^1(X_0, \mathbb{R})$ 
through this identification.

Let $\Gamma$ be a torsion free, finitely generated nilpotent group of step $r$.
Then a $\Gamma$-nilpotent covering graph $X=(V, E)$ 
is defined by the $\Gamma$-covering of $X_0$.
Let $p : E \LA [0, 1]$ and $m : V \LA (0, 1]$ be the $\Gamma$-invariant lifts of  
$p : E_0 \LA [0, 1]$ and $m : V_0 \LA (0, 1]$, respectively. 
Denote by $\widehat{\pi} : G \LA G/[G, G]$ the canonical projection. 
Since $\Gamma$ is a cocompact lattice in $G$,
the subset $\widehat{\pi}(\Gamma) \subset G/[G, G]$
 is also a lattice in $G/[G, G] \cong \g^{(1)}$ (cf.~Malc\'ev \cite{Malcev} and Raghunathan \cite{Rag}). 
Let $\rho_{\mathbb{R}} : \h_1(X_0, \mathbb{R}) 
\LA \widehat{\pi}(\Gamma) \otimes \mathbb{R} \cong \g^{(1)}$
be the canonical surjective linear map induced by the surjective homomorphism
$\rho : \pi_1(X_0) \LA \Gamma$. 
We restrict the inner product $\La \cdot, \cdot \Ra_p$ on $\h^1(X_0, \mathbb{R})$ 
to the subspace $\Hom(\widehat{\pi}(\Gamma), \mathbb{R})$, 
thereafter take it up the dual inner product $\la \cdot , \cdot \ra_{alb}$ 
on $\widehat{\pi}(\Gamma) \otimes \mathbb{R}$. 
Consequently, we induce a flat metric $g_0$ on $\g^{(1)}$
and call it the {\it Albanese metric} on $\g^{(1)}$. 
This procedure can be summarized as follows: 

$$
\xymatrix{ 
(\frak{g}^{(1)}, g_0)  \ar @{<->}[d]^{\mathrm{dual}}
&\hspace{-1.5cm}\cong 
& \hspace{-1.5cm}\qquad\widehat{\pi}(\Gamma) \otimes \mathbb{R}
 \ar @{<<-}[r]^{\rho_{\mathbb{R}}} \ar @{<->}[d]^{\mathrm{dual}} 
 & \h_1(X_0, \mathbb{R}) \ar @{<->}[d]^{\mathrm{dual}} &\\
\Hom(\frak{g}^{(1)}, \mathbb{R}) 
&\hspace{-1.5cm}\cong 
& \hspace{-1.5cm}\quad\Hom(\widehat{\pi}(\Gamma), \mathbb{R})  
\ar @{^{(}->}[r]_{{}^t \rho_{\mathbb{R}}} & \h^1(X_0, \mathbb{R})   
&\hspace{-0.9cm}\cong  \big(\mathcal{H}^1(X_0), \La \cdot , \cdot \Ra_p\big).}
$$



\section{Proof of main results}

\subsection{A one-parameter family of modified harmonic realizations $(\Phi_0^{(\ve)})_{0 \leq \ve \leq 1}$}

Recall that $(p_\ve)_{0 \leq \ve \leq 1}$ is the family of 
transition probabilities defined by (\ref{p_epsilon}). 
We now introduce the family of modified harmonic realizations 
$(\Phi_0^{(\ve)})_{0 \leq \ve \leq 1}$. 
Namely, each $\Phi_0^{(\ve)}$ is the $\Gamma$-equivariant realization of $X$ satisfying 
\begin{equation}\label{p_ve-modified}
\sum_{e \in E_x}p_\ve(e) \Log \Big( \Phi_0^{(\ve)}\big(o(e)\big)^{-1} 
\cdot \Phi_0^{(\ve)}\big(t(e)\Big)\Big|_{\g^{(1)}}
=\ve \rho_{\mathbb{R}}(\gamma_p) \qquad (x \in V).
\end{equation} 
Moreover, we may assume that $\Phi_0^{(0)}(x_*)=\bm{1}_G$ and 
$$
\log \big(\Phi_0^{(\ve)}(x)\big)\big|_{\g^{(k)}}
=\log \big(\Phi^{(\ve)}(x)\big)\big|_{\g^{(k)}} \qquad (0 \leq \ve \leq 1, \, x \in \mathcal{F}, \, k=2, 3, \dots, r)
$$
without loss of generality. 
We define the ($\g^{(1)}$-){\it corrector} $\mathrm{Cor}^{(\ve)}_{\g^{(1)}} : X \LA \g^{(1)}$ by
\begin{equation}\label{corrector}
\mathrm{Cor}^{(\ve)}_{\g^{(1)}}(x):=\log\big(\Phi^{(\ve)}(x)\big)\big|_{\g^{(1)}} - 
\log\big(\Phi_0^{(\ve)}(x)\big)\big|_{\g^{(1)}} \qquad (x \in V, \, 0 \leq \ve \leq 1).
\end{equation}
We note that, thanks to {\bf (A1)}, we have
\begin{equation}\label{corrector}
\sum_{x \in \mathcal{F}}m(x)\mathrm{Cor}_{\g^{(1)}}^{(\ve)}(x)
=\sum_{x \in \mathcal{F}}m(x)\mathrm{Cor}_{\g^{(1)}}^{(0)}(x) \qquad (0 \leq \ve \leq 1).
\end{equation}
In particular, there exists a positive constant $M>0$ independent of $\ve \in [0, 1]$ such that 
$\max_{x \in \mathcal{F}}\|\mathrm{Cor}_{\g^{(1)}}^{(\ve)}(x)\|_{\g^{(1)}} \leq M$ for $0 \leq \ve \leq 1$. 
We also emphasize that, 
if $(\Phi^{(\ve)})_{0 \leq \ve \leq 1}$ satisfies {\bf (A1)} and {\bf (A2)}, then
$(\Phi_0^{(\ve)})_{0 \leq \ve \leq 1}$ also satisfies 
{\bf (A1)} and {\bf (A2)}, respectively. Indeed, 
by combining (\ref{corrector}) and {\bf (A1)}, 
we see that 
$$
\begin{aligned}
&\sum_{x \in \mathcal{F}}m(x) 
 \Log \big( \Phi_0^{(\ve)}(x)^{-1} \cdot \Phi_0^{(0)}(x)\big) \big|_{\g^{(1)}}\\
&=\sum_{x \in \mathcal{F}}m(x)\mathrm{Cor}_{\g^{(1)}}^{(\ve)}(x) 
-\sum_{x \in \mathcal{F}}m(x)\mathrm{Cor}_{\g^{(1)}}^{(0)}(x) +
\sum_{x \in \mathcal{F}}m(x) 
 \Log \big( \Phi^{(\ve)}(x)^{-1} \cdot \Phi^{(0)}(x)\big) \big|_{\g^{(1)}}=0\\
\end{aligned}
$$
for $0 \leq \ve \leq 1$, 
which means that that the family $(\Phi_0^{(\ve)})_{0 \leq \ve \leq 1}$ enjoys the assumption {\bf (A1)}. 
Furthermore, by using {\bf (A2)} and 
$\Phi_0^{(\ve)}(x)^{(i)}=\Phi^{(\ve)}(x)^{(i)}$ for $x \in V, \, 0 \leq \ve \leq 1$ 
and $i=2, 3, \dots, r$, we have
\begin{equation}\label{final bound}
\sup_{0 \leq \ve \leq 1}\max_{x \in \mathcal{F}}
 \big\| \log\big( \Phi_0^{(\ve)}(x)^{-1} \cdot \Phi_0^{(0)}(x)\big) \big|_{\g^{(k)}}\big\|_{\g^{(k)}} \leq C
\end{equation}
for some $C>0$, which implies that 
the family $(\Phi_0^{(\ve)})_{0 \leq \ve \leq 1}$ satisfies the assumption {\bf (A2)}.

We put 
$d\Phi^{(\ve)}_0(e)=\Phi_0^{(\ve)}\big(o(e)\big)^{-1} \cdot \Phi_0^{(\ve)}\big(t(e)\big)$
for $0 \leq \ve \leq 1$ and $e \in E$. 
The aim of this subsection is to study the quantity
$$
\beta_{(\ve)}(\Phi_0^{(\ve)}):=
\sum_{e \in E_0}\widetilde{m}_\ve(e) 
\Log \big( d\Phi_0^{(\ve)}(\widetilde{e})\big)\big|_{\g^{(2)}} \in \g^{(2)} \qquad (0 \leq \ve \leq 1),
$$
where we put $\widetilde{m}_\ve(e)=p_{\ve}(e)m\big(o(e)\big)$ for $e \in E_0$. 
Note that, if the transition probability $p_0$ is $m$-symmetric, 
then $\beta_{(0)}(\Phi_0^{(0)})=\bm{0}_{\g}$.
Loosely speaking, this quantity will appear as a coefficient of the second order term
of the Taylor expansion of $(I-L_{(\ve)}^N )P_\ve f$ in $\ve$,
which is dealt in the proof of Lemma \ref{key-lem2}.
In particular, we are interested in the short time behavior of $\beta_{(\ve)}(\Phi_0^{(\ve)})$
as $\ve \searrow 0$ for later use. 
Intuitively there seems to be little hope of knowing such behavior, 
because $\Phi_0^{(\ve)}$ has the ambiguity in its $\g^{(2)}$-components 
for every $0 \leq \ve \leq 1$.
However, the following proposition asserts that
$\beta_{(\ve)}(\Phi_0^{(\ve)})$ in fact approaches 
$\beta_{(0)}(\Phi_0^{(0)})=\bm{0}_{\g}$
as $\ve \searrow 0$ by imposing only {\bf (A1)}. 

\begin{pr} \label{key-beta}
Under {\bf (A1)}, we have
$$\displaystyle
\lim_{\ve \searrow 0} \beta_{(\ve)}(\Phi_0^{(\ve)}) = \beta_{(0)}(\Phi_0^{(0)})=\bm{0}_{\g}.$$
\end{pr} 

Fix a fundamental domain $\mathcal{F}$ of $X$. 
Set 
$\Psi^{(\ve)}(x)=\Phi_0^{(\ve)}(x)^{-1} \cdot \Phi_0^{(0)}(x)$ for
$0 \leq \ve \leq 1$ and $x \in V$.
Note that the map $\Psi^{(\ve)} : V \LA G$ is $\Gamma$-invariant. 
The following lemma is essential to prove Proposition \ref{key-beta}. 

\begin{lm}\label{lem for beta}
Under {\bf (A1)}, we have
$$
\lim_{\ve \searrow 0} 
\big\| \Log \big( \Psi^{(\ve)}(x)\big) \big|_{\g^{(1)}}\big\|_{\g^{(1)}}=0 \qquad (x \in \mathcal{F}).
$$
In particular, there exists a constant $C$ such that 
$$
\big\| \Log \big( \Psi^{(\ve)}(x)\big) \big|_{\g^{(1)}}\big\|_{\g^{(1)}}
\leq C \qquad (0 \leq \ve \leq 1, \, x \in \mathcal{F}). 
$$
\end{lm}

\vspace{3mm}
\noindent
{\bf Proof.} We set $\ell^2(\mathcal{F}):=\{f : \mathcal{F} \LA \mathbb{C}\}$
 and equip it with the inner product and the corresponding norm defined by
$$
\la f, g \ra_{\ell^2(\mathcal{F})}:=
\sum_{x \in \mathcal{F}}f(x) \ol{g(x)}, \quad 
\|f\|_{\ell^2(\mathcal{F})}:=\Big(\sum_{x \in \mathcal{F}} |f(x)|^2\Big)^{1/2} 
 \qquad \big(f, g \in \ell^2(\mathcal{F})\big),
$$
respectively. 
Since the invariant measure $m|_{\mathcal{F}} : \mathcal{F} \LA (0, 1]$ 
is positive on the finite set $\mathcal{F}$, 
there are positive constants $c$ and $C$
such that 
\begin{equation}\label{cC}
c \Big(\sum_{x \in \mathcal{F}}m(x) |f(x)|^2\Big)^{1/2} 
\leq \|f\|_{\ell^2(\mathcal{F})} 
\leq C \Big(\sum_{x \in \mathcal{F}}m(x) |f(x)|^2\Big)^{1/2} 
\qquad \big( f \in \ell^2(\mathcal{F}) \big).
\end{equation}
It follows from the Perron--Frobenius theorem that 
$\ell^2(\mathcal{F})=\la \phi_0 \ra \oplus \ell_1(\mathcal{F})$,
where $\phi_0 = |\mathcal{F}|^{-1/2}$ is 
the normalized right eigenfunction corresponding to 
the maximal eigenvalue $\alpha_0=1$ of $L$. 
We define
$\ell_1^2(\mathcal{F}):=\big\{f \in \ell^2(\mathcal{F}) \, : \, 
|\mathcal{F}|^{1/2}\la f, m \ra_{\ell^2(\mathcal{F})}=0 \big\}.$
Note that $\ell_1^2(\mathcal{F})$ is preserved by $L$ 
and the transition operator $L_{(\ve)}$ 
maps $\ell_1^2(\mathcal{F})$ to itself for all $0 \leq \ve \leq 1$.
Moreover, the inverse operator of 
$(I-L_{(\ve)})|_{\ell_1^2(\mathcal{F})} : \ell_1^2(\mathcal{F}) \LA \ell_1^2(\mathcal{F})$ 
does exists since
 $L_{(\ve)}$ has a simple eigenvalue $\alpha_0(\ve)=1$ for $0 \leq \ve \leq 1$. 
We define $Q : \ell^2(\mathcal{F}) \LA \ell^2(\mathcal{F})$ by
$$
Qf(x):=\sum_{e \in E_x}q(e) f\big(t(e)\big) 
\qquad \big( f \in \ell^2(\mathcal{F}), \, x \in \mathcal{F}\big).
$$
Then we verify that the transition operator $L_{(\ve)}$ has 
the decomposition of the form $L_{(\ve)}=L_{(0)}+\ve Q$ for every $0 \leq \ve \leq 1$. 

To conclude the claim, it suffices to show 
\begin{equation}\label{lem-want-to-show}
\lim_{\ve \searrow 0} \big\| \Log \big( \Psi^{(\ve)}(\cdot)\big) \big|_{X_i^{(1)}} 
\big\|_{\ell^2(\mathcal{F})}=0 
\qquad (i=1, 2, \dots, d_1)
\end{equation}
by noting (\ref{cC}). We remark that 
$\Log \big( \Psi^{(\ve)}(\cdot)\big) \big|_{X_i^{(1)}}  
\in \ell_1^2(\mathcal{F})$ for $i=1, 2, \dots, d_1$ 
thanks to (\ref{ass-mean}). 
In the following, we fix $i=1, 2, \dots, d_1$. 
The modified harmonicity of $\Phi_0^{(\ve)}$ gives
$$
\begin{aligned}
 (I-L_{(\ve)})\Big( \Log \big( \Psi^{(\ve)}(x)\big) \big|_{X_i^{(1)}} \Big)
&=\ve \Big[ Q\big( \Log\Phi_0^{(0)}(x)\big|_{X_i^{(1)}} \big)
- \rho_{\mathbb{R}}(\gamma_p)\big|_{X_i^{(1)}}\Big] 
\end{aligned}
$$
for $0 \leq \ve \leq 1$ and $x \in \mathcal{F}$. 
This identity implies
\begin{align}\label{estimate1}
&\big\|\Log \big( \Psi^{(\ve)}(\cdot)\big) \big|_{X_i^{(1)}} 
\big\|_{\ell^2(\mathcal{F})} \nn\\
&\leq \ve \big\| (I-L_{(\ve)})\big|_{\ell_1^2(X_0)}^{-1}\big\| \cdot 
          \big\| Q\big( \Log\Phi_0^{(0)}(\cdot)\big|_{X_i^{(1)}} \big)
          - \rho_{\mathbb{R}}(\gamma_p)\big|_{X_i^{(1)}}\big\|_{\ell^2(\mathcal{F})}\nn\\
& \leq \ve \big\| (I-L_{(\ve)})\big|_{\ell_1^2(X_0)}^{-1}\big\| \cdot 
           \Big\{ \big\| \Log\Phi_0^{(0)}(\cdot)\big|_{X_i^{(1)}} \big\|_{\ell^2(\mathcal{F})} 
           + \big\|\rho_{\mathbb{R}}(\gamma_p)\big\|_{\g^{(1)}}\Big\},
\end{align}
where we used $\|Q\| \leq 1$ for the final line. 
By combining (\ref{estimate1}) with the identity
$$
(I-L_{(\ve)})\big|_{\ell_1^2(\mathcal{F})}^{-1} 
= (I-L_{(0)})\big|_{\ell_1^2(\mathcal{F})}^{-1} 
\Big[ I-\ve Q\big|_{\ell_1^2(\mathcal{F})}(I-L_{(0)})\big|_{\ell_1^2(\mathcal{F})}^{-1}\Big],
$$
we obtain
$$
\begin{aligned}
 \big\|\Log \big( \Psi^{(\ve)}(\cdot)\big) \big|_{X_i^{(1)}} \big\|_{\ell^2(\mathcal{F})} 
 &\leq \ve \big\| (I-L_{(0)})\big|_{\ell_1^2(\mathcal{F})}^{-1}\big\| \cdot 
           \Big( 1-\ve \big\|Q\big|_{\ell_1^2(\mathcal{F})}(I-L_{(0)})\big|_{\ell_1^2(\mathcal{F})}^{-1}
           \big\|\Big)^{-1}\\
           &\hspace{1cm} \times \Big\{ \big\| \Log\Phi_0^{(0)}(\cdot)\big|_{X_i^{(1)}}  \big\|_{\ell^2(\mathcal{F})} 
           + \big\|\rho_{\mathbb{R}}(\gamma_p)\big\|_{\g^{(1)}}\Big\}.
\end{aligned}
$$
Here we can choose a sufficiently small constant $\ve_0>0$ such that 
$$
\sup_{0 \leq \ve \leq \ve_0}
\Big( 1-\ve \big\|Q\big|_{\ell_1^2(\mathcal{F})}(I-L_{(0)})
\big|_{\ell_1^2(\mathcal{F})}^{-1}\big\|\Big)^{-1} \leq 2.
$$
Then we have 
$$
\begin{aligned}
 \big\|\Log \big( \Psi^{(\ve)}(\cdot)\big) \big|_{X_i^{(1)}} \big\|_{\ell^2(\mathcal{F})} 
&\leq 2\ve \big\| (I-L_{(0)})\big|_{\ell_1^2(\mathcal{F})}^{-1}\big\| 
            \Big\{ \big\| \Log\Phi_0^{(0)}(\cdot)\big|_{X_i^{(1)}}  \big\|_{\ell^2(\mathcal{F})} 
            + \big\|\rho_{\mathbb{R}}(\gamma_p)\big\|_{\g^{(1)}}\Big\}
\end{aligned}
$$
for sufficiently small $\ve>0$ and this implies (\ref{lem-want-to-show}). \qed

\vspace{2mm}


\vspace{3mm}
\noindent
{\bf Proof of Proposition \ref{key-beta}.} 
By recalling (\ref{p_epsilon}) and that $p_0$ is $m$-symmetric, we have

$$
\begin{aligned}
\beta_{(\ve)}(\Phi_0^{(\ve)}) 
&=\sum_{e \in E_0}\Big\{\frac{1}{2}\big( \widetilde{m}_0(e) - \widetilde{m}_0(\ol{e})\big)
\log \big(d \Phi_0^{(\ve)}(\widetilde{e})\big)\big|_{\g^{(2)}}
+\ve m\big(o(e)\big)q(e)
\log \big(d \Phi_0^{(\ve)}(\widetilde{e})\big)\big|_{\g^{(2)}}\Big\}\\
&=\ve \sum_{e \in E_0}m\big(o(e)\big)q(e)
\log \big(d \Phi_0^{(\ve)}(\widetilde{e})\big)\big|_{\g^{(2)}}. 
\end{aligned}
$$
Then the identity
\begin{equation}\label{identity-dphi}
d\Phi_0^{(\ve)}(e)
=\Psi^{(\ve)}\big(o(e)\big) \cdot d\Phi_0^{(0)}(e) 
\cdot \Psi^{(\ve)}\big(t(e)\big)^{-1} \qquad (0 \leq \ve \leq 1, \, e \in E),
\end{equation}
and (\ref{CBH-formula}) yield

\begin{align}\label{beta-expansion}
\beta_{(\ve)}(\Phi_0^{(\ve)}) &=\ve \sum_{e \in E_0}m\big(o(e)\big)q(e)  
\Big\{ \Log\big( \Psi^{(\ve)}\big(o(\widetilde{e})\big)\big)\big|_{\g^{(2)}} 
-  \Log \big( \Psi^{(\ve)}\big(t(\widetilde{e})\big)\big)\big|_{\g^{(2)}}\Big\}\nn\\
&\hspace{0.8cm}+ \ve \sum_{e \in E_0}m\big(o(e)\big)q(e)   
\Log \big( d\Phi_0^{(0)}(\widetilde{e})\big)\big|_{X_i^{(2)}} \nn\\
&\hspace{0.8cm} -\frac{\ve}{2} \sum_{e \in E_0}m\big(o(e)\big)q(e)   
\Big\{ \mathcal{I}_1^{(\ve)}(\widetilde{e})
+\mathcal{I}_2^{(\ve)}(\widetilde{e})
+\mathcal{I}_3^{(\ve)}(\widetilde{e})\Big\}, 
\end{align}
where 
$$
\begin{aligned}
  \mathcal{I}_1^{(\ve)}(\widetilde{e})=\mathcal{I}_1^{(\ve; \lambda, \nu)}(\widetilde{e})
 &=\Big[ \log \big( \Psi^{(\ve)}(o(\widetilde{e}))\big)\big|_{\g^{(1)}}, 
      \log \big(d\Phi_0^{(0)}(\widetilde{e})\big)\big|_{\g^{(1)}}\Big],\\
  \mathcal{I}_2^{(\ve)}(\widetilde{e})=\mathcal{I}_2^{(\ve; \lambda, \nu)}(\widetilde{e})
 &= \Big[ \log \big( \Psi^{(\ve)}(o(\widetilde{e}))\big)\big|_{\g^{(1)}}, 
      \log \big(\Psi^{(\ve)}(t(\widetilde{e}))^{-1}\big)\big|_{\g^{(1)}}\Big],\\
 \mathcal{I}_3^{(\ve)}(\widetilde{e})=\mathcal{I}_3^{(\ve; \lambda, \nu)}(\widetilde{e})
 &= \Big[\log \big(d\Phi_0^{(0)}(\widetilde{e})\big)\big|_{\g^{(1)}}, 
       \log \big(\Psi^{(\ve)}(t(\widetilde{e}))^{-1}\big)\big|_{\g^{(1)}}\Big]. 
\end{aligned}
$$
Let $\{X_1^{(2)}, X_2^{(2)}, \dots, X_{d_2}^{(2)}\}$ be a basis of $\g^{(2)}$. 
For $i=1, 2, \dots, d_2$, 
we define a function $F_i^{(\ve)} : V \LA \mathbb{R}$ by 
$F_i^{(\ve)}(x):=\log\big(\Psi^{(\ve)}(x)\big)\big|_{X_i^{(2)}}$
for $0 \leq \ve \leq 1$ and $x \in V$. 
Then we see that the function $F_i^{(\ve)}$ is $\Gamma$-invariant. 
Hence, there exists a function $\widehat{F}^{(\ve)} : V_0 \LA \mathbb{R}$
such that $\widehat{F}_i^{(\ve)}\big(\pi(x)\big)=F_i^{(\ve)}(x)$ for $0 \leq \ve \leq 1$ and $x \in V$.
Then, by noting $\del(\gamma_{p_\ve})=0$, we have
$$
\begin{aligned}
&\ve \sum_{e \in E_0}m\big(o(e)\big)q(e)  
\Big\{ \Log\big( \Psi^{(\ve)}\big(o(\widetilde{e})\big)\big)\big|_{\g^{(2)}} 
-  \Log \big( \Psi^{(\ve)}\big(t(\widetilde{e})\big)\big)\big|_{\g^{(2)}}\Big\}\\
&=\sum_{e \in E_0}\big(\widetilde{m}_\ve(e) - \widetilde{m}_0(e)\big)  
\Big\{ \Log\big( \Psi^{(\ve)}\big(o(\widetilde{e})\big)\big)\big|_{\g^{(2)}} 
-  \Log \big( \Psi^{(\ve)}\big(t(\widetilde{e})\big)\big)\big|_{\g^{(2)}}\Big\}\\
&= - {}_{C_1(X_0, \mathbb{R})}\big\la \gamma_{p_\ve}, 
d\widehat{F}_i^{(\ve)}
 \big\ra_{C^1(X_0, \mathbb{R})}
 +\frac{1}{2}\sum_{e \in E_0}\big(\widetilde{m}_0(e) - \widetilde{m}_0(\ol{e})\big)
 d\widehat{F}_i^{(\ve)}(e)\\
&= -{}_{C_0(X_0, \mathbb{R})}\big\la 
\del(\gamma_{p_\ve}), \widehat{F}_i^{(\ve)}
\big\ra_{C^0(X_0, \mathbb{R})}=0.
\end{aligned}
$$
By applying Lemma \ref{lem for beta} and the elementary inequality 
$\|[Z_1, Z_2]\big\|_{\g^{(2)}} \leq C \|Z_1\|_{\g^{(1)}}\|Z_2\|_{\g^{(1)}}$
for $Z_1, Z_2 \in \g^{(1)}$ and some $C>0$, there exists a constant $C>0$ such that 
$\big\|\mathcal{I}_k^{(\ve)}(\widetilde{e})\big\|_{\g^{(2)}} \leq C$ 
for $0 \leq \ve \leq 1$ and $k=1, 2, 3$. 
Summing up the all arguments above and letting $\ve \searrow 0$ in both sides of (\ref{beta-expansion}), 
we obtain the desired convergence. 
This completes the proof of Proposition \ref{key-beta}. \qed

\vspace{2mm}



\subsection{Proof of Theorem \ref{CLT2}}

The following Lemma plays a crucial role in proving Theorem \ref{CLT2}.

\begin{lm}\label{key-lem2}
Let $P_\ve^{H} : C_\infty(G_{(0)}) \to C_\infty(X) \, $ 
be an approximation operator
defined by $P_\ve^H f(x):=f\big(\tau_\ve\Phi_0^{(\ve)}(x)\big)$ for $0 \leq \ve \leq 1$ and $x \in V$. 
Then, for any $f \in C_0^\infty(G_{(0)})$, as $N \to \infty$ and $\ve \searrow 0$ 
with $N^2\ve \searrow 0$, we have
$$
\Big\| \frac{1}{N\ve^2}(I-L_{(\ve)}^N)P^H_\ve f 
- P^H_\ve \A f\Big\|_\infty^X \LA 0,
$$
where $\A$ is the sub-elliptic operator on $C_0^\infty(G_{(0)})$ defined by \normalfont{(\ref{generator2})}. 
\end{lm}

\noindent
{\bf Proof.} 
We apply Taylor's expansion formula (cf.~Alexopoulos \cite[Lemma 5.3]{A2})
for the ($*$)-coordinates of the second kind to $f \in C_0^\infty(G_{(0)})$ at 
$\tau_\ve \big( \Phi_0^{(\ve)}(x)\big) \in G_{(0)}$. Then, recalling that 
$(G_{(0)}, *)$ is a stratified Lie group, we have
\begin{align}\label{keylem2-eq1}
&\frac{1}{N\ve^2}(I - L_{(\ve)}^N)P^H_\ve f(x)\nn\\
&=-\sum_{(i, k)}\frac{\ve^{k-2}}{N}X_{i*}^{(k)}f\Big(\tau_\ve\big(\Phi_0^{(\ve)}(x)\big)\Big)
 \sum_{c \in \Omega_{x, N}(X)}p_\ve(c)\Big(\Phi_0^{(\ve)}(x)^{-1}
 *\Phi_0^{(\ve)}\big(t(c)\big)\Big)_{i*}^{(k)}\nn\\
&\hspace{0.5cm} -\Big(\sum_{(i_1, k_1) \geq (i_2, k_2)}
\frac{\ve^{k_1+k_2-2}}{2N}X_{i_1*}^{(k_1)}X_{i_2*}^{(k_2)}
+ \sum_{(i_2, k_2) > (i_1, k_1)}
\frac{\ve^{k_1+k_2-2}}{2N}X_{i_2*}^{(k_2)}X_{i_1*}^{(k_1)}\big)\Big)
f\Big(\tau_\ve\big(\Phi_0^{(\ve)}(x)\big)\Big)\nn\\
&\hspace{1cm}\times \sum_{c \in \Omega_{x, N}(X)} p_\ve(c)
\Big(\Phi_0^{(\ve)}(x)^{-1}*\Phi_0^{(\ve)}\big(t(c)\big)\Big)_{i_1*}^{(k_1)}
\Big(\Phi_0^{(\ve)}(x)^{-1}*\Phi_0^{(\ve)}\big(t(c)\big)\Big)_{i_2*}^{(k_2)} \nn\\
&\hspace{0.5cm} -\sum_{(i_1, k_1), (i_2, k_2), (i_3, k_3)} \frac{\ve^{k_1+k_2+k_3-2}}{6N}
\frac{\del^3 f}{\del g_{i_1*}^{(k_1)}\del g_{i_2*}^{(k_2)}\del g_{i_3*}^{(k_3)}}(\theta)\nn\\
&\hspace{1cm} \times 
\sum_{c \in \Omega_{x, N}(X)} p_\ve(c)
\Big(\Phi_0^{(\ve)}(x)^{-1}*\Phi_0^{(\ve)}\big(t(c)\big)\Big)_{i_1*}^{(k_1)}\nn\\
&\hspace{1cm} \times \Big(\Phi_0^{(\ve)}(x)^{-1}*\Phi_0^{(\ve)}\big(t(c)\big)\Big)_{i_2*}^{(k_2)}
\Big(\Phi_0^{(\ve)}(x)^{-1}*\Phi_0^{(\ve)}\big(t(c)\big)\Big)_{i_3*}^{(k_3)}
\end{align}
for $x \in V$ and some $\theta \in G_{(0)}$ satisfying 
$$
|\theta_{i_*}^{(k)}| \leq \Big(\Phi_0^{(\ve)}(x)^{-1}*\Phi_0^{(\ve)}\big(t(c)\big)\Big)_{i*}^{(k)}
\qquad (i=1, 2, \dots, d_k, \, k=1, 2, \dots r),
$$
where the summation $\sum_{(i_1, k_1) \geq (i_2, k_2)}$ 
runs over all $(i_1, k_1)$ and $(i_2, k_2)$
with $k_1>k_2$ or $k_1=k_2$ and $i_1 \geq i_2$. 
We denote by $\mathrm{Ord}_\ve(k)$ the terms of 
the right-hand side of (\ref{keylem2-eq1})
whose order of $\ve$ equals just $k$. 
Then, (\ref{keylem2-eq1}) is rewritten as
$$
\frac{1}{N\ve^2}(I - L_{(\ve)}^N)P^H_\ve f(x)=\mathrm{Ord}_\ve(-1)+\mathrm{Ord}_\ve(0)
+\sum_{k \geq 1} \mathrm{Ord}_\ve(k) \qquad (x \in V),
$$
where 
$$
\begin{aligned}
\mathrm{Ord}_\ve(-1) &=
-\frac{1}{N\ve}\sum_{i=1}^{d_1}X_{i*}^{(1)}f\big(\tau_\ve\big(\Phi_0^{(\ve)}(x)\big)\big)
\sum_{c \in \Omega_{x, N}(X)}p_\ve(c)\Big(\Phi_0^{(\ve)}(x)^{-1}
 *\Phi_0^{(\ve)}\big(t(c)\big)\Big)_{i*}^{(1)}
 \end{aligned}
 $$
 and 
 $$
 \begin{aligned}
 \mathrm{Ord}_\ve(0) &=
-\frac{1}{N}\sum_{i=1}^{d_2}X_{i*}^{(2)}f\big(\tau_\ve\big(\Phi_0^{(\ve)}(x)\big)\big)
\sum_{c \in \Omega_{x, N}(X)}p_\ve(c)\Big\{\Big(\Phi_0^{(\ve)}(x)^{-1}
 *\Phi_0^{(\ve)}\big(t(c)\big)\Big)_{i*}^{(2)}\\
 &\hspace{1cm}-\frac{1}{2}\sum_{1 \leq \lambda < \nu \leq d_1}
 [\![X_{\lambda}^{(1)}, X_\nu^{(1)}]\!]\big|_{X_i^{(2)}}\\
 &\hspace{1cm}\times\Big(\Phi_0^{(\ve)}(x)^{-1}
 *\Phi_0^{(\ve)}\big(t(c)\big)\Big)_{\lambda*}^{(1)}\Big(\Phi_0^{(\ve)}(x)^{-1}
 *\Phi_0^{(\ve)}\big(t(c)\big)\Big)_{\nu*}^{(1)}\Big\}\\
 &\hspace{1cm}-\frac{1}{2N}\sum_{1 \leq i, j \leq d_1}X_{i*}^{(1)}X_{j*}^{(1)}f\big(\tau_\ve\big(\Phi_0^{(\ve)}(x)\big)\big)\\
 &\hspace{1cm}\times 
 \sum_{c \in \Omega_{x, N}(X)}p_\ve(c)
 \Big(\Phi_0^{(\ve)}(x)^{-1}
 *\Phi_0^{(\ve)}\big(t(c)\big)\Big)_{i*}^{(1)}\Big(\Phi_0^{(\ve)}(x)^{-1}
 *\Phi_0^{(\ve)}\big(t(c)\big)\Big)_{j*}^{(1)}\\
 &=:\mathcal{I}_1(\ve, N) + \mathcal{I}_2(\ve, N).
\end{aligned}
$$

\noindent
{\bf Step~1.} We first estimate $\mathrm{Ord}_\ve(-1)$. 
By recalling (\ref{rel-1}) and  (\ref{p_ve-modified}),
we have inductively
\begin{align}\label{keylem2-eq3}
&\sum_{c \in \Omega_{x, N}(X)}p_\ve(c)
\Big(\Phi_0^{(\ve)}(x)^{-1}*\Phi_0^{(\ve)}\big(t(c)\big)\Big)_{i*}^{(1)}\nn\\
&=\sum_{c' \in \Omega_{x, N-1}(X)}p_\ve(c') 
\sum_{e \in E_{t(c')}}p_\ve(e)\Big(\Phi_0^{(\ve)}(x)^{-1}
\cdot\Phi_0^{(\ve)}\big(t(c')\big)\cdot\Phi_0^{(\ve)}\big(t(c')\big)^{-1}
\cdot\Phi_0^{(\ve)}\big(t(c)\big)\Big)_{i}^{(1)}\nn\\
&=\sum_{c' \in \Omega_{x, N-1}(X)}p_\ve(c') 
\log\Big(\Phi_0^{(\ve)}(x)^{-1}\cdot\Phi_0^{(\ve)}\big(t(c')\big)\Big)\Big|_{X_{i}^{(1)}}
+\ve \rho_{\mathbb{R}}(\gamma_p)\big|_{X_i^{(1)}}\nn\\
&=N\ve \rho_{\mathbb{R}}(\gamma_p)\big|_{X_i^{(1)}}
\qquad (x \in V, \, i=1, 2, \dots, d_1).
\end{align}

\noindent
{\bf Step~2.} Next we estimate $\mathrm{Ord}_\ve(0)$. 
Let us consider the coefficient of $X_{i_*}^{(2)}f\big(\tau_\ve\big(\Phi_0^{(\ve)}(x)\big)\big)$.
It follows from (\ref{p_ve-modified}) and (\ref{rel-1}) that
\begin{align}\label{keylem2-eq4}
&\frac{1}{N}\sum_{c \in \Omega_{x, N}(X)}p_\ve(c)
\Big\{ \Big(\Phi_0^{(\ve)}(x)^{-1}*\Phi_0^{(\ve)}\big(t(c)\big)\Big)_{i*}^{(2)}\nn\\
&\hspace{0.5cm}  - \frac{1}{2}\sum_{1 \leq \lambda < \nu \leq d_1}
\Big( \Phi_0^{(\ve)}(x)^{-1}*\Phi_0^{(\ve)}\big(t(c)\big)\Big)_{\lambda*}^{(1)}
\Big( \Phi_0^{(\ve)}(x)^{-1}*\Phi_0^{(\ve)}\big(t(c)\big)\Big)_{\nu*}^{(1)}
[\![X_\lambda^{(1)}, X_\nu^{(1)}]\!]\big|_{X_i^{(2)}} \Big\}\nn\\
&=\frac{1}{N}\sum_{c \in \Omega_{x, N}(X)}p_\ve(c) 
\Log \Big(\Phi_0^{(\ve)}(x)^{-1}*\Phi_0^{(\ve)}\big(t(c)\big)\Big)\Big|_{X_{i}^{(2)}}\nn\\
&=\frac{1}{N}\sum_{k=0}^{N-1}\sum_{c' \in \Omega_{x, k}(X)}p_\ve(c')
\sum_{e \in E_{t(c')}}p_\ve(e)\Log \big(d\Phi_0^{(\ve)}(e)\big)\big|_{X_{i}^{(2)}} \qquad (x \in V).
\end{align}
Since the function
$$
M_i^{(\ve)}(x):=\sum_{e \in E_x}p_\ve(e)\Log \big(d\Phi_0^{(\ve)}(e)\big)\big|_{X_{i}^{(2)}}
\qquad (0 \leq \ve \leq 1, \, i=1, 2, \dots, d_2, \, x \in V)
$$
satisfies $M_i^{(\ve)}(\gamma x)=M_i^{(\ve)}(x)$ for $\gamma \in \Gamma$ and 
$x \in V$ due to the $\Gamma$-invariance of $p$ and the $\Gamma$-equivariance
of $\Phi_0$, there exists a function $\mathcal{M}_i^{(\ve)} : V_0 \LA \mathbb{R}$
such that $\mathcal{M}_i^{(\ve)}\big(\pi(x)\big)=M_i^{(\ve)}(x)$ for 
$0 \leq \ve \leq 1, \, i=1, 2, \dots, d_2$ and $x \in V$. 
Moreover, we have
$$
L_{(\ve)}^k\mathcal{M}_i^{(\ve)}\big(\pi(x)\big)=L_{(\ve)}^kM_i(x) 
\qquad (k \in \mathbb{N}, \, 0 \leq \ve \leq 1, \, i=1, 2, \dots, d_2, \, x \in V)
$$
due to the $\Gamma$-invariance of $p$. 
Then, by applying the ergodic theorem 
(cf.~\cite[Theorem 3.4]{IKK}) for the transition operator $L_{(\ve)}$, 
we can choose a sufficiently small $\ve_0>0$ such that 
\begin{equation}\label{keylem2-eq3.5}
\frac{1}{N}\sum_{k=0}^{N-1}L_{(\ve)}^k\mathcal{M}_i^{(\ve)}\big(\pi(x)\big)
=\sum_{x \in V_0}m(x)\mathcal{M}_i^{(\ve)}(x)
+O_{\ve_0}\Big(\frac{1}{N}\Big). 
\end{equation}
By (\ref{keylem2-eq3.5}), we obtain that 
the left-hand side of (\ref{keylem2-eq4}) equals 
$$
\begin{aligned}
&\beta_{(\ve)}(\Phi_0^{(\ve)})\big|_{X_i^{(2)}} + O_{\ve_0}\Big(\frac{1}{N}\Big)
\qquad (0 \leq \ve \leq \ve_0, \, i=1, 2, \dots, d_2)
\end{aligned}
$$
Then Proposition \ref{key-beta} implies that 
$\mathcal{I}_1(\ve, N) \LA 0$
as $N \to \infty$ and $\ve \searrow 0$ with $N^2\ve \searrow 0$. 

We also consider the coefficient of 
$X_{i*}^{(1)}X_{j*}^{(1)}f\big( \tau_\ve \big( \Phi_0^{(\ve)}(x)\big)\big)$.  
We have
\begin{align}\label{keylem2-eq6}
&\frac{1}{2N}\sum_{c \in \Omega_{x, N}(X)}p_\ve(c)
\Big(\Phi_0^{(\ve)}(x)^{-1}*\Phi_0^{(\ve)}\big(t(c)\big)\Big)_{i*}^{(1)}
\Big(\Phi_0^{(\ve)}(x)^{-1}*\Phi_0^{(\ve)}\big(t(c)\big)\Big)_{j*}^{(1)}\nn\\
&=\frac{1}{2N}\Big\{\sum_{c' \in \Omega_{x, N-1(X)}}p_\ve(c') 
\Big(\Phi_0^{(\ve)}(x)^{-1} \cdot \Phi_0^{(\ve)}\big(t(c')\big)\Big)_{i}^{(1)}
\Big(\Phi_0^{(\ve)}(x)^{-1} \cdot \Phi_0^{(\ve)}\big(t(c')\big)\Big)_{j}^{(1)}\nn\\
&\hspace{1cm}+\sum_{e \in E_{t(c')}}p_\ve(e) 
\Log \big(d\Phi_0^{(\ve)}(e)\big)\big|_{X_i^{(1)}} 
\Log \big(d\Phi_0^{(\ve)}(e)\big)\big|_{X_j^{(1)}}\nn\\
&\hspace{1cm}+ 2(N-1)\rho_{\mathbb{R}}(\gamma_{p_\ve})\big|_{X_i^{(1)}}
\rho_{\mathbb{R}}(\gamma_{p_\ve})\big|_{X_j^{(1)}}\Big\}\nn\\
&=\frac{1}{2N}\sum_{k=0}^{N-1}\sum_{c' \in \Omega_{x, k(X)}}
p_\ve(c') \sum_{e \in E_{t(c')}}p_\ve(e)
\Log \big(d\Phi_0^{(\ve)}(e)\big)\big|_{X_i^{(1)}} 
\Log \big(d\Phi_0^{(\ve)}(e)\big)\big|_{X_j^{(1)}}\nn \\
&\hspace{1cm}+ \frac{1}{2}(N-1)\ve^2
\rho_{\mathbb{R}}(\gamma_{p})\big|_{X_i^{(1)}}
\rho_{\mathbb{R}}(\gamma_{p})\big|_{X_j^{(1)}}\nn \\
&=\frac{1}{2N}\sum_{k=0}^{N-1}L_{(\ve)}^k 
N^{(\ve)}_{ij}(x) 
+ \frac{1}{2}(N-1)\ve^2\rho_{\mathbb{R}}(\gamma_{p})\big|_{X_i^{(1)}}
\rho_{\mathbb{R}}(\gamma_{p})\big|_{X_j^{(1)}} \quad (x \in V)
\end{align}
by using  (\ref{p_ve-modified}) and (\ref{rel-2}), where 
the function 
$N_{ij}^{(\ve)} : V \LA \mathbb{R}$ is defined by
$$
N_{ij}^{(\ve)}(x):=\sum_{e \in E_x}p_\ve(e)\Log \big(d\Phi_0^{(\ve)}(e)\big)\big|_{X_i^{(1)}} 
\Log \big(d\Phi_0^{(\ve)}(e)\big)\big|_{X_j^{(1)}}.
$$
for $0 \leq \ve \leq 1, \, i, j =1, 2, \dots, d_1$
and $x \in V$. In the same argument as above, $N_{ij}^{(\ve)}$ is $\Gamma$-invariant
and there exists a function $\mathcal{N}_{ij}^{(\ve)} : V_0 \LA \mathbb{R}$
such that $\mathcal{N}_{ij}^{(\ve)}\big(\pi(x)\big)=N_{ij}^{(\ve)}(x)$ for $x \in V$. 
We also have
$$
L_{(\ve)}^k\mathcal{N}_{ij}^{(\ve)}\big(\pi(x)\big)=L_{(\ve)}^kN_{ij}^{(\ve)}(x) 
\qquad (k \in \mathbb{N}, \, 0 \leq \ve \leq 1, \, i, j=1, 2, \dots, d_2, \, x \in V)
$$
by the $\Gamma$-invariance of $p$. 
Thus, we choose a sufficiently small $\ve_0'>0$ such that
\begin{align}\label{keylem2-eq7}
\frac{1}{2N}\sum_{k=0}^{N-1}L_{(\ve)}^k N_{ij}^{(\ve)}(x)
&=\frac{1}{2N}\sum_{k=0}^{N-1}L_{(\ve)}^k \mathcal{N}_{ij}^{(\ve)}\big(\pi(x)\big)\nn\\
&=\frac{1}{2}\sum_{x \in V_0}m(x)\big(\mathcal{N}(\Phi_0^{(\ve)})_{ij}\big)(x)
+O_{\ve_0'}\Big(\frac{1}{N}\Big)\nn\\
&=\frac{1}{2}\sum_{e \in E_0}\widetilde{m}_\ve(e)
\Log \big(d\Phi_0^{(\ve)}(e)\big)\big|_{X_i^{(1)}} 
\Log \big(d\Phi_0^{(\ve)}(e)\big)\big|_{X_j^{(1)}}\nn\\
&\hspace{0.5cm}+O_{\ve_0'}\Big(\frac{1}{N}\Big) \qquad (0 \leq \ve \leq \ve_0', \, i, j=1, 2, \dots, d_1) 
\end{align}
by the ergodic theorem. 
Recall that $\{V_1, V_2, \dots, V_{d_1}\}$ denotes the orthonormal basis 
in $(\g^{(1)}, g_0^{(0)})$. 
In particular,  put $X_i^{(1)}=V_i$ for $i=1, 2, \dots, d_1$ and 
let $\{\omega_1, \omega_2, \dots, \omega_{d_1}\}$ 
be the dual basis of $\{V_1, V_2, \dots, V_{d_1}\}$. 
Then we have
\begin{align}\label{keylem2-eq8}
&\frac{1}{2}\sum_{e \in E_0}\widetilde{m}_\ve(e)
\Log \big(d\Phi_0^{(\ve)}(\widetilde{e})\big)\big|_{V_i} 
\Log \big(d\Phi_0^{(\ve)}(\widetilde{e})\big)\big|_{V_j}\nn\\
&=\frac{1}{2}\Big( \sum_{e \in E_0}\widetilde{m}_\ve(e)\omega_i^{(\ve)}(e)\omega_j^{(\ve)}(e) 
- \la \gamma_{p_\ve}, \omega_i \ra \la \gamma_{p_\ve}, \omega_j \ra \Big) 
+ \frac{1}{2}\ve^2 \la \gamma_p, \omega_i \ra \la \gamma_p, \omega_j \ra\nn\\
&= \frac{1}{2}\La \omega_i^{(\ve)}, \omega_j^{(\ve)} \Ra_{(\ve)} + \frac{1}{2}\ve^2 
\rho_{\mathbb{R}}(\gamma_p)\big|_{V_i}
\rho_{\mathbb{R}}(\gamma_p)\big|_{V_j}
\qquad (i, j=1, 2, \dots, d_1).
\end{align}
The coefficient of $X_{i*}^{(1)}X_{j*}^{(1)}f\big( \tau_\ve \big( \Phi_0^{(\ve)}(x)\big)\big)$ equals
\begin{equation}\label{keylem2-eq9}
-\Big(\frac{1}{2}\La \omega_i^{(\ve)}, \omega_j^{(\ve)} \Ra_{(\ve)} 
+ \frac{1}{2}N\ve^2 \rho_{\mathbb{R}}(\gamma_p)\big|_{V_i}
\rho_{\mathbb{R}}(\gamma_p)\big|_{V_j}\Big)+O_{\ve_0'}\Big(\frac{1}{N}\Big) \qquad (i, j=1, 2, \dots, d_1)
\end{equation}
by combining (\ref{keylem2-eq6}) with (\ref{keylem2-eq7}) and (\ref{keylem2-eq8}).
Therefore, (\ref{keylem2-eq9}) and the continuity of $\La \cdot , \cdot \Ra_{(\ve)}$
as $\ve \searrow 0$ (cf.~\cite[Lemma 5.2]{IKK}) imply
\begin{equation}\label{keylem2-eq99}
\mathrm{Ord}_\ve(0) =\mathcal{I}_1(\ve, N)+\mathcal{I}_2(\ve, N)= -\frac{1}{2}\sum_{i=1}^{d_1} 
V_{i*}^2 f\Big( \tau_\ve \big( \Phi_0^{(\ve)}(x)\big)\Big)+O(N^2\ve)+O_{\ve_0}\Big(\frac{1}{N}\Big)
\end{equation}
as $N \to \infty$ and $\ve \searrow 0$ with $N^2\ve \searrow 0$. 

We finally discuss the estimate of $\sum_{k \geq 1} \mathrm{Ord}_\ve(k)$. 
At the beginning, we show that the coefficient of $X_{i*}^{(k)}f\big( \tau_\ve \big( \Phi_0^{(\ve)}(x)\big)\big)$
vanishes as $N \to \infty$ and $\ve \searrow 0$ with $N^2\ve \searrow 0$. 
Thanks to
$$
\Big|\Big(\Phi_0^{(\ve)}(x)^{-1} \cdot \Phi_0^{(\ve)}\big(t(c)\big)\Big)_{i}^{(k)}\Big|
\leq CN^k \qquad (0 \leq \ve \leq 1, \, x \in V),
$$
(\ref{p_ve-modified}) and (\ref{rel-5}), we have
$$
\begin{aligned}
&\frac{\ve^{k-2}}{N}\sum_{x \in \Omega_{x, N}(X)}p_\ve(c)
\Big(\Phi_0^{(\ve)}(x)^{-1}*\Phi_0^{(\ve)}\big(t(c)\big)\Big)_{i*}^{(k)} \nn\\
&=\frac{\ve^{k-2}}{N}\sum_{x \in \Omega_{x, N}(X)}p_\ve(c)\Big\{
\Big(\Phi_0^{(\ve)}(x)^{-1} \cdot \Phi_0^{(\ve)}\big(t(c)\big)\Big)_{i}^{(k)}\nn\\
&\hspace{1cm}+\sum_{\substack{|K_1|+|K_2| \leq k-1 \\ |K_2|>0}} C_{K_1, K_2}
\mathcal{P}_*^{K_1}\Big( \Phi_0^{(\ve)}(x)^{-1}\Big) \mathcal{P}^{K_2}
\Big(\Phi_0^{(\ve)}(x)^{-1}\cdot \Phi_0^{(\ve)}\big(t(c)\big)\Big)\Big\}\nn\\
& \leq CM_i^{(k)}\Big( \tau_\ve\big(\Phi_0^{(\ve)}(x)\big)\Big)
\Big( \ve^{k-2}N^{k-1} +\sum_{|K_1| \leq k-2}\ve^{k-1-|K_1|}
+\sum_{\substack{|K_1|+|K_2| \leq k-1 \\ |K_2|\geq 2}}
\ve^{k-2-|K_1|}N^{|K_2|-1}\Big)
\end{aligned}
$$
for $i=1, 2, \dots, d_k$ and some continuous function $M_i^{(k)} : G \LA \mathbb{R}$. 
This converges to zero as $N \to \infty$ and $\ve \searrow 0$ with $N^2\ve \searrow 0$. 
We also observe that the 
coefficient of $X_{i_1*}^{(k_1)}X_{i_2*}^{(k_2)}f\big( \tau_\ve \big( \Phi_0^{(\ve)}(x)\big)\big)$
converges to zero as $N \to \infty$ and $\ve \searrow 0$ with $N^2\ve \searrow 0$
by following the same argument as above.

We also consider the coefficient of 
$(\del^3f/\del g_{i_1*}^{(k_1)}\del g_{i_2*}^{(k_2)}\del g_{i_3*}^{(k_3)})(\theta)$. 
Since $f$ is compactly supported, it is sufficient to show by induction on $k=1, 2, \dots, r$ that,
if $\ve N<1$, then
\begin{equation}\label{keylem2-eq11}
\ve^k\Big(\Phi_0^{(\ve)}(x)^{-1}*\Phi_0^{(\ve)}\big(t(c)\big)\Big)_{i*}^{(k)}
\leq M_i^{(k)}\Big( \tau_\ve\big(\Phi_0^{(\ve)}(x) * \theta\big)\Big) \times \ve N
\end{equation}
for $i=1, 2, \dots, d_k$ and some continuous function $M_i^{(k)} : G \LA \mathbb{R}$. 
The cases $k=1$ and $k=2$ are clear. 
Suppose that (\ref{keylem2-eq11}) holds for less than $k$. 
We have
$$
\begin{aligned}
\ve^k\Big(\Phi_0^{(\ve)}(x)^{-1}*\Phi_0^{(\ve)}\big(t(c)\big)\Big)_{i*}^{(k)}
&=\ve^k\Big\{\Big(\Phi_0^{(\ve)}(x)^{-1} \cdot \Phi_0^{(\ve)}\big(t(c)\big)\Big)_{i}^{(k)}+\sum_{\substack{|K_1|+|K_2| \leq k-1 \\ |K_2|>0}} C_{K_1, K_2}\nn\\
&\hspace{1cm}\times
\mathcal{P}_*^{K_1}\Big( \Phi_0^{(\ve)}(x)^{-1}\Big) \mathcal{P}^{K_2}
\Big(\Phi_0^{(\ve)}(x)^{-1}\cdot \Phi_0^{(\ve)}\big(t(c)\big)\Big)\Big\}.
\end{aligned}
$$
by using (\ref{CBH-formula}) and (\ref{rel-5}).
Then we see that
$$
\begin{aligned}
\Big( \Phi_0^{(\ve)}(x)^{-1}\Big)_{i_1*}^{(k_1)} &= 
\Big( \theta * \big(\Phi_0^{(\ve)}(x)*\theta\big)^{-1}\Big)_{i_1*}^{(k_1)}\nn\\
&=\theta_{i_1*}^{(k_1)} + \Big( \big(\Phi_0^{(\ve)}(x)*\theta\big)^{-1}\Big)_{i_1*}^{(k_1)}\nn\\
&\hspace{1cm}+\sum_{\substack{|L_1|+|L_2|=k_1 \\ |L_1|, |L_2|>0}}C_{L_1, L_2}
\mathcal{P}_*^{L_1}(\theta)\mathcal{P}_*^{L_2}\Big( \big(\Phi_0^{(\ve)}(x)*\theta\big)^{-1}\Big).
\end{aligned}
$$
Thus, we have inductively 
$$
\Big|\Big( \Phi_0^{(\ve)}(x)^{-1}\Big)_{i_1*}^{(k_1)} \Big| 
\leq M\Big( \Phi_0^{(\ve)}(x)*\theta\Big)
$$
for a continuous function $M : G \LA \mathbb{R}$ and $k_1 \leq k-1$.
We then conclude
$$
\begin{aligned}
&\ve^k\Big(\Phi_0^{(\ve)}(x)^{-1}*\Phi_0^{(\ve)}\big(t(c)\big)\Big)_{i*}^{(k)}\nn\\
&\leq C\Big( \ve^k N^k+\sum_{\substack{|K_1|+|K_2| \leq k-1 \\ |K_2|>0}}
M\Big( \tau_\ve\big(\Phi_0^{(\ve)}(x)*\theta\big)\Big)\ve^{k-|K_1|}N^{|K_2|}\Big)\\
&\leq M_i^{(k)}\Big( \tau_\ve\big(\Phi_0^{(\ve)}(x)*\theta\big)\Big) \times \ve N.
\end{aligned}
$$
for some continuous function $M_i^{(k)} : G \LA \mathbb{R}$. 
These estimates  implies that $\sum_{k \geq 1} \mathrm{Ord}_\ve(k)$ converges
to zero as $N \to \infty$ and $\ve \searrow 0$ with $N^2\ve \searrow 0$.

Consequently, we obtain
$$
\Big\|\frac{1}{N\ve^2}\big(I-L_{(\ve)}^N\big)P^H_\ve f(x)
-P^H_\ve \A f(x)\Big\|_\infty^X \LA 0
$$
as $N \to \infty$ and $\ve \searrow 0$ with $N^2\ve \searrow 0$
by combining (\ref{keylem2-eq1}) with (\ref{keylem2-eq3})
and (\ref{keylem2-eq99}).  This completes the proof. \qed

\vspace{2mm}

\noindent
{\bf Proof of Theorem \ref{CLT2}}. 
We partially follow the argument by Kotani \cite[Theorem 4]{Kotani}.
Let $N=N(n)$ be the integer satisfying $n^{1/5} \leq N<n^{1/5}+1$
and let $k_N$ and $r_N$ be the quotient and the remainder of $([nt]-[ns])/N(n)$, respectively.
We put $\ve_N:=n^{-1/2}$ and $h_N:=N\ve_N^2$. 
Then we have
$k_N h_N=\big([nt]-[ns]-r_N\big)\ve_N^2 \to t-s \quad (n \to \infty).$

Since $C_0^\infty(G_{(0)}) \subset \Dom(\A) \subset C_\infty(G_{(0)})$ and 
$C_0^\infty(G_{(0)})$ is dense in $C_\infty(G_{(0)})$, the linear operator $\A$ is 
densely defined in $C_\infty(G_{(0)})$. 
Furthermore, $(\lambda - \A)\big(C_0^\infty(G_{(0)})\big)$ is dense in $C_\infty(G_{(0)})$ 
for some $\lambda>0$
(cf.~Robinson \cite[p.304]{Rob}). 
Hence, by combining Lemma \ref{key-lem2} and Trotter's approximation theorem (cf.~\cite{Trotter}),  
we obtain
\begin{equation}\label{eq7}
\lim_{n \to \infty}\Big\| L_{(n^{-1/2})}^{Nk_N}P^H_{n^{-1/2}}f 
- P^H_{n^{-1/2}}\e^{-(t-s)\A}f\Big\|_\infty^X=0 \qquad \big(  f \in C_0^\infty(G_{(0)})\big).
\end{equation}
On the other hand, Lemma \ref{key-lem2} implies
\begin{equation}\label{eq8}
\lim_{n \to \infty}\Big\| \frac{1}{r_N\ve_N^2}\big(I-L_{(n^{-1/2})}^{r_N}\big)P^H_{n^{-1/2}}f 
- P^H_{n^{-1/2}}\A f \Big\|_\infty^X=0 \qquad \big(  f \in C_0^\infty(G_{(0)})\big).
\end{equation}
Here we have
\begin{align}\label{eq9}
&\Big\| L_{(n^{-1/2})}^{[nt]-[ns]}P^H_{n^{-1/2}}f - P^H_{n^{-1/2}}\e^{-(t-s)\A}f\Big\|_\infty^X \nn\\
&\leq \Big\| \big(I-L_{(n^{-1/2})}^{r_N}\big)P^H_{n^{-1/2}}f\Big\|_\infty^X + 
\Big\| L_{(n^{-1/2})}^{Nk_N}P^H_{n^{-1/2}}f - P^H_{n^{-1/2}}\e^{-(t-s)\A}f\Big\|_\infty^X.
\end{align}
It follows from $\|P^H_{n^{-1/2}}\| \leq 1$ that
\begin{align}\label{eq10}
&\Big\| \big(I - L_{(n^{-1/2})}^{r_N}\big)P^H_{n^{-1/2}}f\Big\|_\infty^X \nn\\
&\leq r_N\ve_N^2 \Big\|\frac{1}{r_N\ve_N^2}\big(I-L_{(n^{-1/2})}^{r_N}\big)P^H_{n^{-1/2}}f
-P^H_{n^{-1/2}}\A f\Big\|_\infty^X + r_N\ve_N^2\big\| P^H_{n^{-1/2}}\A f\big\|_\infty^X \nn \\
&\leq r_N\ve_N^2 \Big\|\frac{1}{r_N\ve_N^2}\big(I-L_{(n^{-1/2})}^{r_N}\big)P^H_{n^{-1/2}}f
-P^H_{n^{-1/2}}\A f\Big\|_\infty^X + r_N\ve_N^2\big\| \A f\big\|_\infty^G.
\end{align}
Then, we obtain   
\begin{equation}\label{semigroup CLT2 harmonic}
\lim_{n \to \infty}\Big\| L_{(n^{-1/2})}^{[nt]-[ns]}P^H_{n^{-1/2}} f 
- P^H_{n^{-1/2}} \e^{-(t-s)\A}f\Big\|_\infty^X=0
\end{equation}
for $f \in C_0^\infty(G_{(0)})$
by combining (\ref{eq7}), (\ref{eq8}), (\ref{eq9}) and (\ref{eq10}) 
with $r_N\ve_N^2 \to 0$ as $n \to \infty$.
For $f \in C_\infty(G_{(0)})$, we also obtain (\ref{semigroup CLT2 harmonic})
by following \cite[Theorem 2.1]{IKK}.

We are now ready to show (\ref{semigroup CLT2}). By the triangular inequality and $\|L\| \leq 1$, we have
\begin{align}
&\Big\| L_{(n^{-1/2})}^{[nt]-[ns]}P_{n^{-1/2}} f 
- P_{n^{-1/2}} \e^{-(t-s)\A}f\Big\|_\infty^X\nn\\
&\leq \Big\| P_{n^{-1/2}} f 
- P^H_{n^{-1/2}} f\Big\|_\infty^X
+\Big\| L_{(n^{-1/2})}^{[nt]-[ns]}P^H_{n^{-1/2}} f 
- P^H_{n^{-1/2}} \e^{-(t-s)\A}f\Big\|_\infty^X\nn\\
&\hspace{1cm} +\Big\| P^H_{n^{-1/2}} \e^{-(t-s)\A}f
 - P_{n^{-1/2}} \e^{-(t-s)\A}f\Big\|_\infty^X \qquad \big(f \in C_\infty(G_{(0)})\big) \label{Pre}
\end{align}
We note that the functions $f$ and $\e^{-(t-s)\A}f$ is uniformly continuous and 
$$
d_{\mathrm{CC}}\big( \tau_{n^{-1/2}}\Phi^{(n^{-1/2})}(x), 
\tau_{n^{-1/2}}\Phi_0^{(n^{-1/2})}(x)\big)=\frac{1}{\sqrt{n}}
d_{\mathrm{CC}}\big( \Phi^{(n^{-1/2})}(x), 
\Phi_0^{(n^{-1/2})}(x)\big) \leq \frac{M}{\sqrt{n}} 
$$
for some $M>0$. Therefore, by using (\ref{semigroup CLT2 harmonic}) 
and by letting $n \to \infty$ in (\ref{Pre}), 
we obtain the desired convergence (\ref{semigroup CLT2}) for $f \in C_\infty(G_{(0)})$. 
The latter part of Theorem \ref{CLT2} is obtained 
in the same way as \cite[Theorem 2.1]{IKK}. 
This completes the proof. \qed

\subsection{Proof of functional CLT under modified harmonicities}
We set 
$\xi_n^{(\ve)}(c):=\Phi_0^{(\ve)}\big(w_n(c)\big)$ for 
$0 \leq \ve \leq 1, \, n=0, 1, 2, \dots$ and $c \in \Omega_{x_*}(X)$.
Then $G_{(0)}$-valued random walk 
$(\Omega_{x_*}(X), \mathbb{P}_{x_*}^{(\ve)}, \{\xi_n^{(\ve)}\}_{n=0}^\infty)$
associated with the transition probability $p_\ve$ is induced. 
For $t \geq 0, \, n=1, 2, \dots$ and $0 \leq \ve \leq 1$, 
we set
$$
\mathcal{X}_t^{(\ve, n)}(c)
:=\tau_{n^{-1/2}}\big( \xi_{[nt]}^{(\ve)}(c)\big) \qquad \big(c \in \Omega_{x_*}(X)\big).
$$
and define 
$$
\mathcal{Y}_{t_k}^{(\ve, n)}(c)
:=\tau_{n^{-1/2}}\big( \xi_{nt_k}^{(\ve)}(c)\big)
=\tau_{n^{-1/2}}\big( \Phi_0^{(\ve)}(w_{k}(c))\big)
 \qquad \big(t_k \in \mathcal{D}_n, \, c \in \Omega_{x_*}(X)\big).
$$
We also define the $G_{(0)}$-valued continuous stochastic process 
$(\mathcal{Y}_t^{(\ve, n)})_{0 \leq t \leq 1} $
given by the $d_{\mathrm{CC}}$-geodesic interpolation of 
$\{ \mathcal{Y}_{t_k}^{(\ve, n)}\}_{k=0}^n$. 
We denote by $\Prob^{(\ve, n)}$ 
the probability measure on $C^{0, \alpha\hol}([0, 1]; G_{(0)})$
induced by the stochastic process $\mathcal{Y}_\cdot^{(\ve, n)}$
for $0 \leq \ve \leq 1$ and $n \in \mathbb{N}$.

The aim of this subsection is to prove the following. 

\begin{tm}\label{FCLT2} 
We assume {\bf (A1)} and {\bf (A2)}. Then the sequence
$\{\mathcal{Y}^{(n^{-1/2}, n)}\}_{n= 1}^\infty $
converges in law to the diffusion process
$(Y_t)_{0 \leq t \leq 1}$ in 
$C^{0, \alpha\text{\normalfont{-H\"ol}}}([0, 1]; G_{(0)})$
as $n \to \infty$ for all $\alpha<1/2$. 
\end{tm}

We put 
$$
\|d\Phi_0^{(\ve)}\|_\infty=\max_{e \in E_0}\max_{k=1, 2, \dots, r}
\big\|\log\big(d\Phi_0^{(\ve)}(\widetilde{e})\big)\big|_{\g^{(k)}}\big\|_{\g^{(k)}}^{1/k}
 \qquad ( 0 \leq \ve \leq 1).
 $$
We describe a relation between $\|d\Phi_0^{(\ve)}\|_\infty$ and 
$\|d\Phi_0^{(0)}\|_\infty$ for every $0 \leq \ve \leq 1$.
Thanks to \cite[Lemma 5.3 (3)]{IKK}, {\bf (A2)} and (\ref{identity-dphi}),
we easily obtain the following:
\begin{lm}\label{dphi-relation}
Under {\bf (A2)}, there exists a positive constant $C$ such that 
$$
\sup_{0 \leq \ve \leq 1}\|d\Phi_0^{(\ve)}\|_\infty \leq C\|d\Phi_0^{(0)}\|_\infty.
$$
\end{lm}

Denote by $G_{(0)}^{(k)}$ the nilpotent Lie group of step $k$ whose Lie algebra coincides with
$(\g^{(1)}, g_0^{(0)}) \oplus \g^{(2)} \oplus \cdots \oplus \g^{(k)}$. 
For the piecewise smooth stochastic process $(\mathcal{Y}_t^{(\ve, n)})_{0 \leq t \leq 1}$, 
we define its truncated process by 
$$
\mathcal{Y}_t^{(\ve, n; \,k)}=
\big( \mathcal{Y}_t^{(\ve, n), 1}, 
 \mathcal{Y}_t^{(\ve, n), 2}, \dots,
  \mathcal{Y}_t^{(\ve, n), k} \big) \in G_{(0)}^{(k)}
  \qquad (k=1, 2, \dots, r)
$$
in the $(*)$-coordinate system.
By Lemma \ref{dphi-relation}, we put 
$$
\sup_{0 \leq \ve \leq 1}\Big\{\|d\Phi_0^{(\ve)}\|_\infty
+\| \rho_{\mathbb{R}}(\gamma_p)\|_{\g^{(1)}}\Big\} \leq C\|d\Phi_0^{(0)}\|_\infty
+\| \rho_{\mathbb{R}}(\gamma_p)\|_{\g^{(1)}}=:M.
$$

As is well-known in probability theory, 
it suffices to show the tightness of $\{\Prob^{(n^{-1/2}, n)}\}_{n=1}^\infty$ and
the convergence of the finite dimensional distribution 
of $\{\mathcal{Y}^{(n^{-1/2}, n)}\}_{n=1}^\infty$ to obtain Theorem \ref{FCLT2}.  
In the former part of this section, we  show the following.

\begin{lm}\label{tightness2}
$\{\Prob^{(n^{-1/2}, n)}\}_{n=1}^\infty$ is tight in
$C^{0, \alpha\text{\normalfont{-H\"ol}}}([0, 1]; G_{(0)}),$
where $\alpha<1/2$. 
\end{lm}

As the first step of the proof of Lemma \ref{tightness2}, 
we need to show the following lemma. 

\begin{lm}\label{sublemma2}
Let $m, n$ be positive integers. Then there exists a constant $C>0$ independent of $n$ 
{\rm(}however, it may depend on $m${\rm)} such that 
\begin{equation}\label{tight2-2}
\mathbb{E}^{\mathbb{P}_{x_*}^{(n^{-1/2})}}
\Big[ d_{\mathrm{CC}}(\mathcal{Y}_{s}^{(n^{-1/2}, n; \, 2)}, 
\mathcal{Y}_{t}^{(n^{-1/2}, n; \, 2)})^{4m}\Big]
\leq C(t-s)^{2m}  \qquad (0 \leq s \leq  t \leq 1).
\end{equation}
\end{lm}
\noindent
{\bf Proof.} 
We split the proof into several steps, based on the argument in \cite{IKN}.

\vspace{2mm}
\noindent
{\bf Step\,1.} 
First we show 
\begin{align}\label{tight2-1}
&\hspace{2cm}\mathbb{E}^{\mathbb{P}_{x_*}^{(n^{-1/2})}}
\Big[ d_{\mathrm{CC}}
\big(\mathcal{Y}_{t_k}^{(n^{-1/2}, n; \, 2)}, 
\mathcal{Y}_{t_\ell}^{(n^{-1/2}, n; \, 2)}\big)^{4m}\Big]
\leq C\Big(\frac{\ell-k}{n}\Big)^{2m}  \nn\\
&\hspace{8cm}\big(n, m \in \mathbb{N}, \, t_k, t_\ell \in \mathcal{D}_n \, (k \leq \ell)\big)
\end{align}
for some $C>0$ which is independent of $n$ (but depending on $m$).
By noting the equivalence mentioned in Section 3.1,  
(\ref{tight2-1}) is equivalent to the existence 
of positive constants $C^{(1)}$ and $C^{(2)}$ independent of $n$ such that 
\begin{equation}\label{expectation2-1}
\mathbb{E}^{\mathbb{P}_{x_*}^{(n^{-1/2})}}
\Big[ \big\| \Log\big((\mathcal{Y}_{t_k}^{(n^{-1/2}, n)})^{-1} \cdot
\mathcal{Y}_{t_{\ell}}^{(n^{-1/2}, n)}\big)
 \big|_{\g^{(1)}}\big\|_{\g^{(1)}}^{4m}\Big] 
 \leq C^{(1)}\Big(\frac{\ell-k}{n}\Big)^{2m},
\end{equation}
\begin{equation}\label{expectation2-2}
\mathbb{E}^{\mathbb{P}_{x_*}^{(n^{-1/2})}}
\Big[ \big\| \Log\big((\mathcal{Y}_{t_k}^{(n^{-1/2}, n)})^{-1}
 \cdot \mathcal{Y}_{t_{\ell}}^{(n^{-1/2}, n)}\big)\big|_{\g^{(2)}}\big\|_{\g^{(2)}}^{2m}\Big]
  \leq C^{(2)}\Big(\frac{\ell-k}{n}\Big)^{2m}.
\end{equation}

\vspace{2mm}
\noindent
{\bf Step\,2.} We here prove (\ref{expectation2-1}). 
We have
\begin{align}\label{step2-2-1}
 &\mathbb{E}^{\mathbb{P}_{x_*}^{(\ve)}}
\Big[ \big\| \Log\big((\mathcal{Y}_{t_k}^{(\ve, n)})^{-1} \cdot 
 \mathcal{Y}_{t_{\ell}}^{(\ve, n)}\big)
 \big|_{\g^{(1)}}\big\|_{\g^{(1)}}^{4m}\Big] \nn\\
&=  \Big( \frac{1}{\sqrt{n}}\Big)^{4m}  \mathbb{E}^{\mathbb{P}_{x_*}^{(\ve)}}
\Big[ \Big( \sum_{i=1}^{d_1} \Log \big((\xi_k^{(\ve)})^{-1} \cdot \xi_\ell^{(\ve)} \big)
 \big|_{X_i^{(1)}}^2\Big)^{2m}\Big]\nn\\
&\leq \Big( \frac{1}{\sqrt{n}}\Big)^{4m} \cdot d_1^{2m} 
\max_{i=1, 2, \dots, d_1} \max_{x \in \mathcal{F}} 
\Big\{\sum_{c \in \Omega_{x, \ell-k}(X)}p_{\ve}(c)\nn\\
&\hspace{1cm}\times
\Log\Big( \Phi_0^{(\ve)}(x)^{-1} \cdot \Phi_0^{(\ve)}\big(t(c)\big)\Big)
\Big|_{X_i^{(1)}}^{4m}\Big\} \qquad (0 \leq \ve \leq 1),
\end{align}
where $\mathcal{F}$ stands for the fundamental domain in $X$ 
containing the reference point $x_* \in V$. 
For $i=1, 2, \dots, d_1, \, x \in \mathcal{F}, \, N \in \mathbb{N}$, $0 \leq \ve \leq 1$
 and $c=(e_1, e_2, \dots, e_N) \in \Omega_{x_*, N}(X)$,  put 
$$
\begin{aligned}
\mathcal{J}_i^{(\ve)}(j)&:=\Log \big( d\Phi_0^{(\ve)}(e_j)\big)\big|_{X_i^{(1)}} 
- \ve\rho_{\mathbb{R}}(\gamma_p) \big|_{X_i^{(1)}},\\
\mathcal{N}_N^{(i, x)}(\Phi_0^{(\ve)}; c)
&:=\Log\Big(\Phi_0^{(\ve)}(x)^{-1} 
\cdot \Phi_0^{(\ve)}\big(t(c)\big)\Big) \Big|_{X_i^{(1)}} 
- N \rho_{\mathbb{R}}(\gamma_{p_{\ve}}) \big|_{X_i^{(1)}}
=\sum_{j=1}^N \mathcal{J}_i^{(\ve)}(j). 
\end{aligned}
$$
Note that 
$|\mathcal{J}_i^{(\ve)}(j)| \leq M$ for 
$0 \leq \ve \leq 1, \, i=1, 2, \dots, d_1$ and $j=1, 2, \dots, N$.
Then we see that 
$\{\mathcal{N}_N^{(i, x)}\}_{N=1}^\infty$ 
is a martingale  
for $i=1, 2, \dots, d_1$ and $x \in \mathcal{F}$
(cf.~\cite[Lemma 3.3]{IKN}).
Hence, we use the Burkholder--Davis--Gundy inequality
with the exponent $4m$ to obtain
\begin{align}\label{martingale2-1}
&\sum_{c \in \Omega_{x, N}(X)}
p_{\ve}(c)\Log\Big( \Phi_0^{(\ve)}(x)^{-1} \cdot \Phi_0^{(\ve)}\big(t(c)\big)\Big)
\Big|_{X_i^{(1)}}^{4m}\nn\\
&\leq 2^{4m-1}\sum_{c \in \Omega_{x, N}(X)}
p_{\ve}(c) \Big\{ \big(\mathcal{N}_N^{(i, x)}(c)\big)^{4m} 
+ \Big(N\ve\rho_{\mathbb{R}}(\gamma_p)\big|_{X_i^{(1)}}\Big)^{4m}\Big\}\nn\\
&\leq2^{4m-1}\mathcal{C}_{(4m)}^{4m}\sum_{c \in \Omega_{x, N}(X)}p_{\ve}(c) 
\Big\{\sum_{j=1}^{N} \mathcal{J}_i^{(\ve)}(j)^2\Big\}^{2m}+2^{4m-1}\ve^{4m}N^{4m}
\big\|\rho_{\mathbb{R}}(\gamma_p)\big\|_{\g^{(1)}}^{4m}\nn\\
&\leq  2^{4m}\mathcal{C}_{(4m)}^{4m} 
M^{2m}N^{2m}
+2^{4m-1}M^{4m}\ve^{4m}N^{4m}\nn\\
&\hspace{6.5cm}(x \in \mathcal{F},\, i=1, 2, \dots, d_1, \, 0 \leq \ve \leq 1, \, N \in \mathbb{N}).
\end{align}
In particular, (\ref{martingale2-1}) implies
\begin{align}\label{step2-2-2}
&\sum_{c \in \Omega_{x, \ell-k}(X)}p_{n^{-1/2}}(c)
\Log\Big( \Phi_0^{(n^{-1/2})}(x)^{-1} 
\cdot \Phi_0^{(n^{-1/2})}\big(t(c)\big)\Big)\Big|_{X_i^{(1)}}^{4m}\nn\\
&\leq \Big\{2^{4m}\mathcal{C}_{(4m)}^{4m} M^{2m}
+2^{4m-1}M^{4m}\Big\}(\ell-k)^{2m}
\end{align}
by putting $\ve=n^{-1/2}$ and $N=\ell-k$, 
where we note $(\ell-k)/n<1$ since $1 \leq k \leq \ell \leq n$. 
We then obtain
$$
\begin{aligned}
&\mathbb{E}^{\mathbb{P}_{x_*}^{(n^{-1/2})}}
\Big[ \big\| \Log\big((\widetilde{\mathcal{Y}}_{t_k}^{(n^{-1/2}, n)})^{-1} \cdot
 \widetilde{\mathcal{Y}}_{t_{\ell}}^{(n^{-1/2}, n)}\big)
 \big|_{\g^{(1)}}\big\|_{\g^{(1)}}^{4m}\Big]\\
 &\leq d_1^{2m}\Big\{2^{4m}\mathcal{C}_{(4m)}^{4m} M^{2m}
+2^{4m-1}M^{4m}\Big\}\Big(\frac{\ell-k}{n}\Big)^{2m}
=:C^{(1)}\Big(\frac{\ell-k}{n}\Big)^{2m}
\end{aligned}
$$
by combining (\ref{step2-2-1}) with (\ref{step2-2-2}), 
which leads to (\ref{expectation2-1}).

\vspace{2mm}
\noindent
{\bf Step\,3.} We show (\ref{expectation2-2}) at this step.  
We also see
\begin{align}\label{step2-3-1}
&\mathbb{E}^{\mathbb{P}_{x_*}^{(\ve)}}
\Big[ \big\| \Log\big((\mathcal{Y}_{t_k}^{(\ve, n)})^{-1} \cdot
\mathcal{Y}_{t_{\ell}}^{(\ve, n)}\big)
 \big|_{\g^{(2)}}\big\|_{\g^{(2)}}^{2m}\Big]\nn\\
&\leq \Big(\frac{1}{n}\Big)^{2m}\cdot d_2^{2m}
\max_{i=1, 2, \dots, d_2} \max_{x \in \mathcal{F}}
\Big\{\sum_{c \in \Omega_{x, \ell-k}(X)} p_{\ve}(c)\nn\\
&\hspace{1cm}\times 
\Log\Big( \Phi_0^{(\ve)}(x)^{-1} \cdot 
\Phi_0^{(\ve)}\big(t(c)\big)\Big)\Big|_{X_{i}^{(2)}}^{2m}\Big\} \qquad (0 \leq \ve \leq 1).
\end{align}
in the similar way to (\ref{step2-2-1}). 
Then we have
\begin{align}\label{step2-3-2}
&\Log\Big( \Phi_0^{(\ve)}(x)^{-1} 
\cdot \Phi_0^{(\ve)}\big(t(c)\big)\Big)\Big|_{X_i^{(2)}}^{2m}\nn\\
&=\Big(\sum_{j=1}^{\ell-k}\Log \big( d\Phi_0^{(\ve)}(e_j) \big)\big|_{X_i^{(2)}}
-\frac{1}{2}\sum_{1 \leq j_1 < j_2 \leq \ell-k}\sum_{1 \leq \lambda < \nu \leq d_1}
[\![X_\lambda^{(1)}, X_\nu^{(1)}]\!]\big|_{X_i^{(2)}} \nn\\
&\hspace{0.8cm}\times\Big\{ 
           \Log \big( d\Phi_0^{(\ve)}(e_{j_1})\big)\big|_{X_\lambda^{(1)}} 
           \Log \big( d\Phi_0^{(\ve)}(e_{j_2})\big)\big|_{X_\nu^{(1)}}\nn\\
&\hspace{2cm}
-  \Log \big( d\Phi_0^{(\ve)}(e_{j_1})\big)\big|_{X_\nu^{(1)}} 
\Log \big( d\Phi_0^{(\ve)}(e_{j_2})\big)\big|_{X_\lambda^{(1)}} \Big\}\Big)^{2m}\nn\\
&\leq 3^{2m-1} \Big\{ \Big( \sum_{j=1}^{\ell-k}
\Log \big( d\Phi_0^{(\ve)}(e_j) \big)\big|_{X_{i}^{(2)}}\Big)^{2m}\nn\\
&\hspace{0.8cm}+L\max_{1 \leq \lambda<\nu \leq d_1}
\Big( \sum_{1 \leq j_1 < j_2 \leq \ell-k}
\Log \big( d\Phi_0^{(\ve)}(e_{j_1})\big)\big|_{X_\lambda^{(1)}} 
\Log \big( d\Phi_0^{(\ve)}(e_{j_2})\big)\big|_{X_\nu^{(1)}} \Big)^{2m} \nn\\
&\hspace{0.8cm}+L \max_{1 \leq \lambda<\nu \leq d_1}
\Big( \sum_{1 \leq j_1 < j_2 \leq \ell-k}\Log \big( d\Phi_0^{(\ve)}(e_{j_1})\big)\big|_{X_\nu^{(1)}} 
\Log \big( d\Phi_0^{(\ve)}(e_{j_2})\big)\big|_{X_\lambda^{(1)}} \Big)^{2m}\Big\},
\end{align}
where we put 
$$
L:=\frac{1}{2}\max_{i=1, 2, \dots, d_2}\max_{1 \leq \lambda < \nu \leq d_1}
\Big|[\![X_\lambda^{(1)}, X_\nu^{(1)}]\!]\big|_{X_i^{(2)}}\Big|.
$$
We fix $i=1, 2, \dots, d_2$. 
Then we have
\begin{align}\label{step2-3-3}
\Big( \sum_{j=1}^{\ell-k}\Log \big( d\Phi_0^{(\ve)}(e_j) \big)\big|_{X_i^{(2)}}\Big)^{2m}
&= (\ell - k)^{2m} \Big( \sum_{j=1}^{\ell-k} \frac{1}{\ell -k}  
\Log \big( d\Phi_0^{(\ve)}(e_j) \big)\big|_{X_i^{(2)}} \Big)^{2m} \nn\\
&\leq (\ell - k)^{2m} \sum_{j=1}^{\ell-k}\frac{1}{\ell-k} 
\Log \big( d\Phi_0^{(\ve)}(e_j) \big)\big|_{X_i^{(2)}}^{2m} \nn\\
&\leq  \|d\Phi_0^{(\ve)}\|_\infty^{4m}(\ell-k)^{2m} \leq M^{4m}(\ell-k)^{2m} .
\end{align}
by applying the Jensen inequality.
For $1 \leq \lambda < \nu  \leq d_1, \, x \in \mathcal{F}, \, 0 \leq \ve \leq 1, \, N \in \mathbb{N}$ 
and $c=(e_1, e_2, \dots, e_N) \in \Omega_{x, N}(X)$, we set
$$
\begin{aligned}
\widetilde{\mathcal{N}}_N^{(\lambda, \nu, x)}(\Phi_0^{(\ve)}; c)&:=
 \sum_{1 \leq j_1 < j_2 \leq N} \mathcal{J}_\lambda^{(\ve)}(j_1)\mathcal{J}_\nu^{(\ve)}(j_2)
=\sum_{j_2=2}^{N} \mathcal{J}_\nu^{(\ve)}(j_2)\sum_{j_1=1}^{j_2-1}\mathcal{J}_\lambda^{(\ve)}(j_1).
\end{aligned}
$$
Then we also see that 
$\{ \widetilde{\mathcal{N}}_N^{(\lambda, \nu, x)}\}_{N=1}^\infty$ 
is an $\mathbb{R}$-valued martingale 
for every $1 \leq\lambda <\nu \leq d$ and $x \in \mathcal{F}$. 
By applying the Burkholder--Davis--Gundy inequality 
with the exponent $2m$, we have
\begin{align}\label{step2-3-4}
&\sum_{c \in \Omega_{x, N}(X)}p_{\ve}(c)
\big(\widetilde{\mathcal{N}}_N^{(\lambda, \nu, x)}(c)\big)^{2m}\nn\\
&\leq \mathcal{C}_{(2m)}^{2m} \sum_{c \in \Omega_{x, N}(X)} 
p_{\ve}(c) \Big\{ \sum_{j_2=2}^N 
 \mathcal{J}_\nu^{(\ve)}(j_2)^2 \times 
 \Big( \sum_{j_1=1}^{j_2-1}\mathcal{J}_\lambda^{(\ve)}(j_1)\Big)^2\Big\}^m\nn\\
& \leq \mathcal{C}_{(2m)}^{2m} \sum_{c \in \Omega_{x, N}(X)} p_{\ve}(c)(N-1)^m 
 \sum_{j_2=2}^N \frac{1}{N-1}\mathcal{J}_\nu^{(\ve)}(j_2)^{2m}
 \Big( \sum_{j_1=1}^{j_2-1}
\mathcal{J}_\lambda^{(\ve)}(j_1)\Big)^{2m}\nn\\
&\leq  \mathcal{C}_{(2m)}^{2m} N^{m} 
\sum_{j_2=2}^{N} \frac{1}{N-1} \Big( \sum_{c \in \Omega_{x, N}(X)} p_{\ve}(c)
\mathcal{J}_\nu^{(\ve)}(j_2)^{4m}\Big)^{1/2}\nn\\
&\hspace{1cm}\times\Big\{  \sum_{c \in \Omega_{x, N}(X)} p_{\ve}(c) 
\Big( \sum_{j_1=1}^{j_2-1}\mathcal{J}_\lambda^{(\ve)}(j_1)\Big)^{4m}\Big\}^{1/2}\nn\\
&\leq  \mathcal{C}_{(2m)}^{2m}M^{2m} N^{m} 
\sum_{j_2=2}^{N} \frac{1}{N-1} \Big\{  \sum_{c \in \Omega_{x, N}(X)} p_{\ve}(c) 
\Big( \sum_{j_1=1}^{j_2-1}\mathcal{J}_\lambda^{(\ve)}(j_1)\Big)^{4m}\Big\}^{1/2},
\end{align}
where we used Jensen's inequality for the third line and Schwarz' inequality for the final line. 
Then the again use of the Burkholder--Davis--Gundy inequality with the exponent $4m$ gives
\begin{align}\label{step2-3-5}
 &\sum_{c \in \Omega_{x, N}(X)} p_{\ve}(c) 
 \Big( \sum_{j_1=1}^{j_2-1}\mathcal{J}_\lambda^{(\ve)}(j_1)\Big)^{4m}\nn\\
 &\leq \mathcal{C}_{(4m)}^{4m}  \sum_{c \in \Omega_{x, N}(X)} p_{\ve}(c)
 \Big(\sum_{j_1=1}^{j_2-1} 
\mathcal{J}_\lambda^{(\ve)}(j_1)^2\Big)^{2m}\nn\\
 &= \mathcal{C}_{(4m)}^{4m} (j_2-1)^{2m} \sum_{c \in \Omega_{x, N}(X)} p_{\ve}(c)
 \Big(\sum_{j_1=1}^{j_2-1}\frac{1}{j_2-1} \mathcal{J}_\lambda^{(\ve)}(j_1)^{2}\Big)^{2m}\nn\\
 &\leq \mathcal{C}_{(4m)}^{4m} j_2^{2m}  \sum_{c \in \Omega_{x, N}(X)} p_{\ve}(c)
 \sum_{j_1=1}^{j_2-1}\frac{1}{j_2-1}\mathcal{J}_\lambda^{(\ve)}(j_1)^{4m}
  \leq \mathcal{C}_{(4m)}^{4m}M^{4m}j_2^{2m}.
\end{align}
It follows from (\ref{step2-3-4}) and (\ref{step2-3-5}) that 
\begin{align}\label{step2-3-6}
\sum_{c \in \Omega_{x, N}(X)}p_{\ve}(c)
\big(\widetilde{\mathcal{N}}_N^{(\lambda, \nu, x)}(c)\big)^{2m} 
&\leq \mathcal{C}_{(2m)}^{2m}M^{2m}N^m
\sum_{j_2=2}^N \frac{1}{N-1}( \mathcal{C}_{(4m)}^{4m}M^{4m}j_2^{2m})^{1/2}\nn\\
&\leq \mathcal{C}_{(2m)}^{2m}\mathcal{C}_{(4m)}^{2m}M^{4m}N^{2m}.
\end{align}
Hence, (\ref{step2-3-6}) implies
\begin{align}\label{step2-3-6'}
&\sum_{c \in \Omega_{x, N}(X)}p_{\ve}(c)
\Big( \sum_{1 \leq j_1 < j_2 \leq N}
\Log \big( d\Phi_0^{(\ve)}(e_{j_1})\big)\big|_{X_\lambda^{(1)}} 
\Log \big( d\Phi_0^{(\ve)}(e_{j_2})\big)\big|_{X_\nu^{(1)}} \Big)^{2m}\nn\\
&\leq 4^{2m-1}\sum_{c \in \Omega_{x, N}(X)}p_{\ve}(c)\Big\{ 
\big(\widetilde{\mathcal{N}}_N^{(\lambda, \nu, x)}(c)\big)^{2m}
+\Big( \ve^2\rho_{\mathbb{R}}(\gamma_p)|_{X_\lambda^{(1)}}
\rho_{\mathbb{R}}(\gamma_p)|_{X_\nu^{(1)}}
 \cdot \frac{N(N-1)}{2} \Big)^{2m}\nn\\
&\hspace{0.3cm}+\Big( \ve\rho_{\mathbb{R}}(\gamma_p)|_{X_\nu^{(1)}}
\sum_{1 \leq j_1 < j_2 \leq N}
\mathcal{J}_{\lambda}^{(\ve)}(j_1)\Big)^{2m}
+\Big(\ve\rho_{\mathbb{R}}(\gamma_p)|_{X_\lambda^{(1)}}
\sum_{1 \leq j_1 < j_2 \leq N}
\mathcal{J}_{\nu}^{(\ve)}(j_2)\Big)^{2m}\Big\}\nn\\
&\leq 4^{2m-1}\Big\{ \mathcal{C}_{(2m)}^{2m}\mathcal{C}_{(4m)}^{2m}M^{4m}N^{2m} 
+ 2^{-2m}M^{4m}\ve^{4m}N^{4m}\nn\\
&\hspace{1cm}+2M^{2m}\ve^{2m}N^{2m}\max_{1 \leq i \leq d_1}
\sum_{c \in \Omega_{x, N}(X)}p_{\ve}(c)
\Big(\sum_{j=1}^N\mathcal{J}_{i}^{(n)}(j)\Big)^{2m}\Big\}\nn\\
&\leq 4^{2m-1}\Big\{ \mathcal{C}_{(2m)}^{2m}\mathcal{C}_{(4m)}^{2m}M^{4m}N^{2m} 
+ 2^{-2m}M^{4m}\ve^{4m}N^{4m}\nn\\
&\hspace{1cm}+2M^{2m}\ve^{2m}N^{2m}\Big(2^{2m}\mathcal{C}_{(2m)}^{2m} 
M^{m}N^{m}
+2^{2m-1}M^{2m}\ve^{2m}N^{2m}\Big)\Big\},
\end{align}
where we used (\ref{martingale2-1}) for the final line. 

We now put $\ve=n^{-1/2}$ and $N=\ell-k$. 
Then we have, for $1 \leq \lambda < \nu \leq d_1$, 
\begin{align}\label{step2-3-7}
&\sum_{c \in \Omega_{x, \ell-k}(X)}p_{\ve}(c)
\Big( \sum_{1 \leq j_1 < j_2 \leq \ell-k}
\Log \big( d\Phi_0^{(\ve)}(e_{j_1})\big)\big|_{X_\lambda^{(1)}} 
\Log \big( d\Phi_0^{(\ve)}(e_{j_2})\big)\big|_{X_\nu^{(1)}} \Big)^{2m}\nn\\
&\leq 4^{2m-1}M^{4m}\Big( \mathcal{C}_{(2m)}^{2m}\mathcal{C}_{(2m)}^{4m}
+2^{-2m}+2^{2m+1}\mathcal{C}_{(2m)}^{2m}M^{-m}+2^{2m}\Big)(\ell-k)^{2m}
\end{align}
due to (\ref{step2-3-6'})
and $(\ell-k)/n<1$.  
We obtain
$$
\begin{aligned}
\mathbb{E}^{\mathbb{P}_{x_*}^{(n^{-1/2})}}
\Big[ \big\| \Log\big((\widetilde{\mathcal{Y}}_{t_k}^{(n^{-1/2}, n)})^{-1} \cdot
 \widetilde{\mathcal{Y}}_{t_{\ell}}^{(n^{-1/2}, n)}\big)
 \big|_{\g^{(2)}}\big\|_{\g^{(2)}}^{2m}\Big] 
\leq C^{(2)}\Big( \frac{\ell-k}{n}\Big)^{2m}.
\end{aligned}
$$
by combining (\ref{step2-3-1}) with (\ref{step2-3-2}), (\ref{step2-3-3}) 
and (\ref{step2-3-7}),
where 
$$
C^{(2)}:=d_2^{2m}3^{2m-1}\Big\{ M^{4m}+2L \cdot 
4^{2m-1}M^{4m}\Big( \mathcal{C}_{(2m)}^{2m}\mathcal{C}_{(2m)}^{4m}
+2^{-2m}+2^{2m+1}\mathcal{C}_{(2m)}^{2m}M^{-m}+2^{2m}\Big)\Big\}.
$$
This means (\ref{expectation2-2}) and  
we thus obtain (\ref{tight2-1}).

\vspace{2mm}
\noindent
{\bf Step\,4.} We show  (\ref{tight2-2}) at the last step.
Suppose that $t_k \leq s \leq t_{k+1}$ and $t_{\ell} \leq t \leq t_{\ell+1}$ for some
$1 \leq k \leq \ell \leq n$. 
Then we have
$$
\begin{aligned}
d_{\mathrm{CC}}(\mathcal{Y}_s^{(n^{-1/2}, n; \,2)}, 
\mathcal{Y}_{t_{k+1}}^{(n^{-1/2}, n; \,2)}) 
&= (k-ns)d_{\mathrm{CC}}(\mathcal{Y}_{t_k}^{(n^{-1/2}, n; \,2)}, 
\mathcal{Y}_{t_{k+1}}^{(n^{-1/2}, n; \,2)}),\\
d_{\mathrm{CC}}(\mathcal{Y}_{t_\ell}^{(n^{-1/2}, n; \,2)}, 
\mathcal{Y}_{t}^{(n^{-1/2}, n; \,2)}) 
&= (nt-\ell)d_{\mathrm{CC}}(\mathcal{Y}_{t_\ell}^{(n^{-1/2}, n; \,2)}, 
\mathcal{Y}_{t_{\ell+1}}^{(n^{-1/2}, n; \,2)})
\end{aligned}
$$
by noting that the piecewise smooth stochastic process $\mathcal{Y}_\cdot^{(n^{-1/2}, n)}$ is given
by the $d_{\mathrm{CC}}$-geodesic interpolation. 
Hence, (\ref{tight2-1}) and the triangle inequality yield 
$$
\begin{aligned}
&\mathbb{E}^{\mathbb{P}_{x_*}^{(n^{-1/2})}}
\Big[ d_{\mathrm{CC}}(\mathcal{Y}_{s}^{(n^{-1/2}, n; \,2)}, 
\mathcal{Y}_{t}^{(n^{-1/2}, n; \,2)})^{4m}\Big]\\
&\leq 3^{4m-1} \Big\{ (k+1-ns)^{4m} 
\cdot C\Big(\frac{1}{n}\Big)^{2m} + C \Big(\frac{\ell-k-1}{n}\Big)^{2m}
+  (nt-\ell)^{4m}\cdot C\Big(\frac{1}{n}\Big)^{2m} \Big\}\\
&\leq C\Big\{ (t_{k+1}-s)^{2m} + (t_\ell-t_{k+1})^{2m} 
+ (t-t_{\ell})^{2m}\Big\} \leq C(t-s)^{2m}.
\end{aligned}
$$  
This completes the proof of Lemma \ref{sublemma2}. \qed

\vspace{2mm}
In what follows, we write
$$
\begin{aligned}
d\mathcal{Y}_{s, t}^{(\ve, n)*}&:=
(\mathcal{Y}_s^{(\ve, n)})^{-1} * \mathcal{Y}_t^{(\ve, n)} 
\qquad (0 \leq \ve \leq 1, \, n \in \mathbb{N}, \, 0 \leq s \leq t \leq 1)
 \end{aligned}
$$
for brevity. 
We now show the following by using Lemma \ref{sublemma2}.

\begin{lm}\label{sublemma3}
For $m, n \in \mathbb{N}$, $k=1, 2, \dots, r$ and  $\alpha<\frac{2m-1}{4m}$, 
there exist a measurable set $\Omega_k^{(n)} \subset \Omega_{x_*}(X)$ and  
a non-negative random variable 
$\mathcal{K}_k^{(n)} \in L^{4m}(\Omega_{x_*}(X)
 \to \mathbb{R}; \,\mathbb{P}_{x_*}^{(n^{-1/2})})$ 
 such that 
$\mathbb{P}_{x_*}^{(n^{-1/2})}(\Omega_{k}^{(n)})=1$ and 
\begin{equation}\label{Holder-estimate}
d_{\mathrm{CC}}\big( \mathcal{Y}_s^{(n^{-1/2}, n; \, k)}(c),  
\mathcal{Y}_t^{(n^{-1/2}, n; \, k)}(c)\big)
\leq \mathcal{K}_k^{(n)}(c)(t-s)^\alpha 
\quad (c \in \Omega_{k}^{(n)}, \, 0 \leq s \leq t \leq 1). 
\end{equation}
\end{lm}

\noindent
{\bf Proof.} We partially follow Lyons' original proof 
(cf.~\cite[Theorem 2.2.1]{Lyons}) of the extension theorem  
in rough path theory. 
We prove (\ref{Holder-estimate}) by induction on the step number 
$k=1, 2, \dots, r$. 

\vspace{2mm}
\noindent
{\bf Step~1.} 
In the cases $k=1, 2$, we have already obtained (\ref{Holder-estimate})
in Lemma \ref{sublemma2}. 
In fact, 
 (\ref{Holder-estimate}) for $k=1, 2$ 
are  obtained by a simple application of
the Kolmogorov--Chentsov criterion 
with the bound
\begin{equation}\label{bound}
\|\mathcal{K}_k^{(n)}\|_{L^{4m}(\mathbb{P}_{x_*}^{(n^{-1/2})})} \leq \frac{5C}{(1-2^{-\theta})(1-2^{\alpha-\theta})} 
\qquad (n, m \in \mathbb{N}, \, k=1, 2),
\end{equation}
where $\theta=(2m-1)/4m$ and $C$ is a constant 
independent of $n$ which appears in the right-hand side of (\ref{tight2-2}). 
See e.g., Stroock \cite[Theorem 4.3.2]{Stroock} for details. 

\vspace{2mm}
\noindent
{\bf Step~2.} We now fix $n \in \mathbb{N}$.
 Assume that (\ref{Holder-estimate}) holds up to step $k$. 
We note that this assumption is equivalent to the existence of 
measurable sets $\{\widehat{\Omega}_j^{(n)}\}_{j=1}^k$ and 
non-negative random variables $\{\widehat{\mathcal{K}}_j^{(n)}\}_{j=1}^k$ such that
$\mathbb{P}_{x_*}^{(n^{-1/2})}(\widehat{\Omega}_{j}^{(n)})=1$ and 
\begin{align}\label{Holder-estimate2}
&\big\| \big( d\mathcal{Y}_{s, t}^{(n^{-1/2}, n)*}(c) \big)^{(j)}\big\|_{\mathbb{R}^{d_j}}
\leq \widehat{\mathcal{K}}_{j}^{(n)}(c)(t-s)^{j\alpha} 
\qquad (c \in \widehat{\Omega}_j^{(n)}, \, 0 \leq s \leq t \leq 1)
\end{align}
with $\widehat{\mathcal{K}}_j^{(n)} \in L^{4m/j}(\Omega_{x_*}(X) \to \mathbb{R}; 
\,  \mathbb{P}_{x_*}^{(n^{-1/2})})$ for $n, m \in \mathbb{N}$ and $j=1, 2, \dots, k$.

We fix $0 \leq s \leq t \leq 1, \, n \in \mathbb{N}$ and 
write $\widehat{\Omega}_{k+1}^{(n)}=\bigcap_{j=1}^k \widehat{\Omega}_j^{(n)}$. 
We denote by $\Delta$ the partition $\{s=t_0<t_1<\cdots<t_N=t\}$  
of the time interval $[s, t]$ independent of $n \in \mathbb{N}$. 
We now define two $G_{(0)}^{(k+1)}$-valued random variables $\mathcal{Z}_{s, t}^{(n)}$
and $\mathcal{Z}(\Delta)_{s, t}^{(n)}$ by
$$
\begin{aligned}
\big(\mathcal{Z}_{s, t}^{(n)}\big)^{(j)}&:=\begin{cases}
\big( d\mathcal{Y}_{s, t}^{(n^{-1/2}, n)*}\big)^{(j)}, & (j=1, 2, \dots, k), \\
\bm{0} & (j=k+1),
\end{cases}\\
\mathcal{Z}(\Delta)_{s, t}^{(n)}&
:=\mathcal{Z}_{t_0, t_1}^{(n)} * \mathcal{Z}_{t_1, t_2}^{(n)}
* \cdots * \mathcal{Z}_{t_{N-1}, t_N}^{(n)},
\end{aligned}
$$
respectively. For $i=1, 2, \dots, d_{k+1}$, 
(\ref{CBH-formula}) and (\ref{Holder-estimate2}) imply
$$
\begin{aligned}
&\Big| \big(\mathcal{Z}(\Delta)_{s, t}^{(n)}(c)\big)^{(k+1)}_{i*}
- \big(\mathcal{Z}(\Delta \setminus \{t_{N-1}\})_{s, t}^{(n)}(c)\big)^{(k+1)}_{i*}\Big|\\
&=\Big| \big( \mathcal{Z}_{t_{N-2}, t_{N-1}}^{(n)}(c) * 
\mathcal{Z}_{t_{N-1}, t_{N}}^{(n)}(c)\big)^{(k+1)}_{i*}
 - \big(\mathcal{Z}_{t_{N-2}, t_{N}}^{(n)}(c)\big)^{(k+1)}_{i*}\Big|\\
 &=\Bigg| \sum_{\substack{|K_1|+|K_2| = k+1 \\ |K_1|, |K_2|>0}}C_{K_1, K_2}
 \mathcal{P}^{K_1}_*\big(\mathcal{Z}_{t_{N-2}, t_{N-1}}^{(n)}(c)\big)
 \mathcal{P}^{K_2}_*\big(\mathcal{Z}_{t_{N-1}, t_{N}}^{(n)}(c)\big)\Bigg|\\
 &\leq C \sum_{\substack{|K_1|+|K_2| = k+1 \\ |K_1|, |K_2|>0}}
 \Big|\mathcal{P}^{K_1}_* \big(d\mathcal{Y}_{t_{N-2}, t_{N-1}}^{(n^{-1/2}, n)*}(c)\big)\Big| 
 \Big| \mathcal{P}^{K_2}_*\big(d\mathcal{Y}_{t_{N-1}, t_N}^{(n^{-1/2}, n)*}(c)\big)\Big|   \\
 &\leq \widehat{\mathcal{K}}_{k+1}^{(n)}(c)(t_N-t_{N-2})^{(k+1)\alpha} \leq \widehat{\mathcal{K}}_{k+1}^{(n)}(c)
 \Big( \frac{2}{N-1}(t-s)\Big)^{(k+1)\alpha} \qquad (c \in \widehat{\Omega}_{k+1}^{(n)}),
\end{aligned}
$$
where the 
random variable $\widehat{\mathcal{K}}_{k+1}^{(n)} : \Omega_{x_*}(X) \LA \mathbb{R}$
is given by
$$
\begin{aligned}
\widehat{\mathcal{K}}_{k+1}^{(n)}(c)&:=C \sum_{\substack{|K_1|+|K_2| = k+1 \\ |K_1|, |K_2|>0}}
\mathcal{Q}^{(n, K_1)}(c)\mathcal{Q}^{(n, K_2)}(c),\\
\mathcal{Q}^{(n, K)}(c)&:=\widehat{\mathcal{K}}_{k_1}^{(n)}(c) \cdot \cdots \cdot 
\widehat{\mathcal{K}}_{k_\ell}^{(n)}(c) \qquad \big( K=\big((i_1, k_1), (i_2, k_2), \dots, (i_\ell, k_\ell)\big) \big).
\end{aligned}
$$
We emphasize that $\widehat{\mathcal{K}}_{k+1}^{(n)}$ is non-negative and has the following integrability:
$$
\begin{aligned}
\mathbb{E}^{\mathbb{P}_{x_*}^{(n^{-1/2})}}\big[(\widehat{\mathcal{K}}_{k+1}^{(n)})^{4m/(k+1)}\big]
&\leq C\sum_{\substack{k_1, \dots, k_\ell>0 \\ k_1+\dots +k_\ell=k+1}}
\mathbb{E}^{\mathbb{P}_{x_*}^{(n^{-1/2})}}\Big[
\big(\widehat{\mathcal{K}}_{k_1}^{(n)} \cdot \cdots \cdot 
\widehat{\mathcal{K}}_{k_\ell}^{(n)}\big)^{4m/(k+1)}\Big]\\
&\leq C\sum_{\substack{k_1, \dots, k_\ell>0 \\ k_1+\dots +k_\ell=k+1}}
\prod_{\lambda=1}^\ell \mathbb{E}^{\mathbb{P}_{x_*}^{(n^{-1/2})}}
\Big[ \big( \widehat{\mathcal{K}}_{k_\lambda}^{(n)}\big)^{4m/k_\lambda}\Big]^{k_\lambda/(k+1)}<\infty,
\end{aligned}
$$
where we used the generalized H\"older inequality for the second line. 
We then have
\begin{align}\label{Holder2}
&\Big| \big(\mathcal{Z}(\Delta)_{s, t}^{(n)}(c)\big)_{i_*}^{(k+1)}\Big|\nn\\
&\leq\Big| \big(\mathcal{Z}(\Delta \setminus \{t_{N-1}\})_{s, t}^{(n)}(c)\big)_{i_*}^{(k+1)}\Big|
+\widehat{\mathcal{K}}_{k+1}^{(n)}(c)
 \Big( \frac{2}{N-1}(t-s)\Big)^{(k+1)\alpha} \nn\\
&\leq \Big| \big(\mathcal{Z}(\{s, t\})_{s, t}^{(n)}(c)\big)_{i_*}^{(k+1)}\Big|
+\sum_{\ell=1}^{N-2}\widehat{\mathcal{K}}_{k+1}^{(n)}(c)
\Big( \frac{2}{N-\ell}\Big)^{(k+1)\alpha}(t-s)^{(k+1)\alpha}\nn\\
&\leq  \Big| \big(\mathcal{Z}_{s, t}^{(n)}(c)\big)_{i_*}^{(k+1)}\Big|
+\widehat{\mathcal{K}}_{k+1}^{(n)}(c)  2^{(k+1)\alpha} \zeta\big( (k+1)\alpha\big) 
 (t-s)^{(k+1)\alpha} \nn\\
&\leq \widehat{\mathcal{K}}_{k+1}^{(n)}(c)(t-s)^{(k+1)\alpha} 
\qquad (i=1, 2, \dots, d_{k+1}, \, c \in \widehat{\Omega}_{k+1}^{(n)})
\end{align}
by successively removing points until the partition $\Delta$ coincides with $\{s, t\}$,
where $\zeta(z)$ denotes the Riemann zeta function 
$\zeta(z):=\sum_{n=1}^\infty (1/n^z)$ for $z \in \mathbb{R}.$

We now show that the family $\{\mathcal{Z}(\Delta)_{s, t}^{(n)}\}$ is a 
Cauchy sequence. 
Let $\delta>0$ and we take two partitions $\Delta=\{s=t_0<t_1< \cdots <t_N=t\}$ 
and $\Delta'$ of $[s, t]$
independent of $n \in \mathbb{N}$ satisfying $|\Delta|, |\Delta'|<\delta$. 
We set $\widehat{\Delta}:=\Delta \cup \Delta'$ and write
$$
\widehat{\Delta}_\ell=\widehat{\Delta} \cap [t_\ell, t_{\ell+1}]
=\{t_\ell=s_{\ell0}<s_{\ell1}<\cdots<s_{\ell L_\ell}=t_{\ell+1}\} \qquad (\ell=0, 1, \dots, N-1).
$$
Then (\ref{CBH-formula}) and (\ref{Holder2}) give
$$
\begin{aligned}
&\Big| \big(\mathcal{Z}(\Delta)_{s, t}^{(n)}(c)\big)_{i_*}^{(k+1)}
- \big(\mathcal{Z}(\widehat{\Delta})_{s, t}^{(n)}(c)\big)_{i_*}^{(k+1)}\Big|\\
&=\Big| \big(\mathcal{Z}_{t_0, t_1}^{(n)}(c) 
* \cdots * \mathcal{Z}_{t_{N-1}, t_N}^{(n)}(c)\big)^{(k+1)}_{i_*}
 - \big( \mathcal{Z}(\widehat{\Delta}_0)_{t_0, t_1}^{(n)}(c)
* \cdots * \mathcal{Z}(\widehat{\Delta}_{N-1})_{t_{N-1}, t_N}^{(n)}(c)\big)^{(k+1)}_{i_*}\Big|\\
&=\Big| \big(\mathcal{Z}_{t_0, t_1}^{(n)}(c)\big)^{(k+1)}_{i_*}
+\big(\mathcal{Z}_{t_1, t_2}^{(n)}(c) 
* \cdots * \mathcal{Z}_{t_{N-1}, t_N}^{(n)}(c)\big)^{(k+1)}_{i_*}\\
&\hspace{1cm} -\big(\mathcal{Z}(\widehat{\Delta}_0)_{t_0, t_1}^{(n)}(c)\big)^{(k+1)}_{i_*}
-\big(\mathcal{Z}(\widehat{\Delta}_1)_{t_1, t_2}^{(n)} (c)
* \cdots * \mathcal{Z}(\widehat{\Delta}_{N-1})_{t_{N-1}, t_N}^{(n)}(c)\big)^{(k+1)}_{i_*}\Big|\\
&\leq \widehat{\mathcal{K}}_{k+1}^{(n)}(c)(t_1-t_0)^{(k+1)\alpha}
+\Big| \big(\mathcal{Z}_{t_1, t_2}^{(n)}(c) 
* \cdots * \mathcal{Z}_{t_{N-1}, t_N}^{(n)}(c)\big)^{(k+1)}_{i_*}\\
&\hspace{1cm} - \big( \mathcal{Z}(\widehat{\Delta}_0)_{t_1, t_2}^{(n)}(c)
* \cdots * \mathcal{Z}(\widehat{\Delta}_{N-1})_{t_{N-1}, t_N}^{(n)}(c)\big)^{(k+1)}_{i_*}\Big|
\qquad (i=1, 2, \dots, d_{k+1}, \, c \in \widehat{\Omega}_{k+1}^{(n)}).
\end{aligned}
$$
By repeating this kind of estimate and noting $(k+1)\alpha>1$, we obtain
\begin{align}\label{Holder3}
&\Big| \big(\mathcal{Z}(\Delta)_{s, t}^{(n)}(c)\big)_{i_*}^{(k+1)}
- \big(\mathcal{Z}(\widehat{\Delta})_{s, t}^{(n)}(c)\big)_{i_*}^{(k+1)}\Big|\nn\\
&\leq \sum_{\ell=0}^{N-1}\widehat{\mathcal{K}}_{k+1}^{(n)}(c)(t_{\ell+1}-t_{\ell})^{(k+1)\alpha}\nn\\
&\leq \widehat{\mathcal{K}}_{k+1}^{(n)}(c) \Big( \max_{\Delta}(t_{\ell+1}-t_{\ell})^{(k+1)\alpha-1}\Big)
\sum_{\ell=0}^{N-1}(t_{\ell+1}-t_{\ell})\nn\\
&\leq \widehat{\mathcal{K}}_{k+1}^{(n)}(c)(t-s) \cdot \delta^{(k+1)\alpha-1} 
\qquad (i=1, 2, \dots, d_{k+1}, \, c \in \widehat{\Omega}_{k+1}^{(n)}).
\end{align}
Thus, (\ref{Holder3})  leads to
$$
\begin{aligned}
&\Big| \big(\mathcal{Z}(\Delta)_{s, t}^{(n)}(c)\big)_{i_*}^{(k+1)}
- \big(\mathcal{Z}(\Delta')_{s, t}^{(n)}(c)\big)_{i_*}^{(k+1)}\Big|\\
&\leq \Big| \big(\mathcal{Z}(\Delta)_{s, t}^{(n)}(c)\big)_{i_*}^{(k+1)}
- \big(\mathcal{Z}(\widehat{\Delta})_{s, t}^{(n)}(c)\big)_{i_*}^{(k+1)}\Big|
+\Big| \big(\mathcal{Z}(\widehat{\Delta})_{s, t}^{(n)}(c)\big)_{i_*}^{(k+1)}
- \big(\mathcal{Z}(\widetilde{\Delta})_{s, t}^{(n)}(c)\big)_{i_*}^{(k+1)}\Big|\\
&\leq 2\widehat{\mathcal{K}}_{k+1}^{(n)}(c)(t-s) \cdot \delta^{(k+1)\alpha-1} \LA 0 
\qquad (i=1, 2, \dots, d_{k+1}, \, c \in \widehat{\Omega}_{k+1}^{(n)})
\end{aligned}
$$
as $\delta \searrow 0$ uniformly in $0 \leq s \leq t \leq 1$. 
Therefore, there exists, for $0 \leq s \leq t \leq 1$,  
$$
\ol{\mathcal{Z}}_{s, t}^{(n)}(c):=\begin{cases}
\dis \lim_{|\Delta| \searrow 0}\mathcal{Z}(\Delta)_{s, t}^{(n)}(c) &
 (c \in \widehat{\Omega}_{k+1}^{(n)}), \\
\bm{1}_G &  (c \in \Omega_{x_*}(X) \setminus \widehat{\Omega}_{k+1}^{(n)}).
\end{cases}
$$
satisfying
$$
\big\| \big( \ol{\mathcal{Z}}_{s, t}^{(n)}(c) \big)^{(k+1)}\big\|_{\mathbb{R}^{d_{k+1}}}
\leq \widehat{\mathcal{K}}_{k+1}^{(n)}(c)(t-s)^{(k+1)\alpha} \qquad (c \in \widehat{\Omega}_{k+1}^{(n)}),
$$
due to (\ref{Holder2}). 
We will show 
$$
\ol{\mathcal{Z}}_{s, t}^{(n)}(c)
=\mathcal{Y}_{s}^{(n^{-1/2}, n; \, k+1)}(c)^{-1} *
\mathcal{Y}_{t}^{(n^{-1/2}, n; \, k+1)}(c)
\qquad (0 \leq s \leq t \leq 1, \, c \in \widehat{\Omega}_{k+1}^{(n)})
$$
as the last step. 
For this,  it is sufficient to check that
\begin{align}\label{last aim}
\big( \ol{\mathcal{Z}}_{s, t}^{(n)}(c) \big)^{(k+1)}
&=\big(d\mathcal{Y}_{s, t}^{(n^{-1/2}, n)*}(c)\big)^{(k+1)}
\qquad (0 \leq s \leq t \leq 1, \, c  \in \widehat{\Omega}_{k+1}^{(n)})
\end{align}
by the definition of $\ol{\mathcal{Z}}_{s, t}^{(n)}$.
We fix $i=1, 2, \dots, d_{k+1}$ and $c \in \widehat{\Omega}_{k+1}^{(n)}$. Put 
$$
\Psi^i_{s, t}(c):=\big(d\mathcal{Y}_{s, t}^{(n^{-1/2}, n)*}(c)\big)^{(k+1)}_{i*} - 
\big( \ol{\mathcal{Z}}_{s, t}^{(n)}(c) \big)^{(k+1)}_{i_*} \quad (0 \leq s \leq t \leq 1).
$$
Then we easily see that $\Psi^i_{s, t}(c)$ is additive in the sense that 
\begin{equation}\label{additivity}
\Psi^i_{s, t}(c)=\Psi^i_{s, u}(c)+\Psi^i_{u, t}(c) \qquad (0 \leq s \leq u \leq t \leq 1).
\end{equation}
Since the piecewise smooth stochastic process 
$(\mathcal{Y}_t^{(n^{-1/2}, n)})_{0 \leq t \leq 1}$ is given by the 
$d_{\mathrm{CC}}$-geodesic interpolation of $\{\mathcal{X}^{(n^{-1/2}, n)}_{t_k}\}_{k=0}^n$,  
we have 
$$
\big\| \big(d\mathcal{Y}_{s, t}^{(n^{-1/2}, n)*}(c)\big)^{(k+1)}\big\|_{\mathbb{R}^{d_{k+1}}}
\leq \widetilde{\mathcal{K}}^{(n)}_{k+1}(c)(t-s)^{(k+1)\alpha} 
\qquad (c \in \widetilde{\Omega}_{k+1}^{(n)})
$$
for some set $\widetilde{\Omega}_{k+1}^{(n)}$ 
with $\mathbb{P}_{x_*}^{(n^{-1/2})}(\widetilde{\Omega}_{k+1}^{(n)})=1$ and 
random variable $\widetilde{\mathcal{K}}^{(n)}_{k+1} : \Omega_{x_*}(X) \LA \mathbb{R}$. 
Thus, we have
$$
\big| \Psi^i_{s, t}(c)\big| \leq \big( \widetilde{\mathcal{K}}^{(n)}_{k+1}(c)
+\widehat{\mathcal{K}}^{(n)}_{k+1}(c)\big)(t-s)^{(k+1)\alpha}
 \qquad (0 \leq s \leq t \leq 1, \, 
 c \in \widetilde{\Omega}_{k+1}^{(n)} \cap \widehat{\Omega}_{k+1}^{(n)}).
$$
We may write $\widehat{\Omega}_{k+1}^{(n)}$ instead of
$\widetilde{\Omega}_{k+1}^{(n)} \cap \widehat{\Omega}_{k+1}^{(n)}$ by abuse of notation.
Because its probability equals one.
For any small $\ve>0$, there is a sufficiently large $N \in \mathbb{N}$ such that 
$1/N < \ve$. 
We then obtain as $\ve \searrow 0$, 
$$
\begin{aligned}
\Big|\Psi^i_{0, t}(c)\Big|
&= \Big| \Psi^i_{0, 1/N}(c) + \Psi^i_{1/N, 2/N}(c) + 
\cdots + \Psi^i_{[Nt]/N, t}(c)\Big|\\
&\leq \big( \widetilde{\mathcal{K}}^{(n)}_{k+1}(c)
+\widehat{\mathcal{K}}^{(n)}_{k+1}(c)\big)\ve^{(k+1)\alpha-1}
\Big\{ \underbrace{\frac{1}{N}+\cdots 
+ \frac{1}{N}}_{[Nt]\text{-times}}+\Big(t-\frac{[Nt]}{N}\Big)\Big\}\\
&= \big( \widetilde{\mathcal{K}}^{(n)}_{k+1}(c)
+\widehat{\mathcal{K}}^{(n)}_{k+1}(c)\big)\ve^{(k+1)\alpha-1}t \LA 0 
\qquad (0 \leq t \leq 1, \, c \in \widehat{\Omega}_{k+1}^{(n)}) 
\end{aligned}
$$
by  (\ref{additivity}) and $(k+1)\alpha-1>0$.
This implies that $\Psi^i_{0, t}(c)=0$ for 
$0 \leq t \leq 1$ and $c \in \widehat{\Omega}_{k+1}^{(n)}$. 
Hence, it follows from (\ref{additivity}) that 
$
\Psi^i_{s, t}(c)
=\Psi^i_{0, t}(c)-\Psi^i_{0, s}(c)=0$ for
$0 \leq s \leq t \leq 1$ and $c \in \widehat{\Omega}_{k+1}^{(n)})$,
which leads to (\ref{last aim}). 
Consequently, we know that there are a measurable 
set $\Omega_{k+1}^{(n)} \subset \Omega_{x_*}(X)$
with probability one and a non-negative random variable 
$\mathcal{K}_{k+1}^{(n)} \in L^{4m}(\Omega_{x_*}(X) \to \mathbb{R}; \,\mathbb{P}_{x_*}^{(n^{-1/2})})$
satisfying
$$
\begin{aligned}
d_{\mathrm{CC}}\big( \mathcal{Y}_s^{(n^{-1/2}, n; \, k+1)}(c),  
\mathcal{Y}_t^{(n^{-1/2}, n; \, k+1)}(c)\big)
&\leq \mathcal{K}_{k+1}^{(n)}(c)(t-s)^\alpha \quad (c \in \Omega_{k+1}^{(n)}, \, 0 \leq s \leq t \leq 1).
\end{aligned}
$$
This completes the proof of Lemma \ref{sublemma3}. \qed

\vspace{2mm}
\noindent
{\bf Proof of Lemma \ref{tightness2}}. 
For $m, n \in \mathbb{N}$ and $\widehat{\alpha}<\frac{2m-1}{4m}$,  
Lemma \ref{sublemma3} implies
$$
\begin{aligned}
\mathbb{E}^{\mathbb{P}_{x_*}^{(n^{-1/2})}}\Big[d_{\mathrm{CC}}\big( \mathcal{Y}_s^{(n^{-1/2}, n; \, r)},  
\mathcal{Y}_t^{(n^{-1/2}, n; \, r)}\big)^{4m}\Big]
&\leq \mathbb{E}^{\mathbb{P}_{x_*}^{(n^{-1/2})}}
\big[\big(\mathcal{K}_r^{(n)}\big)^{4m}\big](t-s)^{4m\widehat{\alpha}}. \\
\end{aligned}
$$
for $0 \leq s \leq t \leq 1$.
We thus have, by (\ref{bound}),
$$
\mathbb{E}^{\mathbb{P}_{x_*}^{(n^{-1/2})}}\Big[d_{\mathrm{CC}}\big( \mathcal{Y}_s^{(n^{-1/2}, n; \, r)},  
\mathcal{Y}_t^{(n^{-1/2}, n; \, r)}\big)^{4m}\Big]
\leq C(t-s)^{4m\widehat{\alpha}} \qquad (0 \leq s \leq t \leq 1). 
$$ 
for some constant $C>0$ independent of $n \in \mathbb{N}$. 
Furthermore, thanks to {\bf (A2)} and $\Phi_0^{(0)}(x_*)=\bm{1}_G$, 
there is a sufficiently large constant $C>0$ such that 
$$
\sup_{n \in \mathbb{N}} \big\| \log\big(\Phi_0^{(n^{-1/2})}(x_*)\big)\big|_{\g^{(k)}}\big\|_{\g^{(k)}} \leq C
\qquad (k=1, 2, \dots, r). 
$$
By applying
the Kolmogorov tightness criterion, it follows that 
the family $\{\Prob^{(n^{-1/2}, n)}\}_{n=1}^\infty$ is tight in 
$C^{0, \alpha\text{\rm{-H\"ol}}}([0, 1]; G_{(0)})$ for $\alpha <\frac{4m\widehat{\alpha}-1}{4m}<\frac{1}{2}-\frac{1}{2m}$. 
By letting $m \to \infty$, we complete the proof. \qed

\vspace{2mm}

By using Theorem \ref{CLT2}, Lemma \ref{sublemma3} 
and repeating the same argument as in \cite[Lemma 4.8]{IKN}, 
we easily obtain the convergence
of finite dimensional distribution of $(\mathcal{Y}^{(n^{-1/2}, n)})_{0 \leq t \leq 1}$, 
that is, for every $0 \leq s_1 < s_2 < \cdots < s_\ell \leq 1 \, (\ell \in \mathbb{N})$,
$$
(\mathcal{Y}_{s_1}^{(n^{-1/2}, n)}, 
\mathcal{Y}_{s_2}^{(n^{-1/2}, n)}, 
\dots, \mathcal{Y}_{s_\ell}^{(n^{-1/2}, n)})
\overset{(d)}{\LA} (Y_{s_1}, Y_{s_2}, \dots, Y_{s_\ell}) \text{ as }n \to \infty.
$$

Finally, by combining this with Lemma \ref{tightness2}, 
we complete the proof of Theorem \ref{FCLT2}.

\subsection{Proof of Theorem \ref{FCLT2-general}}

In this subsection, we give a proof of Theorem \ref{FCLT2-general}. 
We show that 
the same pathwise H\"older estimate as Lemma \ref{sublemma3}
also holds for the stochastic process $\{\ol{\mathcal{Y}}^{(n^{-1/2}, n)}\}_{n=1}^\infty$
by applying the corrector method.

\begin{lm}\label{sublemma-general}
For $m, n \in \mathbb{N}$ and  $\alpha<\frac{2m-1}{4m}$,
there exist an $\mathcal{F}_\infty$-measurable set $\ol{\Omega}_r^{(n)} \subset \Omega_{x_*}(X)$ and 
a non-negative random variable 
$\ol{\mathcal{K}}_r^{(n)} \in L^{4m}\big(\Omega_{x_*}(X) \to \mathbb{R}; \mathbb{P}_{x_*}\big)$
such that 
$\mathbb{P}_{x_*}^{(n^{-1/2})}(\ol{\Omega}_r^{(n)})=1$ and 
\begin{equation}\label{pathwise-general}
d_{\mathrm{CC}}\big(\ol{\mathcal{Y}}_{s}^{(n^{-1/2}, n)}(c), \ol{\mathcal{Y}}_{t}^{(n^{-1/2}, n)}(c)\big)
\leq \ol{\mathcal{K}}_r^{(n)}(c)(t-s)^\alpha 
\qquad (c \in \ol{\Omega}_{r}^{(n)}, \, 0 \leq s <t \leq 1). 
\end{equation}
\end{lm}

\noindent
{\bf Proof.} Fix $n \in \mathbb{N}$ and $1 \leq k \leq \ell \leq n$. 
By triangular inequality, we have
$$
\begin{aligned}
&d_{\mathrm{CC}}(\ol{\mathcal{Y}}_{t_k}^{(n^{-1/2}, n)}, \ol{\mathcal{Y}}_{t_{\ell}}^{(n^{-1/2}, n)})\\
&\leq d_{\mathrm{CC}}(\ol{\mathcal{Y}}_{t_k}^{(n^{-1/2}, n)}, \mathcal{Y}_{t_k}^{(n^{-1/2}, n)})
+d_{\mathrm{CC}}(\mathcal{Y}_{t_k}^{(n^{-1/2}, n)}, \mathcal{Y}_{t_\ell}^{(n^{-1/2}, n)})
+d_{\mathrm{CC}}(\mathcal{Y}_{t_\ell}^{(n^{-1/2}, n)}, \ol{\mathcal{Y}}_{t_\ell}^{(n^{-1/2}, n)}). 
\end{aligned}
$$
We set 
$\mathcal{Z}^{(n)}_t:=(\mathcal{Y}^{(n^{-1/2}, n)}_t)^{-1} * \ol{\mathcal{Y}}^{(n^{-1/2}, n)}_t$ 
for $0 \leq t \leq 1$ and $n \in \mathbb{N}$. 
By definition, we have
$$
\log\big(\mathcal{Z}_{t_k}^{(n)})|_{\g^{(1)}}=\frac{1}{\sqrt{n}}\mathrm{Cor}_{\g^{(1)}}^{(n^{-1/2})}(w_{k}) 
\qquad (n \in \mathbb{N}, \, k=0, 1, \dots, n)
$$
so that there is a constant $C>0$ such that 
$\big\| \log\big(\mathcal{Z}_{t_k}^{(n)})|_{\g^{(1)}} \big\|_{\g^{(1)}} \leq Cn^{-1/2}$
for $n \in \mathbb{N}$ and $k=0, 1, 2, \dots, n$. 
Moreover, it follows from the choice of the components of $\Phi_0^{(\ve)}(x) \, (0 \leq \ve \leq 1, \, x \in V)$ that 
$\big\| \log\big(\mathcal{Z}_{t_k}^{(n)})|_{\g^{(i)}} \big\|_{\g^{(i)}} \leq Cn^{-i/2}$
for $n \in \mathbb{N}$ and $k=0, 1, 2, \dots, n$. 
By noting the equivalence of homogeneous norms, we have
\begin{equation}\label{corrector-est-1}
d_{\mathrm{CC}}(\ol{\mathcal{Y}}_{t_k}^{(n^{-1/2}, n)}, \mathcal{Y}_{t_k}^{(n^{-1/2}, n)})
\leq C\big\|\mathcal{Z}_{t_k}^{(n)}\big\|_{\Hom}
=C\sum_{i=1}^r \big\| \log\big(\mathcal{Z}_{t_k}^{(n)})|_{\g^{(i)}} \big\|_{\g^{(i)}}^{1/i}
\leq \frac{C}{\sqrt{n}}
\end{equation}
for $n \in \mathbb{N}$ and $k=0, 1, 2, \dots, n$. 
Then it follows from Lemma \ref{sublemma3} and (\ref{corrector-est-1}) that there exist
an $\mathcal{F}_\infty$-measurable set $\ol{\Omega}_r^{(n)} \subset \Omega_{x_*}(X)$ and  
a non-negative random variable 
$\ol{\mathcal{K}}_r^{(n)} \in L^{4m}\big(\Omega_{x_*}(X) \to \mathbb{R}; \mathbb{P}_{x_*}\big)$
such that 
$\mathbb{P}_{x_*}^{(n^{-1/2})}(\ol{\Omega}_r^{(n)})=1$ and 
\begin{align}
d_{\mathrm{CC}}\big(\ol{\mathcal{Y}}_{t_k}^{(n^{-1/2}, n)}(c), \ol{\mathcal{Y}}_{t_\ell}^{(n^{-1/2}, n)}(c)\big)
&\leq \frac{C}{\sqrt{n}} + \mathcal{K}_r^{(n)}(c)\Big(\frac{\ell-k}{n}\Big)^\alpha + \frac{C}{\sqrt{n}}\nn\\
&\leq \ol{\mathcal{K}}_r^{(n)}(c)\Big(\frac{\ell-k}{n}\Big)^\alpha \label{general-est-partial}
\end{align}
for $c \in \ol{\Omega}_r^{(n)}$ and $0 \leq k \leq \ell \leq n$. 
For $0 \leq s < t \leq 1$, take $0 \leq k \leq \ell \leq n$ such that
$k/n \leq s < (k+1)n$ and $\ell/n \leq t < (\ell+1)/n$. 
Since the stochastic process $(\ol{\mathcal{Y}}_t^{(n^{-1/2}, \,n)})_{0 \leq t \leq 1}$
is also give by the $d_{\mathrm{CC}}$-geodesic interpolation, we have
$$
\begin{aligned}
d_{\mathrm{CC}}\big(\ol{\mathcal{Y}}_{s}^{(n^{-1/2}, n)}, \ol{\mathcal{Y}}_{t_{k+1}}^{(n^{-1/2}, n)}\big)
&=(k-ns)d_{\mathrm{CC}}\big(\ol{\mathcal{Y}}_{t_k}^{(n^{-1/2}, n)}, \ol{\mathcal{Y}}_{t_{k+1}}^{(n^{-1/2}, n)}\big),\\
d_{\mathrm{CC}}\big(\ol{\mathcal{Y}}_{t_\ell}^{(n^{-1/2}, n)}, \ol{\mathcal{Y}}_{t}^{(n^{-1/2}, n)}\big)
&=(nt-\ell)d_{\mathrm{CC}}\big(\ol{\mathcal{Y}}_{t_\ell}^{(n^{-1/2}, n)}, \ol{\mathcal{Y}}_{t_{\ell+1}}^{(n^{-1/2}, n)}\big).
\end{aligned}
$$
We then use the triangular inequality and (\ref{general-est-partial}) to obtain
$$
\begin{aligned}
&d_{\mathrm{CC}}\big(\ol{\mathcal{Y}}_{s}^{(n^{-1/2}, n)}(c), \ol{\mathcal{Y}}_{t}^{(n^{-1/2}, n)}(c)\big)\\
&\leq (k-ns)  \ol{\mathcal{K}}_r^{(n)}(c)\Big(\frac{1}{n}\Big)^\alpha
+\ol{\mathcal{K}}_r^{(n)}(c)\Big(\frac{\ell-k-1}{n}\Big)^\alpha + (nt-\ell)\ol{\mathcal{K}}_r^{(n)}(c)
\Big(\frac{1}{n}\Big)^\alpha\\
&\leq \ol{\mathcal{K}}_r^{(n)}(c) \Big\{ \Big(\frac{k+1}{n}-s\Big)^\alpha + \Big(\frac{\ell-k-1}{n}\Big)^\alpha
+\Big(t - \frac{\ell}{n}\Big)^\alpha\Big\} 
\leq \ol{\mathcal{K}}_r^{(n)}(c) (t-s)^\alpha \qquad (c \in \ol{\Omega}_r^{(n)}),
\end{aligned}
$$
which completes the proof. \qed

\vspace{2mm}
\noindent
{\bf Proof of Theorem \ref{FCLT2-general}.}
We split the proof into two steps.

\vspace{2mm}
\noindent
{\bf Step~1.} We show that 
the sequence $\{\ol{\mathcal{Y}}^{(n^{-1/2}, n)}\}_{n=1}^\infty$ converges in law
to $(Y_t)_{0 \leq t \leq 1}$ in 
$C([0, 1]; G_{(0)})$ as $n \to \infty$. 
For $0 \leq t \leq 1$, take an integer $0 \leq k \leq n$ such that $k/n \leq t <(k+1)/n$. 
Then, by the triangular inequality, 
(\ref{Holder-estimate}), (\ref{pathwise-general}) and (\ref{corrector-est-1}),
we have, $\mathbb{P}_{x_*}^{(n^{-1/2})}$-almost surely, 
\begin{align}
&d_{\mathrm{CC}}(\mathcal{Y}_t^{(n^{-1/2}, n)}, \ol{\mathcal{Y}}_t^{(n^{-1/2}, n)})\nn\\
&\leq d_{\mathrm{CC}}(\mathcal{Y}_{t_k}^{(n^{-1/2}, n)}, \mathcal{Y}_t^{(n^{-1/2}, n)})
+d_{\mathrm{CC}}(\mathcal{Y}_{t_k}^{(n^{-1/2}, n)}, \ol{\mathcal{Y}}_{t_k}^{(n^{-1/2}, n)})
+d_{\mathrm{CC}}(\ol{\mathcal{Y}}_{t_k}^{(n^{-1/2}, n)}, \ol{\mathcal{Y}}_t^{(n^{-1/2}, n)})\nn\\
&\leq \mathcal{K}_r^{(n)}\Big(t-\frac{k}{n}\Big)^\alpha + 
\frac{C}{\sqrt{n}} + \ol{\mathcal{K}}_r^{(n)}\Big(t-\frac{k}{n}\Big)^\alpha\nn\\
&\leq \big\{ \mathcal{K}_r^{(n)}+\ol{\mathcal{K}}_r^{(n)}+C\big\}\Big(\frac{1}{\sqrt{n}}\Big)^\alpha
\qquad \Big(m \in \mathbb{N}, \, \alpha<\frac{2m-1}{4m}\Big). \label{difference-y}
\end{align}
Let $\rho$ be a metric on $C([0, 1]; G_{(0)})$ defined by
$$
\rho(w^{(1)}, w^{(2)}):=
\max_{0 \leq t \leq 1}d_{\mathrm{CC}}(w_t^{(1)}, w_t^{(2)}) + d_{\mathrm{CC}}(\bm{1}_G, (w_0^{(1)})^{-1} * w_0^{(2)}).
$$
By applying the Chebyshev inequality, (\ref{final bound}) and (\ref{difference-y}), we have, for 
$\ve >0$ and $m \in \mathbb{N}$, 
$$
\begin{aligned}
&\mathbb{P}_{x_*}^{(n^{-1/2})}\Big(\rho(\mathcal{Y}^{(n^{-1/2}, n)}, \ol{\mathcal{Y}}^{(n^{-1/2}, n)})>\ve\Big)\\
&\leq \Big(\frac{1}{\ve}\Big)^{4m}
\mathbb{E}^{\mathbb{P}_{x_*}^{(n^{-1/2})}}
\Big[\rho(\mathcal{Y}^{(n^{-1/2}, n)}, \ol{\mathcal{Y}}^{(n^{-1/2}, n)})^{4m}\Big]\\
&\leq \Big(\frac{2}{\ve}\Big)^{4m}\Big\{
\mathbb{E}^{\mathbb{P}_{x_*}^{(n^{-1/2})}}\Big[\max_{0 \leq t \leq 1}
d_{\mathrm{CC}}(\mathcal{Y}_t^{(n^{-1/2}, n)}, \ol{\mathcal{Y}}_t^{(n^{-1/2}, n)})^{4m}\Big]
+C \cdot \Big(\frac{1}{\sqrt{n}}\Big)^{4m}\Big\}\\
&\leq  3^{4m-1}\Big(\frac{2}{\ve}\Big)^{4m}\Big(\frac{1}{\sqrt{n}}\Big)^{4m\alpha}
\Big\{ \mathbb{E}^{\mathbb{P}_{x_*}^{(n^{-1/2})}}\big[(\mathcal{K}_r^{(n)})^{4m}\big]
+\mathbb{E}^{\mathbb{P}_{x_*}^{(n^{-1/2})}}\big[(\ol{\mathcal{K}}_r^{(n)})^{4m}\big]
+\mathbb{E}^{\mathbb{P}_{x_*}^{(n^{-1/2})}}\big[C^{4m}\big]\Big\}\\
&\hspace{1cm}+C \cdot \Big(\frac{2}{\ve}\Big)^{4m} \cdot \Big(\frac{1}{\sqrt{n}}\Big)^{4m} \to 0
\quad (n \to \infty).
\end{aligned}
$$
Therefore, by Slutzky's theorem, we obtain
the convergence in law of $\{\ol{\mathcal{Y}}^{(n^{-1/2}, n)}\}_{n=1}^\infty$ to the 
diffusion process $(Y_t)_{0 \leq t \leq 1}$ in $C([0, 1]; G_{(0)})$ as $n \to \infty$. 

\vspace{2mm}
\noindent
{\bf Step~2.}
The previous step  implies the convergence of 
the finite-dimensional distribution of $\{\ol{\mathcal{Y}}^{(n^{-1/2}, n)}\}_{n=1}^{\infty}$. 
On the other hand, we can show
 that the sequence of probability measures 
$\{\ol{\Prob}^{(n^{-1/2}, n)}:=\mathbb{P}_{x_*}^{(n^{-1/2})} \circ (\ol{\mathcal{Y}}^{(n^{-1/2}, n)})^{-1}\}_{n=1}^\infty$
is tight
in $C^{0, \alpha{{\normalfont \hol}}}([0, 1]; G_{(0)})$, 
by noting Lemma \ref{sublemma-general} and
 by following the same argument as the proof of Lemma \ref{tightness2}. 
Therefore, we complete the proof by combining these two. \qed 

\section{Example}

In this section, we discuss an example of $\Gamma$-nilpotent covering graph
in the case where $\Gamma$ is the 3-dimensional discrete
Heisenberg group $\mathbb{H}^3(\mathbb{Z})$ defined by
$$
\mathbb{H}^3(\mathbb{Z}):=\Bigg\{ 
\begin{bmatrix}
1 & x & z \\
0 & 1 & y \\
0 & 0 & 1
\end{bmatrix} \, \Bigg| \, x, y, z \in \mathbb{Z}\Bigg\}.
$$
It is known that $\Gamma=\mathbb{H}^3(\mathbb{Z})$ is 
a cocompact lattice in the 3-dimensional Heisenberg group $G=\mathbb{H}^3(\mathbb{R})=(\mathbb{R}^3, \star)$,
where the product $\star$ on $\mathbb{R}^3$ is given by
$$
(x, y, z)\star (x', y', z')=(x+x', y+y', z+z'+xy').
$$
The corresponding Lie algebra $\g$ is given by $(\mathbb{R}^3, [\cdot, \cdot])$
generated by $X_1=(1, 0, 0), X_2=(0, 1, 0)$ and $X_3=(0, 0, 1)$ with
$[X_1, X_2]=X_3$ and $[X_1, X_3]=[X_2, X_3]
=\bm{0}_{\g}$ under the matrix bracket $[X, Y]:=XY-YX$ for $X, Y \in \g$. 
We then see that the Lie algebra $\g$ is decomposed as
$\g=\g^{(1)} \oplus \g^{(2)}$, where
$\g^{(1)}:=\Span\{X_1, X_2\}$ and $\g^{(2)}:=\Span\{X_3\}$.

Suppose that $\Gamma=\mathbb{H}^3(\mathbb{Z})$ is generated by two elements
$\gamma_1=(1, 0, 0)$ and $\gamma_2=(0, 1, 0)$. We put 
$\bg:=(1/3, 1/3, 1/3) \in G=\mathbb{H}^3(\mathbb{R})$ and define
$$
\begin{aligned}
V_{1}&:= \big\{ g=\gamma_{i_1}^{\ve_1} \star \dots \star \gamma_{i_\ell}^{\ve_\ell} 
                                \star \bm{1}_G \, \big| \, i_k \in \{1, 2\}, \, \ve_k=\pm1 \, 
                                (1 \leq k \leq \ell), \, \ell \in \mathbb{N} \cup \{0\}\big\},\\
V_{2}&:= \big\{ g=\gamma_{i_1}^{\ve_1} \star\dots \star \gamma_{i_\ell}^{\ve_\ell} 
                                \star \bg \, \big| \, i_k \in \{1, 2\}, \, \ve_k=\pm1 \, 
                                (1 \leq k \leq \ell), \, \ell \in \mathbb{N} \cup \{0\}\big\}.
\end{aligned}
$$
Consider an $\mathbb{H}^3(\mathbb{Z})$-nilpotent covering graph $X=(V, E)$ defined by
$V=V_1 \sqcup V_2$ and 
$$
\begin{aligned}
E &:= \big\{ (g, h) \in V_1 \times V_2\, | \, g^{-1} \star h=\bg_1, \,
                                  \gamma_1^{-1} \star \bg_1, \, \gamma_2^{-1}  \star \bg_1 \big\}.
\end{aligned}
$$
We note that  $X$ is invariant under the actions $\gamma_1$ and $\gamma_2$. 
We call the graph $X$ a 3-{\it dimensional Heisenberg hexagonal lattice}, 
which is a generalization of the classical hexagonal lattice to the nilpotent setting (see Figure 1).  
We have discussed two $\mathbb{H}^3(\mathbb{Z})$-nilpotent covering graphs called
3-{\it dimensional Heisenberg triangular lattice} and 
3-{\it dimensional Heisenberg dice lattice}, respectively.
It should be noted that a 3-dimensional Heisenberg dice lattice is regarded as 
a hybrid of 3-dimensional Heisenberg triangular lattice
and a 3-dimensional Heisenberg hexagonal lattice. 
See \cite[Section 6]{IKN} for details of them. 
The quotient graph $X_0=(V_0, E_0)=\Gamma \backslash X$ is given by $V_0=\{\x, \y\}$ 
and $E_0=\{e_i, \ol{e}_i \, | \, 1 \leq i \leq 3\}$ (see Figure 2).

\begin{figure}[htbp]
\begin{center}
\includegraphics[width=15cm]{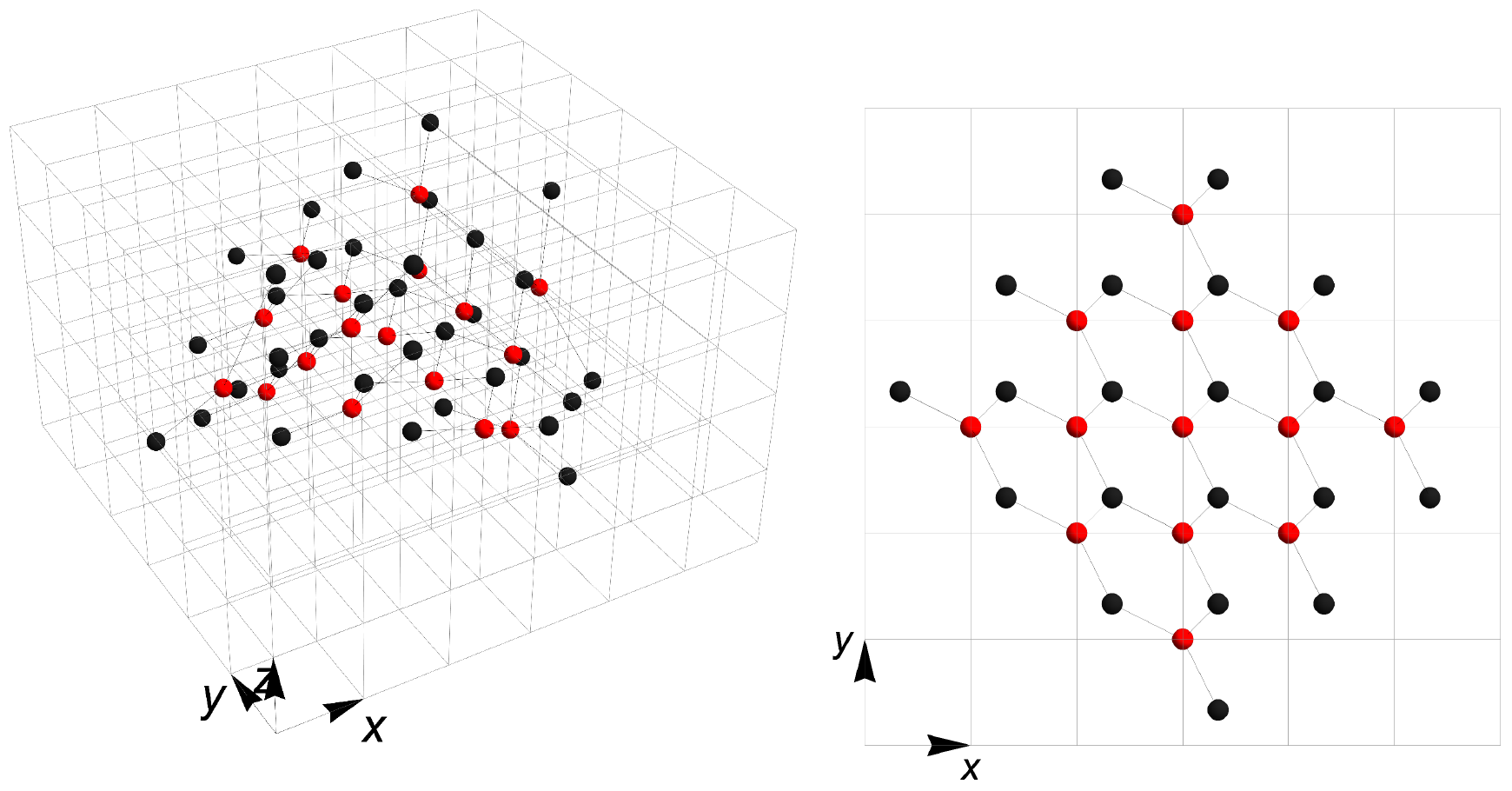}
\caption{A part of 3-dimensional Heisenberg hexagonal lattice and 
the projection of it on the $xy$-plane}
\label{dices}
\end{center}
\vspace{-0.5cm}
\end{figure}

\begin{figure}[htbp]
\begin{center}
\includegraphics[width=6cm]{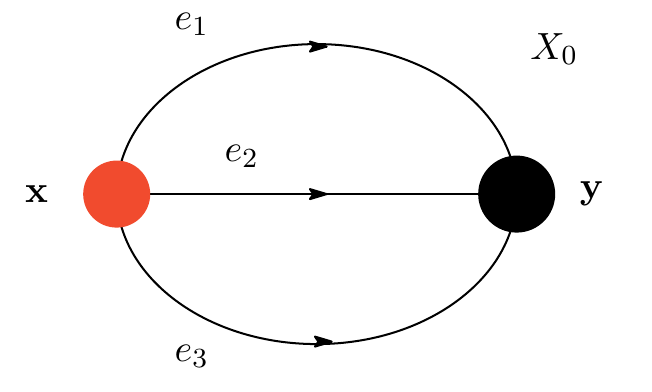}
\caption{The quotient graph $X_0=\mathbb{H}^3(\mathbb{Z}) \backslash X$}
\label{dices}
\end{center}
\vspace{-0.5cm}
\end{figure}

We now discuss a non-symmetric random walk on $X$
with a transition probability $p : E \LA (0, 1]$ given by
\begin{align}
p\big((g, g \star \bg)\big)&=\alpha, & 
p\big((g, g \star \gamma_1^{-1}  \star \bg)\big)&=\beta, & 
p\big((g, g \star \gamma_2^{-1}  \star \bg)\big)&=\gamma, \nn\\
p\big(\ol{(g, g \star \bg)}\big)&=\alpha', &
p\big(\ol{(g, g \star \gamma_1^{-1}  \star \bg)}\big)&=\beta', &
p\big(\ol{(g, g \star \gamma_2^{-1}  \star \bg)}\big)&=\gamma', \nn
\end{align}
for every $g \in V_1$, where $\alpha, \beta, \gamma, \alpha', \beta', \gamma'>0$,
$\alpha+\beta+\gamma=1$ and $\alpha'+\beta'+\gamma'=1$.
Then, the transition probability induces a non-symmetric random walk
on $X_0$ with
$$
p(e_1)=\alpha, \quad p(e_2)=\beta, \quad p(e_3)=\gamma, \quad
p(\ol{e}_1)=\alpha', \quad p(\ol{e}_2)=\beta', \quad p(\ol{e}_3)=\gamma'.
$$
The invariant measure $m : V_0=\{\x, \y\} \LA (0, 1]$ is given by
$m(\x)=1/2$ and $m(\y)=1/2$.
The first homology group $\h_1(X_0, \mathbb{R})$ is spanned 
by two 1-cycles
$[c_1]:=[e_1 * \ol{e}_2]$ and $[c_2]:=[e_3*\ol{e}_2]$.
Then the homological direction is calculated as
$
\gamma_p=(\check{\alpha}/2)[c_1]+(\check{\gamma}/2)[c_2].
$
Here, we use the notations
$$
\hat{\alpha}=\alpha+\alpha', \quad \hat{\beta}=\beta+\beta', \quad
\hat{\gamma}=\gamma+\gamma', \quad
\check{\alpha}=\alpha-\alpha', \quad
\check{\beta}=\beta-\beta', \quad
\check{\gamma}=\gamma-\gamma' \quad
$$
for simplicity. The canonical surjective linear map
$\rho_{\mathbb{R}} : \h_1(X_0, \mathbb{R}) \LA \g^{(1)}$ is given by
$\rho_{\mathbb{R}}([c_1])=X_1$ and $\rho_{\mathbb{R}}([c_2])=X_2$
so that we obtain
$\rho_{\mathbb{R}}(\gamma_p)=(\check{\alpha}/2)X_1+(\check{\gamma}/2)X_2.$

Let $\{\omega_1, \omega_2\} \subset \Hom(\g^{(1)}, \mathbb{R})$ be the dual basis of 
$\{X_1, X_2\}$. By identifying each $\omega_i$ with ${}^t\rho_{\mathbb{R}}(\omega_i) \in \h^1(X_0, \mathbb{R})$
for $i=1, 2$, we can see that $\{\omega_1, \omega_2\}$ is regarded as 
the dual basis of $\{[c_1], [c_2]\} \subset \h_1(X_0, \mathbb{R})$. 

We now discuss the family of transition probabilities $(p_\ve)_{0 \leq \ve \leq 1}$.
By definition, we have
\begin{align}
p_\ve(e_1) &= \frac{1}{2}(\hat{\alpha} + \ve \check{\alpha}), &
p_\ve(e_2) &= \frac{1}{2}(\hat{\beta} + \ve \check{\beta}), &
p_\ve(e_3) &= \frac{1}{2}(\hat{\gamma} + \ve \check{\gamma}), \nn\\
p_\ve(\ol{e}_1) &= \frac{1}{2}(\hat{\alpha} - \ve \check{\alpha}), &
p_\ve(\ol{e}_2) &= \frac{1}{2}(\hat{\beta} - \ve \check{\beta}), &
p_\ve(\ol{e}_3) &= \frac{1}{2}(\hat{\gamma} - \ve \check{\gamma})\nn
\end{align}
for $0 \leq \ve \leq 1$.
Direct computations give us
\begin{align}
\La \omega_1, \omega_1 \Ra_{p_\ve}
&=\frac{\hat{\alpha}(\hat{\beta}+\hat{\gamma})-\ve^2\check{\alpha}^2}{4}, &
\La \omega_1, \omega_2 \Ra_{p_\ve}
&=\frac{\hat{\alpha}\hat{\gamma}+\ve^2\check{\alpha}\check{\gamma}}{4}, \nn\\
\La \omega_2, \omega_2 \Ra_{p_\ve}
&=\frac{(\hat{\alpha}+\hat{\beta})\hat{\gamma} - \ve^2\check{\gamma}^2}{4}. &\nn
\end{align}
We define the $\Gamma$-Albanese torus $\mathrm{Alb}^{\Gamma}$ by 
$(\g^{(1)}/\mathcal{L}, g_0)$, where 
$\mathcal{L}=\{aX_1+bX_2 \, | \, a, b \in \mathbb{Z}\}$. 
Then we have
$$
\begin{aligned}
v_\ve^{-1} \equiv \mathrm{vol}_{g_0^{(\ve)}}(\mathrm{Alb}^{\Gamma})^{-1}
&=\sqrt{\det \big( \La \omega_i, \omega_j \Ra_{p_\ve}\big)_{i, j=1}^2}\\
&=\frac{1}{4}\sqrt{2\hat{\alpha}\hat{\beta}\hat{\gamma}
+\ve^2\hat{\alpha}\check{\alpha}(\hat{\beta}\check{\gamma}+\check{\beta}\hat{\gamma})
+\ve^2\hat{\beta}\check{\beta}(\hat{\alpha}\check{\gamma}+\check{\alpha}\hat{\gamma})
+\ve^2\hat{\gamma}\check{\gamma}(\hat{\alpha}\check{\beta}+\check{\alpha}\hat{\beta})}
\end{aligned}
$$
and
\begin{align}
\la X_1, X_1 \ra_{g_0^{(\ve)}}&=\frac{v_\ve^2}{4}\big\{ (\hat{\alpha}+\hat{\beta})\hat{\gamma} -\ve^2 \check{\gamma}^2\big\}, &
\la X_1, X_2 \ra_{g_0^{(\ve)}}&=\frac{v_\ve^2}{4}(\hat{\alpha}\hat{\gamma}+\ve^2\check{\alpha}\check{\gamma}), \nn\\
\la X_2, X_2 \ra_{g_0^{(\ve)}}&=\frac{v_\ve^2}{4}\big\{\hat{\alpha}(\hat{\beta}+\hat{\gamma})-\ve^2\check{\alpha}^2\big\}. &&\nn
\end{align}

We here give the family of modified harmonic realizations $(\Phi_0^{(\ve)})_{0 \leq \ve \leq 1}$. 
Introduce four arbitrary families $(x_\ve)_{0 \leq \ve \leq 1}$, $(y_\ve)_{0 \leq \ve \leq 1}$,
$(z_\ve)_{0 \leq \ve \leq 1}$ and $(\kappa_\ve)_{0 \leq \ve \leq 1}$ in $\mathbb{R}$
such that $x_0=y_0=z_0=0$. 
Then it follows from (\ref{p_ve-modified}) that the realization given by
$$
\Phi_0^{(\ve)}(\bm{1}_G)=(x_\ve, y_\ve, z_\ve), \qquad
\Phi_0^{(\ve)}(\bg)=\Big(x_\ve+\frac{\hat{\alpha}}{2}, y_\ve+\frac{\hat{\gamma}}{2}, z_\ve+\kappa_\ve\Big)
$$
is $(p_\ve$-)modified harmonic for $0 \leq \ve \leq 1$. 
We can see that the family $(\Phi_0^{(\ve)})_{0 \leq \ve \leq 1}$ satisfies {\bf (A1)} if and only if
$x_\ve+y_\ve=0$ for $0 \leq \ve \leq 1$. 
In addition, if $(z_\ve)_{0 \leq \ve \leq 1}$ and $(\kappa_\ve)_{0 \leq \ve \leq 1}$ satisfy
$
\sup_{0 \leq \ve \leq 1}|z_\ve+\kappa_\ve| \leq C
$
for a positive constant $C$, then the family $(\Phi_0^{(\ve)})_{0 \leq \ve \leq 1}$ 
satisfies {\bf (A2)} as well as {\bf (A1)}. 

Let $\{V_1, V_2\}$ be the Gram--Schmidt orthogonalization of $\{X_1, X_2\}$
with respect to $g_0^{(0)}$, that is, 
$$
V_1=\frac{\sqrt{\hat{\alpha}(\hat{\beta}+\hat{\gamma})}}{2}X_1 
-\frac{\hat{\alpha}\hat{\beta}}{2\sqrt{\hat{\alpha}(\hat{\beta}+\hat{\gamma})}}X_2,
\qquad 
V_2=\frac{2v_0}{\sqrt{\hat{\alpha}(\hat{\beta}+\hat{\gamma})}}X_2.
$$
Then the infinitesimal generator $\A$ defined by (\ref{generator2}) is given by 
$$
\begin{aligned}
\A&=-\frac{1}{2}(V_1^2+V_2^2) - \Big(\frac{\check{\alpha}}{2}X_1+\frac{\check{\gamma}}{2}X_2\Big)\\
&=-\frac{1}{2}(V_1^2+V_2^2) -\frac{\check{\alpha}}{\sqrt{\hat{\alpha}(\hat{\beta}+\hat{\gamma})}}V_1
-\frac{v_0(\check{\alpha}\hat{\alpha}\hat{\gamma}
+\check{\gamma}\hat{\alpha}\hat{\beta}+\check{\gamma}\hat{\alpha}\hat{\gamma})}
{4\sqrt{\hat{\alpha}(\hat{\beta}+\hat{\gamma})}}V_2. 
\end{aligned}
$$

\vspace{4mm}
\noindent
{\bf Acknowledgement.} 
The authors are grateful to Professor Shoichi Fujimori for making pictures 
of the 3-dimensional Heisenberg hexagonal lattice and kindly allowing them 
to use these pictures in the present paper.
They would like to thank Professor Seiichiro Kusuoka for reading of
our manuscript carefully and for giving valuable comments. 
A part of this work was done during the stay of the third named author 
at Hausdorff Center for Mathematics, Universit\"at Bonn in March 2017 
with the support of research fund of Research Institute for Interdisciplinary Science, 
Okayama University. 
He would like to thank Professor Massimiliano Gubinelli 
for warm hospitality and helpful discussions.





\begin{thebibliography}{99}
\bibitem{A1} G. Alexopoulos: 
{\it Convolution powers on discrete groups of polynomial volume growth},
Canad. Math. Soc. Conf. Proc. {\bf 21} (1997), pp. 31--57.
\bibitem{A2} G. Alexopoulos:
{\it Random walks on discrete groups of polynomial growth},
Ann. Probab. {\bf{30}} (2002), pp. 723--801.
\bibitem{A3} G. Alexopoulos:
{\it Sub-Laplacians with drift on Lie groups of polynomial volume growth},
Mem. Amer. Math. Soc. {\bf 155} (2002), no. 739. 
\bibitem{BLP} A. Bensoussan, J. L. Lions and G. Papanicolaou: 
{\it Asymptotic analysis for periodic structures}, 
Studies in Mathematics and its Applications {\bf 5}, 
North-Holland Publishing Co., Amsterdam-New York, 1978.
\bibitem{Biskup} M. Biskup:
{\it{Recent progress on the random conductance model}},
Probab. Surv. {\bf{8}} (2011), pp. 294--373. 
%
\bibitem{FV}
P.K. Friz and N.B. Victoir: Multidimensional Stochastic Processes as Rough Paths, 
Theory and Applications, Cambridge Studies in Advanced Mathematics {\bf{120}}, 
Cambridge Univ. Press, Cambridge, 2010.
\bibitem{Gromov} M. Gromov: 
{\it Groups of polynomial growth and expanding maps}, 
IHES. Publ. Math. {\bf 53} (1981), pp. 53--73.
\bibitem{Ishiwata}
S. Ishiwata:
{\it A central limit theorem on a covering graph with a transformation group 
of polynomial growth},
J. Math. Soc. Japan {\bf 55} (2003), pp. 837--853.
\bibitem{IKK} S. Ishiwata, H. Kawabi and M. Kotani:
{\it{Long time asymptotics of non-symmetric random walks on crystal lattices}},
J. Funct. Anal. {\bf 272} (2017), pp. 1553--1624.
\bibitem{IKN} S. Ishiwata, H. Kawabi and R. Namba:
{\it{Central limit theorems for non-symmetric 
random walks on nilpotent covering graphs: Part I}},
preprint (2018). Available at
{\tt{arXiv:1806.03804}}.
\bibitem{KL} C. Kipnis and C. Landim:
Scaling Limits of Interacting Particle Systems, 
Grundlehren der mathematischen Wissenschaften {\bf 320}, 
Springer-Verlag, Berlin, 1999. 
\bibitem{Kotani} M. Kotani:
{\it A central limit theorem for magnetic transition operators on a crystal lattice},
J. London Math.\ Soc. {\bf{65}} (2002), pp. 464--482.
\bibitem{Kotani contemp} M. Kotani:
{\it An asymptotic of the large deviation for random walks on a crystal lattice},
Contemp.\ Math.\ {\bf{347}} (2004), pp. 141--152.
\bibitem{KS00-CMP} M. Kotani and T. Sunada:
{\it Albanese maps and off diagonal long time asymptotics for the heat kernel},
Comm.\ Math.\ Phys. {\bf{209}} (2000), pp. 633--670.
\bibitem{KS00-TAMS}
M. Kotani and T. Sunada:
{\it Standard realizations of crystal lattices via harmonic maps}, 
Trans. Amer. Math. Soc. {\bf{353}} (2000), pp. 1--20. 
\bibitem{KS06} M. Kotani and T. Sunada:
{\it Large deviation and the tangent cone at infinity of a crystal lattice},
Math. Z. {\bf{254}} (2006), pp. 837--870.
\bibitem{KSS} M. Kotani, T. Shirai and T. Sunada:
{\it Asymptotic behavior of the transition probability of a random walk on an infinite graph},
J. Funct. Anal. {\bf{159}} (1998), pp. 664--689.
\bibitem{Kozlov} S.M. Kozlov: 
{\it{The averaging method and walks in inhomogeneous environments}},
Russian Math. Surveys {\bf{40}} (1985), pp. 73--145.
\bibitem{Kumagai}T. Kumagai:
Random Walks on Disordered Media and their Scaling Limits, 
\'Ecole d'\'Et\'e de Probabilit\'es de Saint-Flour XL-2010, LNM {\bf{2101}}, Springer, Cham, 2014.
\bibitem{LCL} T. Lyons, M.~Caruana and T.~L\'evy: Differential Equations Driven by Rough Paths,
\'Ecole d'\'Et\'e de Probabilit\'es de Saint-Flour XXXIV-2004, 
LNM \ {\bf{1908}}, Springer, Berlin,  2007.
\bibitem{Lyons} T. Lyons: {\it{Differential equations driven by rough signals}},
Rev. Math. Iberoamericana\ {\bf{14}} (1998), pp. 215--310.
\bibitem{LQ} T. Lyons and Z. Qian: {System Control and Rough Paths},
Oxford Mathematical Monographs. Oxford Univ. Press, Oxford, 2002.
\bibitem{Malcev} A. I. Mal\'cev: {\it On a class of homogeneous spaces}, 
Amer. Math. Soc. Transl. {\bf 39} (1951), pp. 276--307.
\bibitem{Pap} G. Pap: {\it{Central limit theorems on stratified Lie groups}}, 
Probability theory and mathematical statistics (Vilnius, 1993), pp. 613--627,
TEV, Vilnius, 1994.
%
\bibitem{PV}G.C. Papanicolaou and S.R.S. Varadhan:
{\it{Boundary value problem with rapidly oscillating random coefficients}},
In: Random Fields, Vol. I, II (Esztergom, 1979), pp. 835--873, 
Colloq. Math. Soc. J\'anos Bolyai, {\bf{27}}, North-Holland, Amsterdam-New York, 1981. 
\bibitem{Rag} M. S. Raghunathan: 
Discrete Subgroups of Lie Groups, Springer-Verlag Berlin, 1972.
\bibitem{Raugi} A. Raugi: {\it Th\`eor\'eme de la limite centrale sur les groupes nilpotents}, 
Z. Wahrsch. Verw. Gebiete. {\bf 43} (1978), pp. 149--172.
\bibitem{Rob} D. W. Robinson: 
Elliptic Operators and Lie Groups, 
Oxford Mathematical Mono- graphs, Oxford Univ. Press, New York, 1991.  
\bibitem{Stroock} D.W. Stroock:
Probability Theory. An Analytic View, Cambridge Univ. Press, 1993. 
\bibitem{S0} T. Sunada:
{\textit{Discrete geometric analysis}}, in ``Analysis on Graphs and its Applications", pp. 51--83, 
Proc. Sympos. Pure Math. {\textbf{77}}, Amer. Math. Soc., Providence, RI, 2008. 
\bibitem{S} T. Sunada:
Topological Crystallography with a View Towards Discrete  Geometric 
Analysis, Surveys and Tutorials in the Applied Mathematical Sciences {\bf{6}},
Springer Japan, 2013.
\bibitem{S2} T. Sunada:
{\it Topics in mathematical crystallography}, in the preceedings
of the symposium ``Groups, graphs and random walks'', 
London Math. Soc., Lecture Note Series 
{\bf{436}}, 
Cambridge Univ. Press, 2017, pp. 473--513.
\bibitem{Tanaka}
R. Tanaka: 
{\it{Hydrodynamic limit for weakly asymmetric simple exclusion processes in crystal lattices}}, 
Comm. Math. Phys. {\bf{315}} (2012), pp. 603--641.
\bibitem{Trotter}H.F. Trotter:
{\it Approximation of semi-groups of operators}, Pacific J. Math. 
{\bf{8}} (1958), pp. 887--919.
\bibitem{VSC}
N.T. Varopoulos, L. Saloff-Coste and T. Coulhon: 
Analysis and Geometry on Groups,
Cambridge Tracts in Mathematics {\bf{100}}. 
Cambridge Univ. Press, Cambridge, 1992.
\end{thebibliography}
\end{document}